\date{}

\documentclass[12pt]{amsart}
\usepackage{latexsym,amsmath,amsfonts,amscd,amssymb}
\usepackage{graphics}
\theoremstyle{remark}

\newcommand{\cZ}{{\mathcal Z}}

\newcommand{\CC}{{\mathbb C}}

\newcommand{\QQ}{{\mathbb Q}}
\newcommand{\RR}{{\mathbb R}}

\newcommand{\ZZ}{{\mathbb Z}}

\renewcommand{\a}{\alpha}
\renewcommand{\b}{\beta}
\renewcommand{\d}{\delta}
\newcommand{\g}{\gamma}

\newcommand{\eps}{\epsilon}

\begin{document}
\title[STATISTICS ON RIEMANN ZEROS]{STATISTICS ON RIEMANN ZEROS}


\author[R. P\'{e}rez Marco]{Ricardo P\'{e}rez Marco}
\address{CNRS, LAGA UMR 7539, Universit\'e Paris XIII \\ 
99, Avenue J.-B. Cl\'ement, 93430-Villetaneuse, France}

\email{ricardo@math.univ-paris13.fr}

\maketitle

{ \centerline{\sc Abstract}}

\bigskip

\begin{minipage}{14cm}
\noindent

We numerically study the 
statistical properties of differences of
zeros of Riemann zeta function and $L$-functions predicted by the
theory of the e\~ne product. In particular, this provides a 
simple algorithm that computes any non-real Riemann zeros from very 
large ones (''self-replicating property of Riemann zeros''). Also 
the algorithm computes the full 
sequence of non-real zeros of Riemann zeta function from the
sequence of non-real zeros of any Dirichlet $L$-function (''zeros of $L$-functions know 
about Riemann zeros''). We also check that the first error to the convergence to the 
classical GUE statistic near $0$ is a Fresnel distribution.
\end{minipage}


\tableofcontents

\newpage

\section{Introduction.} \label{sec:introduction}

The goal of this article is to check numerically the 
predictions for the statistics of differences of zeros of Riemann
zeta function and other $L$-functions anticipated by the theory of
the e\~ne ring structure (see \cite{PM1} and \cite{PM2}). We carry out a detailled 
statistical analysis of the differences of zeros confirming 
the predicted limit distributions.

\bigskip

The study of the distribution for the differences of zeros is not
new. In 1972 H.L. Montgomery (\cite{Mo1}, \cite{Mo2}) investigated these
statistics in relation with the class number problem. Although he
states at the beginning of his article \cite{Mo2} that

\bigskip

{\it "Our goal is to investigate the distribution of the
differences $\gamma -\gamma'$ between the zeros."}
\bigskip

he only studies this distribution for differences of nearby
(consecutive or semi-locally close) zeros. Montgomery formulated
the "Pair correlation conjecture" according to which, when $T\to +\infty$, 
$$
\frac{1}{N(T)} \left | \left \{ (\g ,\g') ; \frac{2\pi \a}{\log T}<
\g-\g'<\frac{2\pi \b} {\log T} \right \} \right | \sim \int_\a^\b
\left ( 1- \left (\frac{\sin (\pi t)} {\pi t} \right)^2\right ) \
dt \ ,
$$
where $N(T)\sim \frac{1}{2\pi} T \log T$ is the asymptotic number
of zeros with positive imaginary part less than $T$, $0<\a <\b$, 
and the sum runs over the imaginary part of non-trivial Riemann
zeros. Note the normalizing factor $\frac{1}{2\pi} \log T$ that
restricts the statistics to nearby zeros. According to the well
known story, F. Dyson recognized the pair correlation distribution
for eigenvalues of large random hermitian matrices used by
physicists (the Gaussian Unitary Ensemble or GUE). This attracted
much attention since it was believed to give some support to the
spectral approach to the Riemann Hypothesis proposed by Hilbert
and Polya (indeed first by Polya according to \cite{Od}), which asks 
to identify an hermitian operator whose
spectrum is composed by the non-real Riemann zeros (this belief is not 
justified as we will discuss later on).

Then in the 80's A.M. Odlyzko performed intensive computations of
the distribution of normalized differences of consecutive Riemann
zeros confirming numerically the GUE distribution conjectured by
Montgomery (\cite{Od}). Later N. Katz and P. Sarnak pursued the numerical
exploration of zeros for other $L$-functions (\cite{KS}). Their work was
continued more recently by many others (see \cite{Co}).


These authors restrict their study to normalized differences of
consecutive zeros. Our goal is to study the distribution of global
differences at large. These global differences appear on compact sets uniformly
distributed in first approximation. At the second order they
present very significant discrepancies with the uniform
distribution (section \ref{sec:self}). The precise location of
these discrepancies is highly significant: We notice a deficiency
of differences of non-real Riemann zeros exactly and precisely at
the same locations as the Riemann zeros. In other words, the statistics 
of differences of Riemann zeros pinpoints the exact location of these zeros.

\bigskip

In particular, this property indicates that large Riemann zeros do
know the location of all Riemann zeros, since the
statistics remain unchanged by removing a finite number (or a
density $0$ subset) of Riemann zeros. It is indeed checked numerically in
section \ref{sec:large} that the statistics of very large zeros do find the 
location the first Riemann zeros.

\bigskip

We also study the differences of the zeros of an arbitrary
Dirichlet $L$-function (section \ref{sec:L-functions_know}). We observe the same
phenomenon. The differences have a uniform distribution except for
deficits located exactly at the precise location of Riemann
zeros. This fact provides a simple algorithm which computes the sequence
of Riemann zeros from the sequence of zeros of any Dirichlet
$L$-function. This result confirms what has been part of the
folklore intuition among the community of specialists (see \cite{Co}):
Zeros of $L$-functions do know about Riemann zeros.

\bigskip

We perform a statistics in sections  \ref{sec:L-functions_replicate}
and \ref{sec:L-functions_mating} which is apparently new
in the literature. We study the differences of the zeros of one
Dirichlet $L$-function $L_1$ against the zeros of another
Dirichlet $L$-function $L_2$. A similar result is observed: The
discrepancies from the uniform distribution are located at the
zeros of another $L$-function which can be computed explicitly. 
When $L_2$ is the Riemann zeta function (section \ref{sec:L-functions_mating}),
the discrepancies are located at the zeros of $L_1$, i.e.
Riemann's zeta function plays the role of the identity for this
mating operation. In section \ref{sec:Euler_factors}, we mate local Euler factors with Riemann 
zeta function. The discrepancies pinpoint to the arithmetic location of the 
poles of the Euler factor.

\bigskip

We also make a numerical study of the distribution predicted by
Montgomery's conjecture (section \ref{sec:fine}). The theory of the e\~ne
product provides a refinement of Montgomery's conjecture. More precisely, it provides, and we
verify numerically, an asymptotics for the error from
Montgomery's GUE limit distribution. At the first order we observe
the predicted Fresnel (or sine cardinal) distribution. These
refinements of Montgomery's conjecture seem also new.

\bigskip

 The numerical results presented in this article are unexpected without the theory 
of the e\~ne product. Even more surprising is that such elementary statistics have been 
unnoticed so far. We
encourage the skeptic readers who want to "put their hand upon the
wound" to check the numerical results by themselves in their own
personal computer (any modern personal computer can do the job). It is
indeed extremely simple as we indicate below. All statistics
presented in this article were performed with public software and
public data on a low profile regular laptop (of 2005, with 500 Mb of RAM memory). 
They are fully reproducible and we provide all necessary information to reproduce
them.

\medskip

 Together with the numerical results we give the simple code
for the computations. These were performed using the public domain
statistical software "R". One can download and install the program
from the site www.r-project.org where tutorials are also
available. This statistical software runs under Linux and Windows
and the source code is open. Some source files for zeros of
the Riemann zeta function were obtained from A. M. Odlyzko's web
site \cite{Od}. A large source file for the 35 first million
non-trivial Riemann zeros and other files with zeros of
$L$-functions are from M. Rubinstein public site \cite{Ru}. The author is very
grateful to A.M. Odlyzko and M. Rubinstein for making available their data.

\bigskip

We have restricted our statistics to Dirichlet $L$-functions only
because we have only access to large lists of zeros for these
functions. These results
hold in general for more general $L$ and zeta functions. We encourage
the readers with access to zero data for other more general $L$-functions 
to perform similar statistics.

\bigskip

From now on we refer to non-real zeros of Riemann's zeta function
simply as "zeros of Riemann's zeta function" or "Riemann zeros".
The same simplified terminology is used for non-real zeros of
general $L$-functions. We also refer to differences as "deltas". We
use the notation $\rho$ for a non-real zero and the notation $\g$
for $\rho=1/2+i\g $. Then the sequence of $\g$'s is the sequence
of non-trivial zeros of Riemann real analytic function on the
vertical line $\{ \Re s=1/2 \}$. We eventually refer to them also
as Riemann zeros.

\bigskip

\textbf{Acknowledgments.} 

These numerical computations were conducted during 2004 and 2005, and this manuscript was 
written in 2005. It was circulated among some close friends. I am very 
grateful for their support during these years. These results 
were announced in the conference in Honor of the 200th birthday of Galois, 
"Differential Equations and Galois Theory", held at I.H.E.S. in October 2011.

\bigskip

\section{Self-replicating property of Riemann zeros.} \label{sec:self}

 We consider the sequence $\cZ=(\rho_i)_ {i\geq 1}$ of Riemann's zeros with positive imaginary part.
We label them in increasing order, $\Im \rho_i > \Im \rho_j$ for $i>
j$. We write $\rho_j=1/2+i \gamma_j$. The first few zeros were
computed by B. Riemann  with high accuracy as was found in
Riemann's Nachlass ([Ri2]). The tabulation of the first $80$
zeros, i.e. all zeros less than $200$ follows:

\begin{align*}
\g_1 &=14.134725142\ldots \\ 
\g_2 &=21.022039639\ldots \\
\g_3 &=25.010857580\ldots \\ 
\g_4 &=30.424876126\ldots \\ 
\g_5 &=32.935061588\ldots \\
\g_6 &=37.586178159 \ldots \\ 
\g_7 &=40.918719012 \ldots \\ 
\g_8 &=43.327073281 \ldots \\ 
\g_9 &=48.005150881 \ldots \\
\g_{10} &=49.773832478 \ldots \\ 
\g_{11} &=52.970321478 \ldots \\ 
\g_{12} &=56.446247697 \ldots \\ 
\g_{13} &=59.347044003 \ldots \\
\g_{14} &=60.831778525 \ldots \\ 
\g_{15} &=65.112544048 \ldots \\ 
\g_{16} &=67.079810529 \ldots \\ 
\g_{17} &=69.546401711 \ldots \\ 
\g_{18} &=72.067157674 \ldots \\ 
\g_{19} &=75.704690699 \ldots \\ 
\g_{20} &=77.144840069 \ldots \\
\g_{21} &=79.337375020 \ldots \\ 
\g_{22} &=82.910380854 \ldots \\ 
\g_{23} &=84.735492981 \ldots \\ 
\g_{24} &=87.425274613 \ldots \\ 
\g_{25} &=88.809111208 \ldots \\ 
\g_{26} &=92.491899271 \ldots \\ 
\g_{27} &=94.651344041 \ldots \\
\g_{28} &=95.870634228 \ldots \\ 
\g_{29} &=98.831194218 \ldots \\ 
\g_{30} &=101.317851006 \ldots \\
\g_{31} &=103.725538040 \ldots \\ 
\g_{32} &=105.446623052 \ldots \\ 
\g_{33} &=107.168611184 \ldots \\ 
\g_{34} &=111.029535543 \ldots \\ 
\g_{35} &=111.874659177 \ldots \\ 
\g_{36} &=114.320220915 \ldots \\ 
\g_{37} &=116.226680321 \ldots 
\end{align*}

\begin{align*}
\g_{38} &=118.790782866 \ldots \\ 
\g_{39} &=121.370125002 \ldots \\ 
\g_{40} &=122.946829294 \ldots \\ 
\g_{41} &=124.256818554 \ldots \\ 
\g_{42} &=127.516683880 \ldots \\ 
\g_{43} &=129.578704200 \ldots \\ 
\g_{44} &=131.087688531 \ldots \\
\g_{45} &=133.497737203 \ldots \\ 
\g_{46} &=134.756509753 \ldots \\ 
\g_{47} &= 138.116042055 \ldots \\ 
\g_{48} &=139.736208952 \ldots \\ 
\g_{49} &=141.123707404 \ldots \\ 
\g_{50} &=143.111845808 \ldots \\ 
\g_{51} &=146.000982487 \ldots \\
\g_{52} &=147.422765343 \ldots \\ 
\g_{53} &=150.053520421 \ldots \\ 
\g_{54} &=150.925257612 \ldots \\ 
\g_{55} &=153.024693811 \ldots \\ 
\g_{56} &=156.112909294 \ldots \\
\g_{57} &=157.597591818 \ldots \\ 
\g_{58} &=158.849988171 \ldots \\ 
\g_{59} &=161.188964138 \ldots \\ 
\g_{60} &=163.030709687 \ldots \\ 
\g_{61} &=165.537069188 \ldots \\ 
\g_{62} &=167.184439978 \ldots \\ 
\g_{63} &=169.094515416 \ldots \\
\g_{64} &=169.911976479 \ldots \\ 
\g_{65} &=173.411536520 \ldots \\
\g_{66} &=174.754191523 \ldots \\ 
\g_{67} &=176.441434298 \ldots \\
\g_{68} &=178.377407776 \ldots \\ 
\g_{69} &=179.916484020 \ldots \\ 
\g_{70} &=182.207078484 \ldots \\ 
\g_{71} &=184.874467848 \ldots \\ 
\g_{72} &=185.598783678 \ldots \\ 
\g_{73} &=187.228922584 \ldots \\ 
\g_{74} &=189.416158656 \ldots 
\end{align*}

\begin{align*}
\g_{75} &=192.026656361 \ldots \\ 
\g_{76} &=193.079726604 \ldots \\ 
\g_{77} &=195.265396680 \ldots \\ 
\g_{79} &=196.876481841 \ldots \\ 
\g_{80} &=198.015309676 \ldots \\ 
\vdots
\end{align*}

We use the notation, for $1\leq j<i$,
$$
\delta_{i,j}=\g_i-\g_j \ .
$$

We fix $N>0$. Given $T>0$ we consider the subset of all deltas of
zeros:

$$
\Delta_{T,N}=\{ 0<\delta_{i,j} \leq T; 1\leq j<i\leq N  \}\subset
]0,T]
$$

\bigskip

We denote by $T_0=\g_N$ and $N=N(T_0)$. We are interested in the
numerical study of the distribution of the elements of
$\Delta_{T,N}$ when $N\to +\infty$ and $T>0$ is kept fixed.

\bigskip

We represent the histogram of values in $\Delta_{T,N}$ by ticks of
$\eps=10^{-m} $ for some $T>0$ and $N$ large. That is, for each
integer $1\leq k\leq 1+T\eps^{-1}$ we count how many deltas
$\delta$ yield $k=[\eps^{-1} \delta]$ (here the brackets denote
the integer part). We denote by $x_k\geq 0$ this number. The
histograms represent the sequence $(x_k)_{1\leq k\leq
1+T\eps^{-1}}$. We represent the figures obtained from two
statistics. One very fast and the second more intensive and
precise. The first one, named (a), computes the deltas of $N=10^5$
zeros with precision $\eps=10^{-1}$ and $T=100$. The computation
takes about 15 minutes on the author's laptop. The statistics
named (b) computes the deltas for $N=5.10^6$ with precision
$\eps=10^{-2}$ and $T=200$. This computation took about 2 days on
the author's laptop. Obviously, reducing $T$ or increasing $\eps$
greatly diminishes the computation time.

\bigskip

Figure 1.a represents the histogram for statistics (a). Figure 1.b
represents the histogram for statistic (b), restricted to the
range $[0,100]$. Figure 1.c represents the histogram for statistic
(b) in the full range $[0,200]$.


\begin{center}
  \resizebox{6cm}{!}{\includegraphics{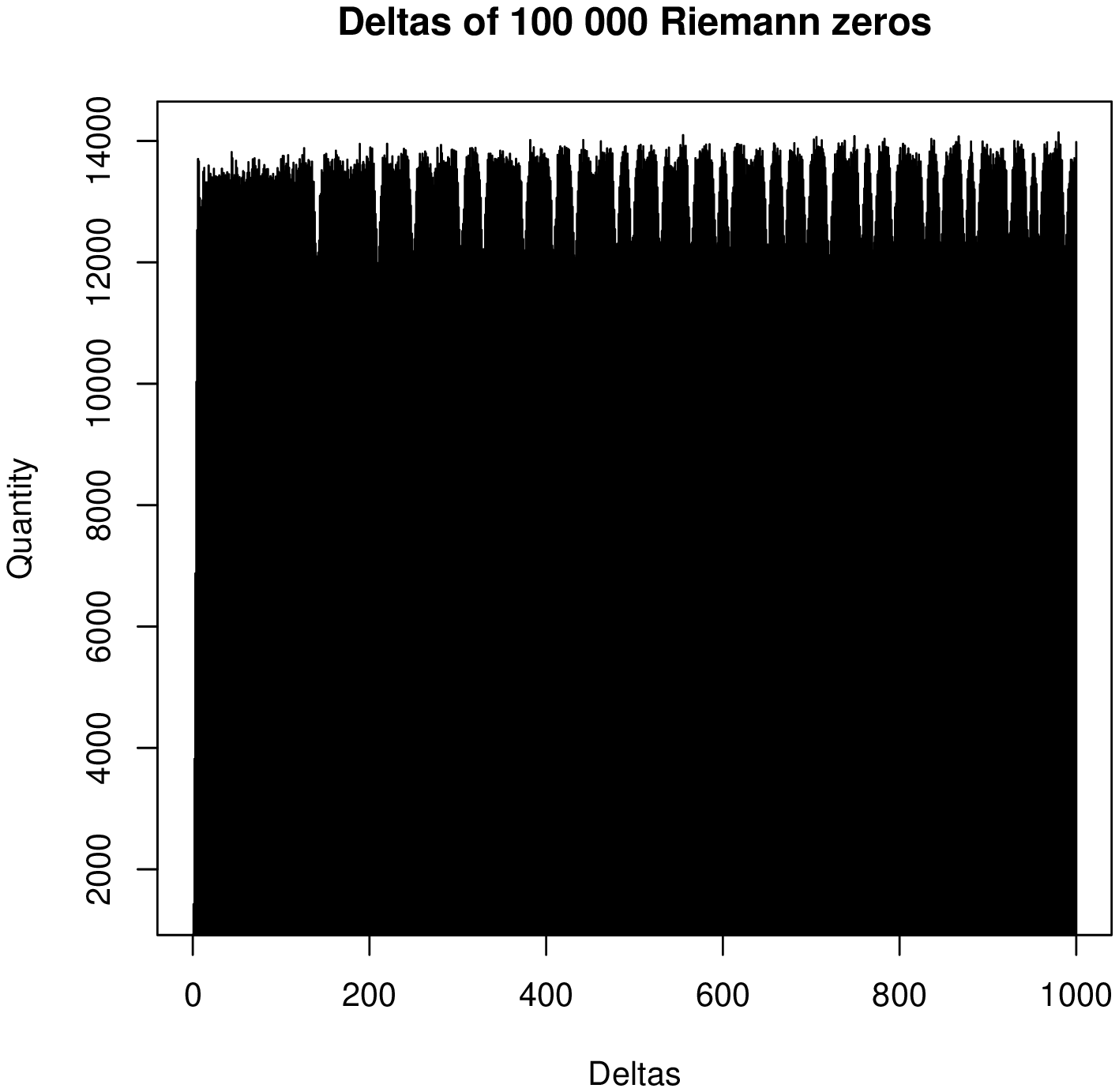}}    
  \resizebox{6cm}{!}{\includegraphics{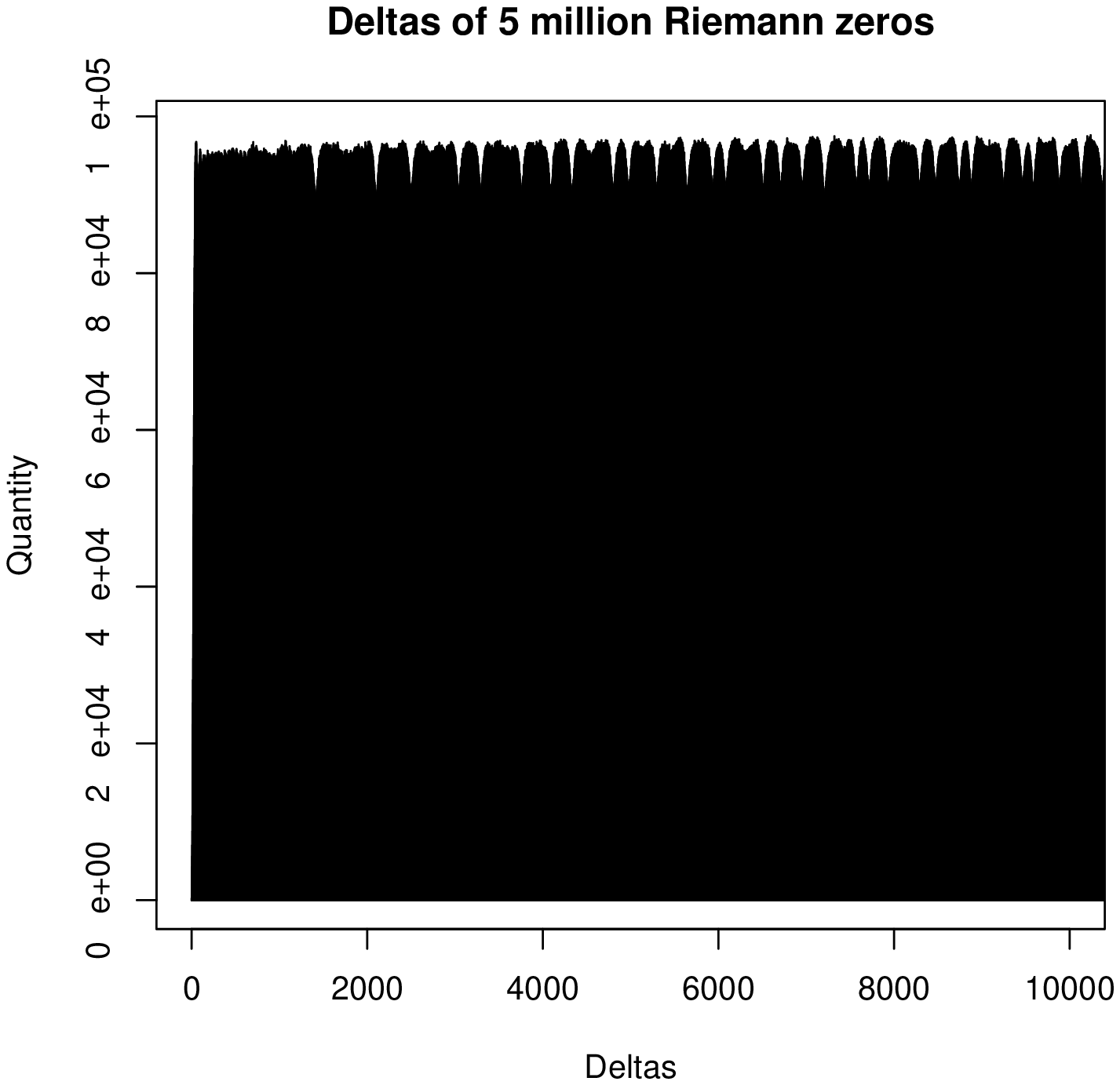}}
  \resizebox{6cm}{!}{\includegraphics{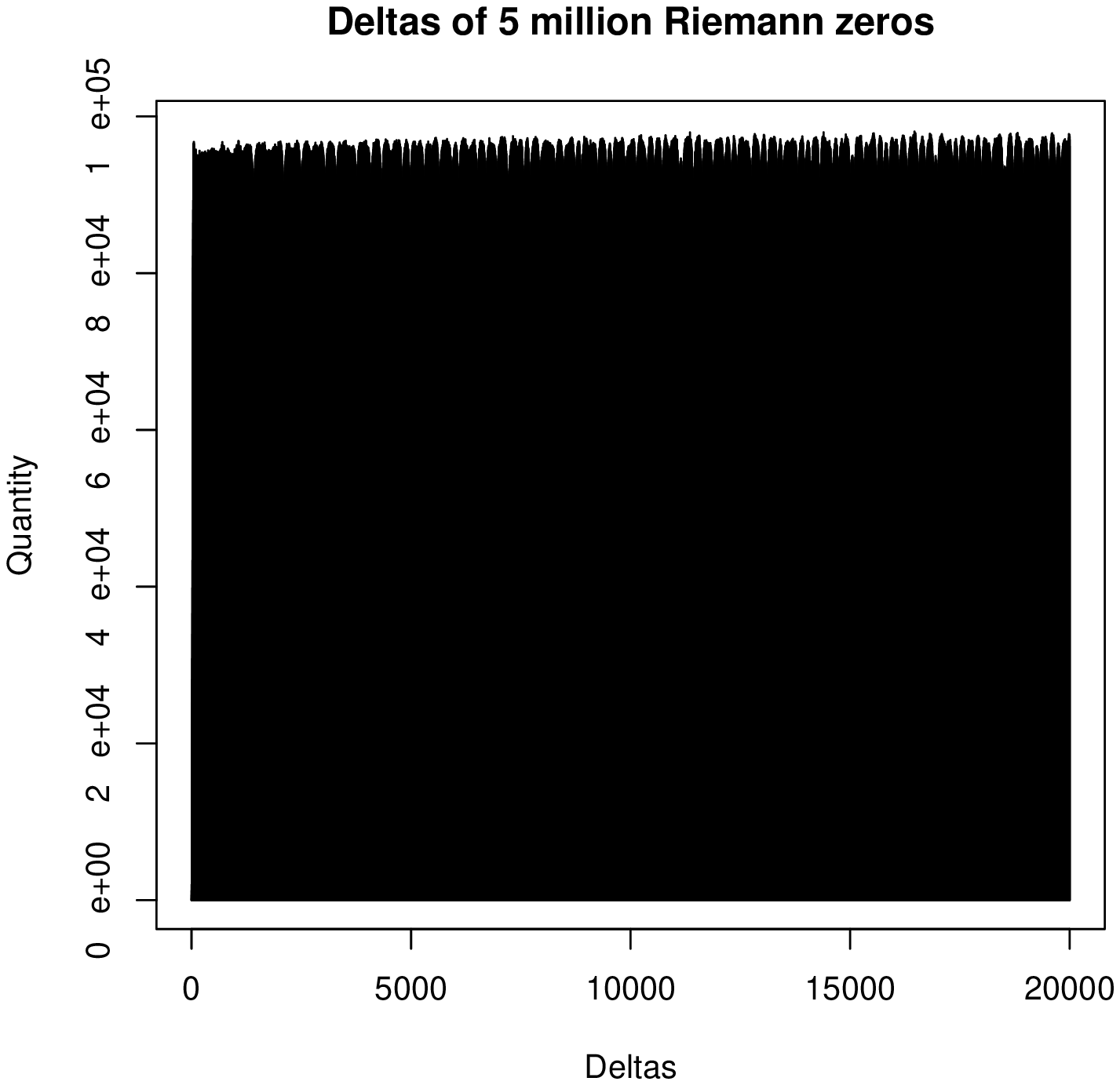}}
\end{center}

\centerline{Figures 1.a, 1.b and 1.c.}

\bigskip

At first sight this distribution seems to converge (once properly normalized) weakly to the
uniform distribution on $[0,T]$. This uniformity for the
statistics of global deltas of zeros (i.e. non-consecutive) does
not seem to have been noted explicitly in the literature. It is
hinted by the tail asymptotic to density $1$ of Montgomery's pair
correlation distribution, but it doesn't follow from that due to the
semi-locality of the differences. We discuss this later in more
detail.

\medskip

But a closer look to the histogram shows some divergence
to the uniform distribution. The convergence does not appear to be
uniform on $[0,T]$. We can notice one major deficiency for small
deltas. This appears when we zoom in the picture near $0$ (see
figures 2.a and 2.b for the range of deltas $[0,2]$).


\begin{center}
  \resizebox{6cm}{!}{\includegraphics{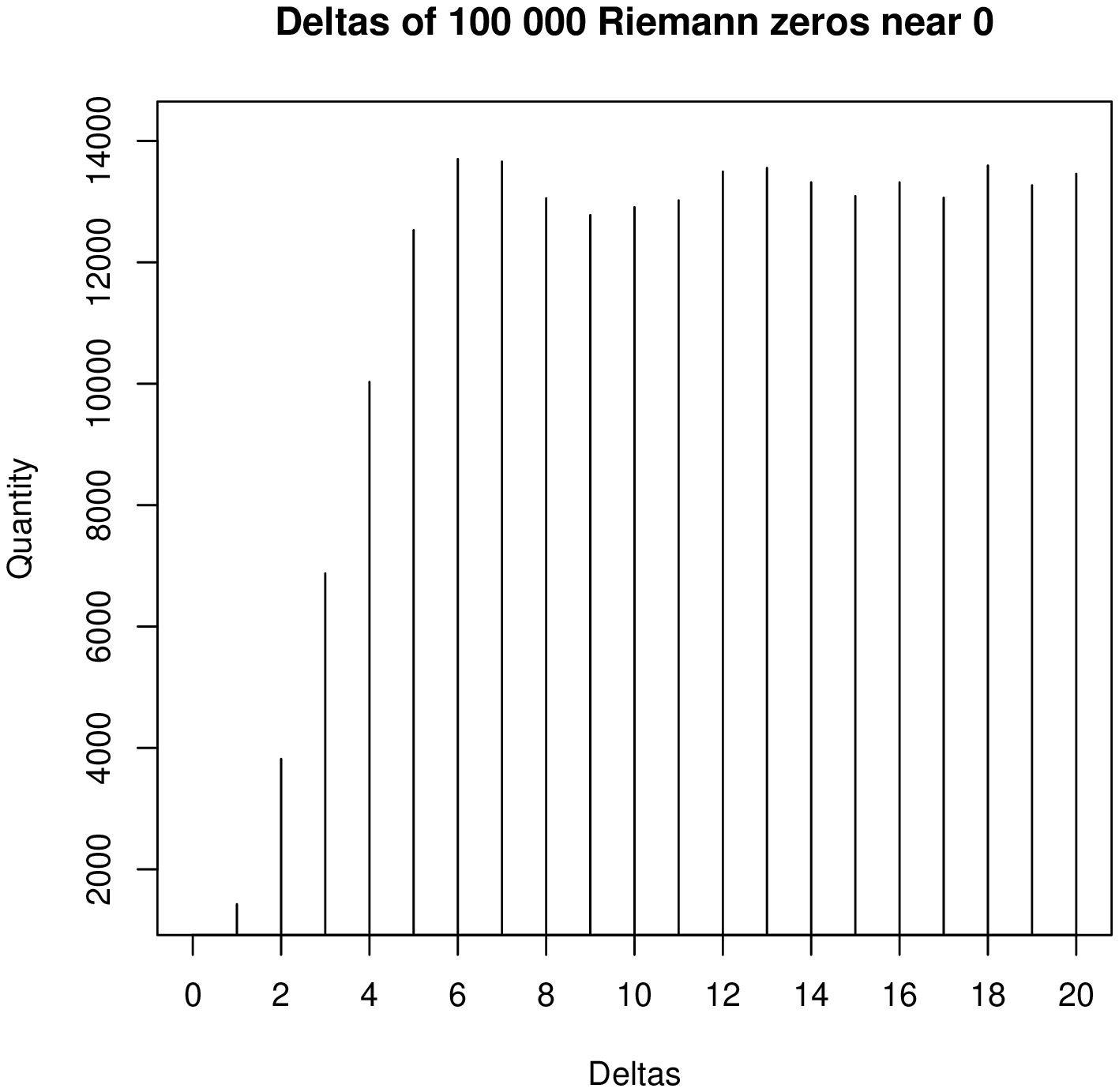}}    
  \resizebox{6cm}{!}{\includegraphics{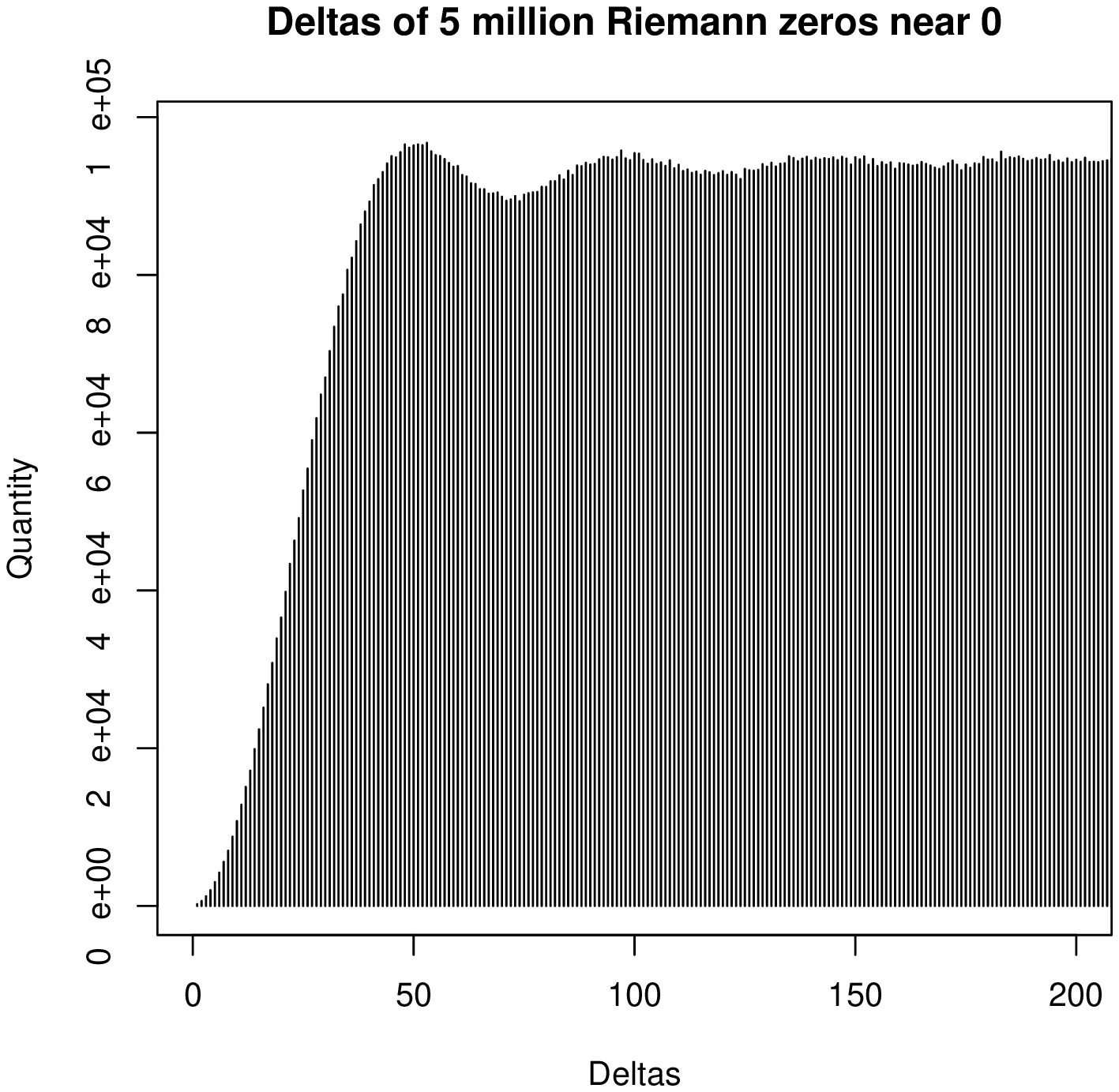}}
\end{center}

\centerline{Figures 2.a and 2.b.}

\bigskip

This deficiency is related to the observed fact that the GUE pair
correlation distribution implies that consecutive zeros tend to
repel each other. A closer look at the figures reveals a
compressed scaled GUE pair correlation distribution as expected.
The factor of compression is ${1\over 2\pi} \log T_0$ as one
should expect from Montgomery's conjecture. Notice that scaling
the figures by the factor ${1\over 2\pi }\log T_0$ (as done by those
authors studying numerically Montgomery's conjectures) pushes away
to $+\infty$ (when $T_0\to +\infty$ and $N\to +\infty$) all the
other interesting irregularities of the histogram that are the focus 
of our study.

\bigskip

Indeed, other divergences to the uniform distribution appear at
some special places distinct from $0$. We notice a remarkable
deficit of deltas at certain locations. This can be seen clearly
by zooming in at several places. Figures 3.a and 3.b shows zooms
at the interval $[10.00,30.00]$. Figures 4.a and 4.b are centered
at the interval $[30.00,50.00]$. Figure 5.a and 5.b at the
interval $[80.00,100.00]$. In all these pictures we observe at
certain precise locations noticeable negative spikes, i.e. a well
localized deficit of deltas.



\begin{center}
  \resizebox{6cm}{!}{\includegraphics{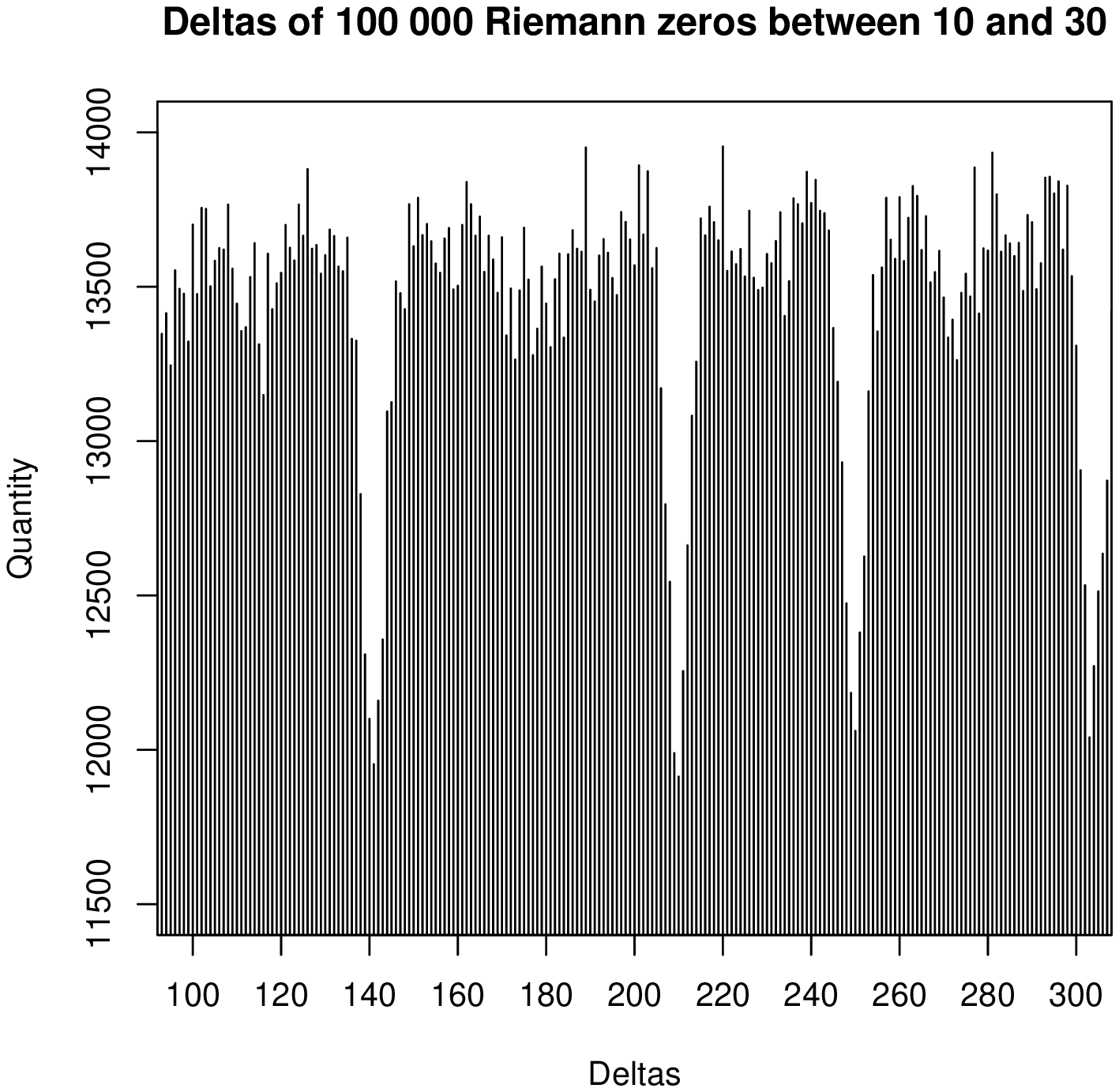}}    
  \resizebox{6cm}{!}{\includegraphics{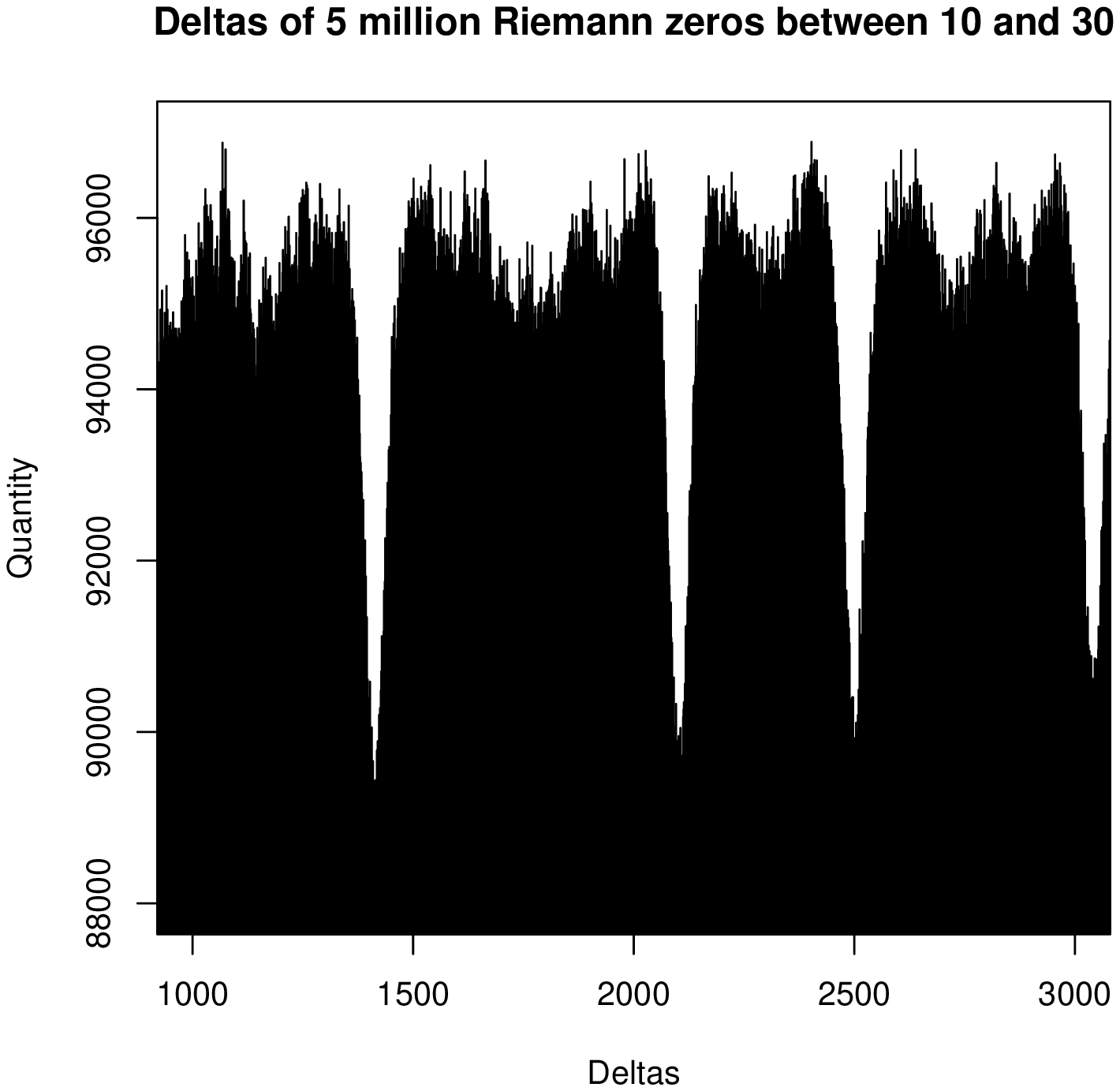}}
\end{center}

\centerline{Figures 3.a and 3.b.}



\begin{center}
  \resizebox{6cm}{!}{\includegraphics{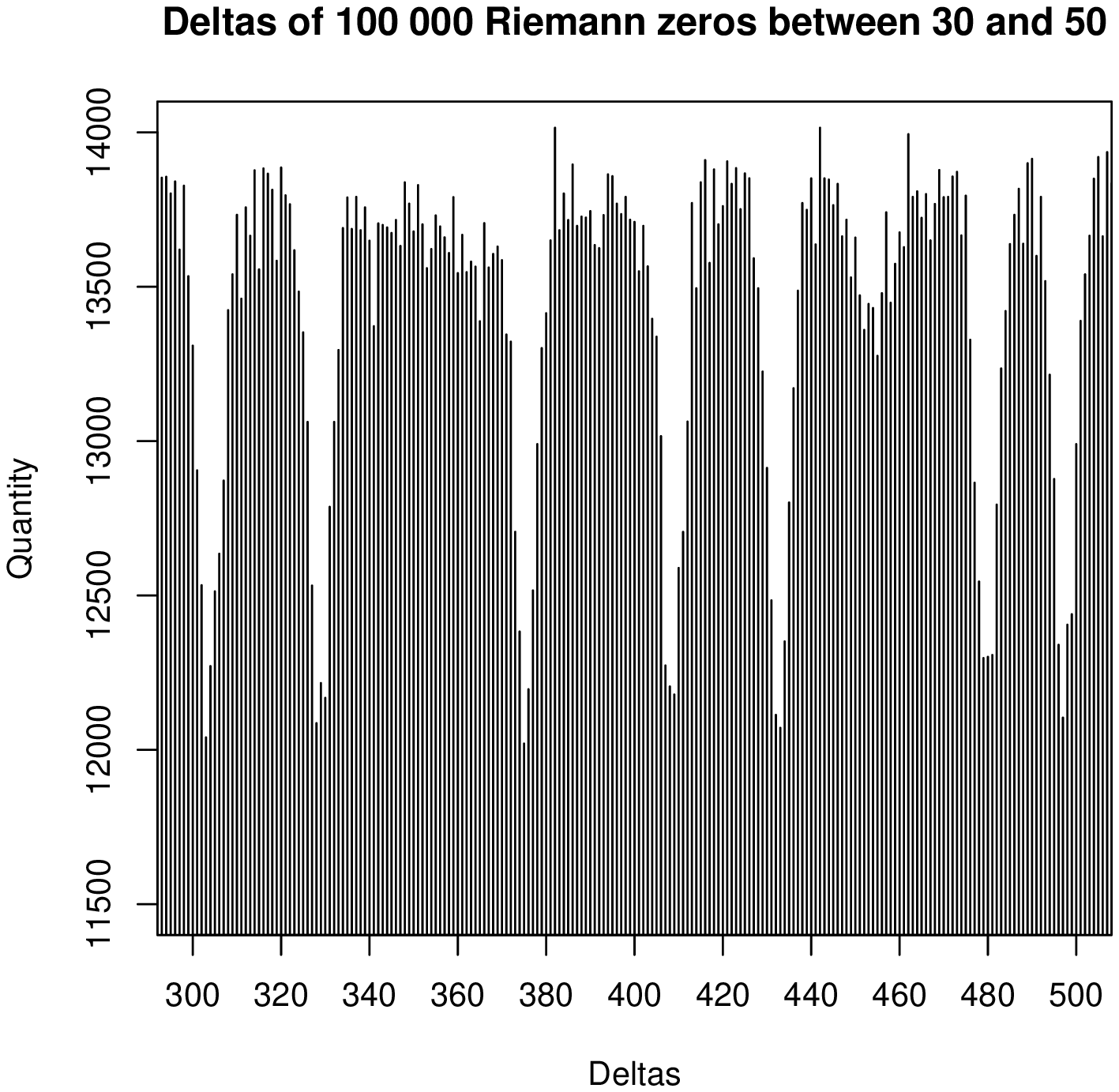}}    
  \resizebox{6cm}{!}{\includegraphics{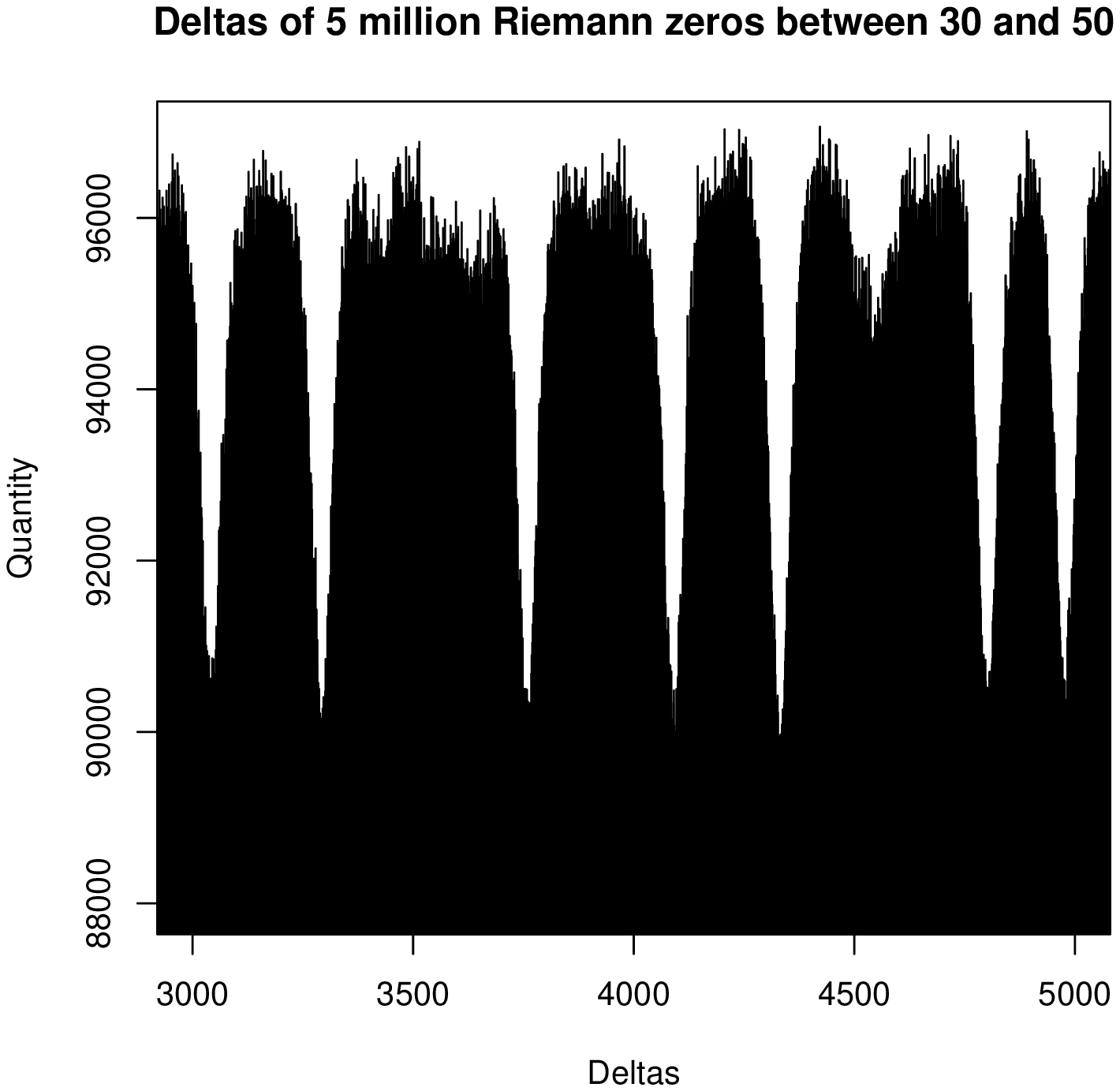}}
\end{center}

\centerline{Figures 4.a and 4.b.}



\begin{center}
  \resizebox{6cm}{!}{\includegraphics{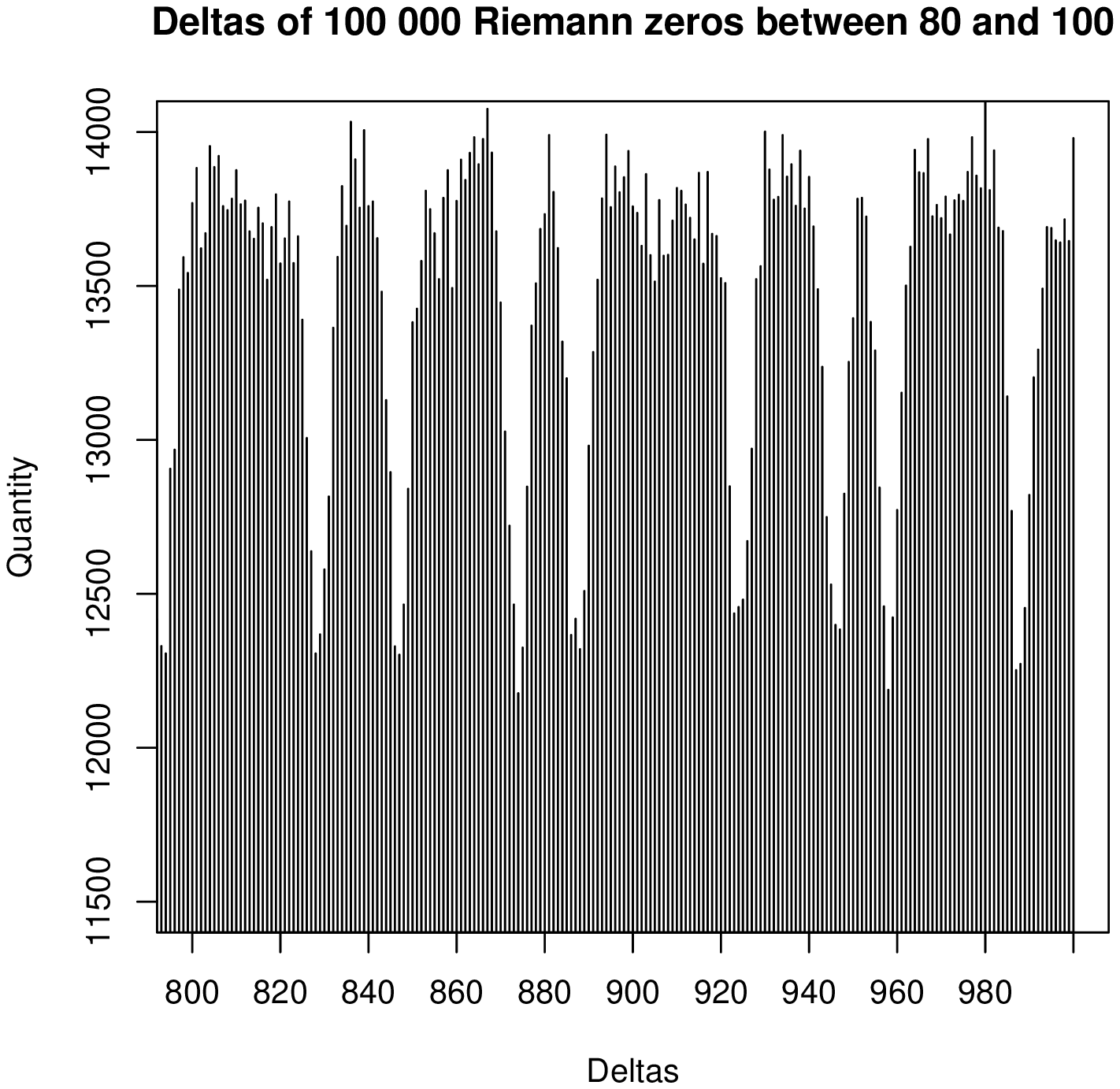}}    
  \resizebox{6cm}{!}{\includegraphics{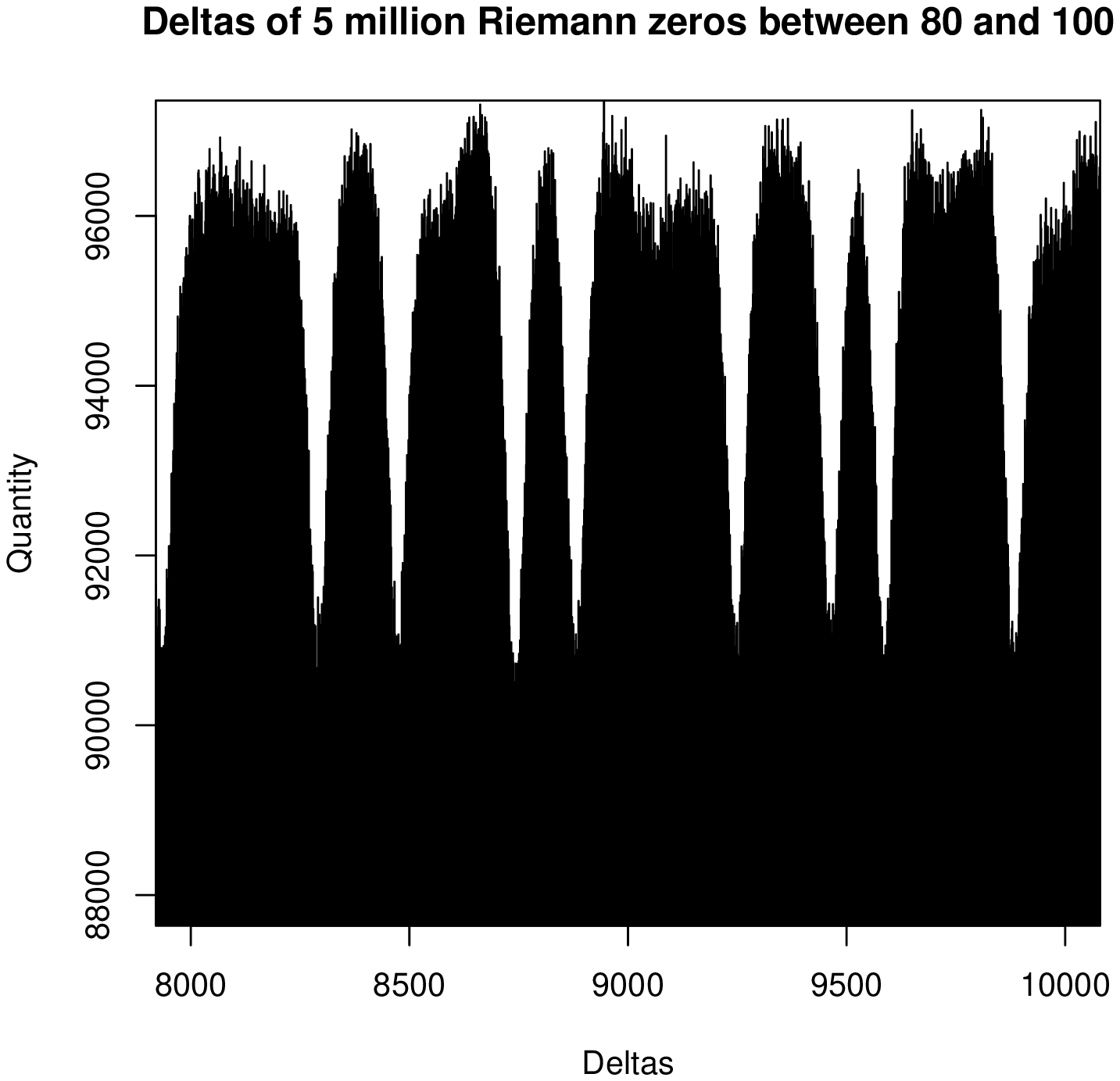}}
\end{center}

\centerline{Figures 5.a and 5.b.}

\medskip

For statistics (b) with $T=200$ we can check larger intervals.
Figures 6 and 7 are centered around the intervals
$[100.00,120.00]$ and $[190.00,200.00]$.

\begin{center}
  \resizebox{6cm}{!}{\includegraphics{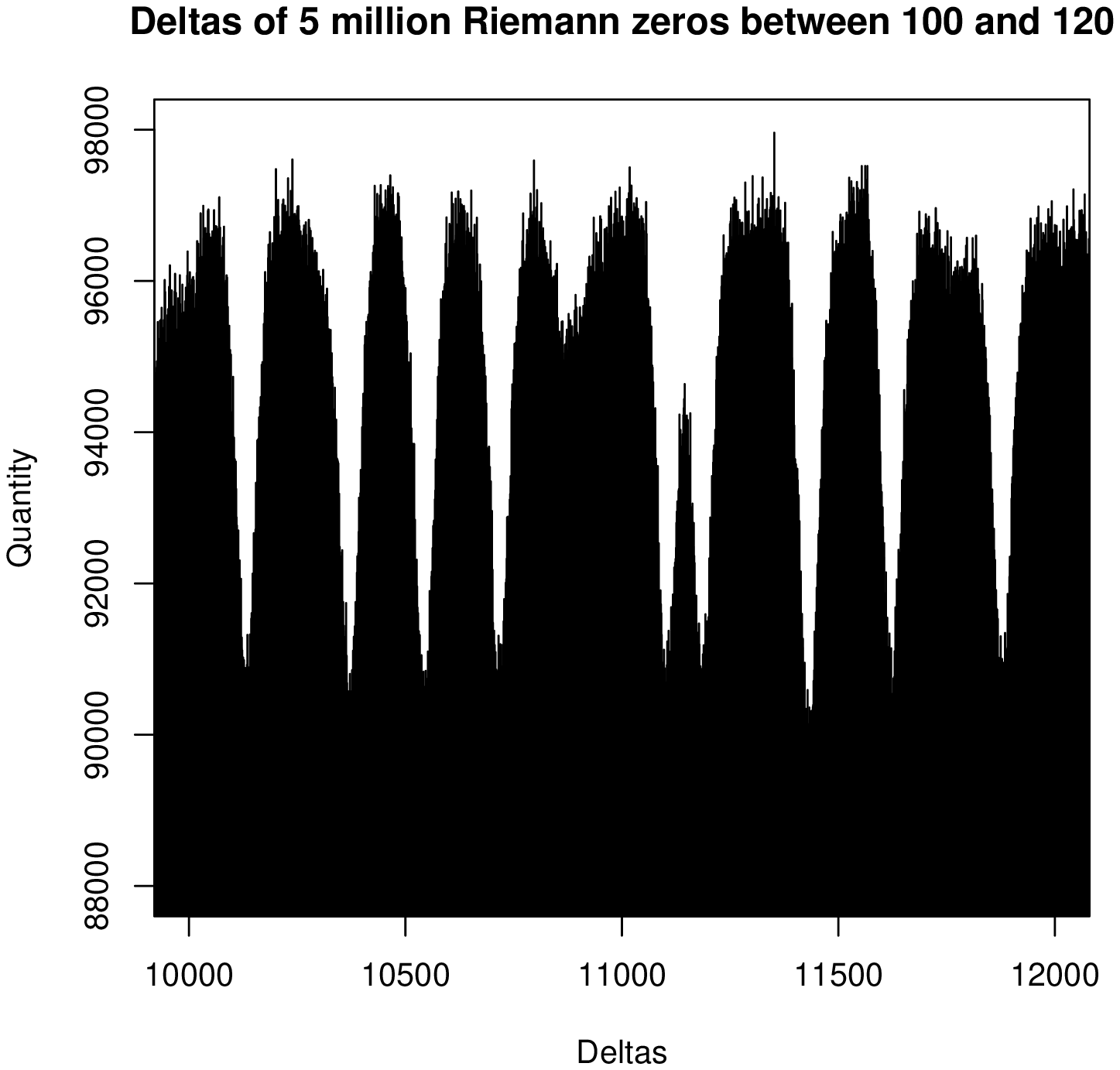}}    
  \resizebox{6cm}{!}{\includegraphics{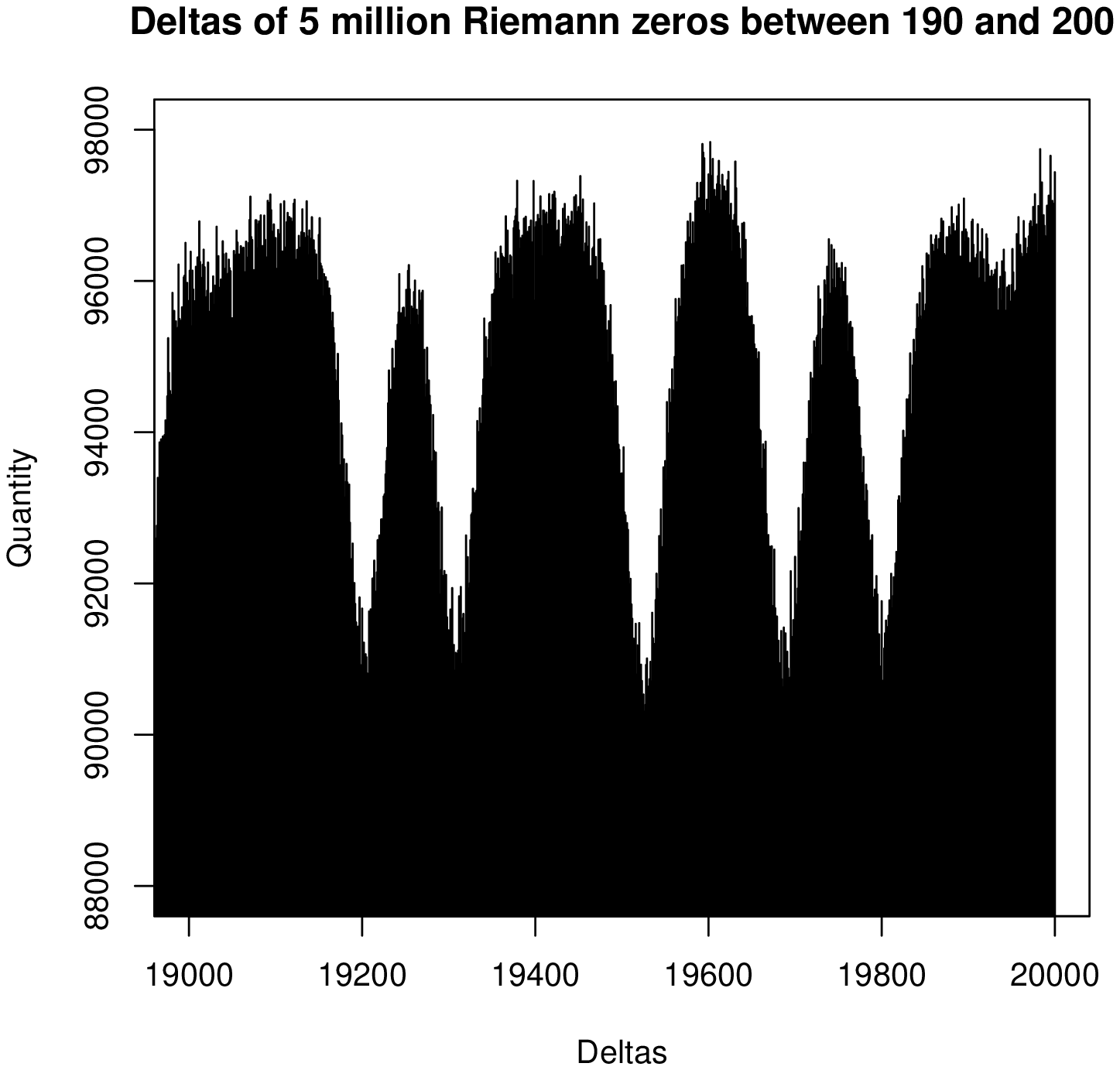}}
\end{center}

\centerline{Figures 6 and 7. Statistics (b).}

%
%
%

At this point the reader should take a moment
and compare these pictures, and in particular the location of the
deficiencies, with the tabulated list of Riemann zeros.

The key observation now is that the location of these
negatives spikes is truly special. These locations are precisely
at the very same location of the Riemann zeros. We recognize in
figures 3.a and 3.b the locations of the 4 first Riemann zeros. In
all the Figures 3a, 3b, 4a, 4b, 5a, 5b, 6 and 7 we recognize the location of the
zeros in the corresponding intervals. Note in particular in figure
6 the two nearby zeros near the value $111$,

$$
\g_{34} =111.029535543\ldots \ \ \ \ \ \g_{35}=111.874659177\ldots
$$

\medskip

We conclude that

\medskip

\centerline{\textbf{Riemann zeros do repel their deltas.}}

\medskip

This property of the sequence of Riemann zeros is even more
surprising considering the fact that it is not invariant by
translation, i.e. by a global translation of the sequence. The set
of deltas is independent of such a translation, but obviously not this
property. The location of each zero is well determined. Any
variation on the location of a single zero is obviously irrelevant
for the distribution of deltas, but the zero will then miss the
location of the negative spike. Therefore, only the statistics of
the deltas determines the precise location of the zeros. This
implies that \textbf{any subsequence of density $1$ of Riemann zeros
does determine the whole sequence}. For this reason we name this
property \textbf{the self replicating property of the zeros}.


The self replicating property of the Riemann zeros is
completely mysterious without the motivation that lies behind this 
numerical study: The theory of the e\~ne product.

\medskip

We can confirm numerically these observations (in statistics (a)
for simplicity) by noticing that in the histogram all the deficit
values with cumulative count inferior to $12\ 500$ fall near a
Riemann zero, and conversely any Riemann zero yields a group of
deficit values. The list of the values of $k$ for which $x_k< 12 \
500$ in statistics (a) is the following:  1, \ 2,\    3,\   4,\
5,\ 139,\ 140,\ 141,\ 142,\ 143,\ 208,\ 209,\ 210,\ 211,\ 212,\
248,\ 249,\
 250,\ 251, \ 252, 302,\ 303,\ 304,\ 305,\ 306,\  327,\ 328,\ 329,\ 330,\ 374,\ 375,\
376,\ 377,\ 407,\ 408,\ 409, \ 410,\ 431, 432, \ 433,\ 434,\ 478,\
479,\ 480,\ 481,\ 496,\ 497,\ 498,\ 499,\ 528,\ 529,\ 530,\ 531,\
563,\ 564,\ 565, 566,\ 592,\ 593,\ 594,\ 606,\ 607,\ 608,\ 609,\
649,\ 650,\ 651,\ 652,\ 669,\ 670,\ 671,\ 694,\ 695,\ 696, 718,\
719,\ 720,\ 721,\ 722,\ 755,\ 756,\ 757,\ 758,\ 770,\ 771,\ 772,\
792,\ 793,\ 794,\ 827,\ 828,\ 829, 830,\ 846,\ 847,\ 848,\ 873,\
874,\ 875,\ 886,\ 887,\ 888,\ 889,\ 923,\ 924,\ 925,\ 945,\ 946,\
947,\ 957, 958,\ 959,\ 987,\ 988,\ 989.

\medskip

To make the main observation more precise, we can, for example,
average out all deficit values in each group. We discard the first
group of values near $0$ that corresponds to the deficit at $0$
(we will come back to this). Then we find out as many
groups as Riemann zeros and their averages are denoted by $\bar
\g_1, \bar \g_2, \bar\g_3,\ldots$ They are all very close to the
corresponding zero. Table 6 compares the sequence of zeros $(\g_i
)$ with the sequence of averages $(\bar \g_i)$ for all $29$ zeros
less than $100$. We rounded up the averages to the first decimal.
The matching of the averages $\bar \g$ with the zeros $\g$ is
striking.

$$
\vbox{ \offinterlineskip \halign{ \strut \vrule $#$ &\vrule $#$
&\vrule $#$ \vrule \cr \noalign{\hrule} {\ \hbox { i}} & \ \g_i &
\ \bar \g_i \cr \noalign{\hrule} \ 1 \ & \ 14.134725142 \ & \ 14.1
\cr \noalign{\hrule} \ 2 \ & \ 21.022039639\ & \ 21.0 \cr
\noalign{\hrule} \ 3 \ & \ 25.010857580\ & \ 25.0 \cr
\noalign{\hrule} \ 4 \ & \ 30.424876126\ & \ 30.4 \cr
\noalign{\hrule} \ 5 \ & \ 32.935061588\ & \ 32.9 \cr
\noalign{\hrule} \ 6 \ & \ 37.586178159\ & \ 37.6 \cr
\noalign{\hrule} \ 7 \ & \ 40.918719012 \ & \ 40.9 \cr
\noalign{\hrule} \ 8 \ & \ 43.327073281\ & \ 43.3 \cr
\noalign{\hrule} \ 9 \ & \ 48.005150881\ & \ 48.0 \cr
\noalign{\hrule} \ 10 \ & \ 49.773832478 \ & \ 49.8 \cr
\noalign{\hrule} \ 11 \ & \ 52.970321478 \ & \ 53.0 \cr
\noalign{\hrule} \ 12 \ & \ 56.446247697 \ & \ 56.5 \cr
\noalign{\hrule} \ 13 \ & \ 59.347044003 \ & \ 59.3 \cr
\noalign{\hrule} \ 14 \ & \ 60.831778525 \ & \ 60.8 \cr
\noalign{\hrule} \ 15 \ & \ 65.112544048 \ & \ 65.1 \cr
\noalign{\hrule} \ 16 \ & \ 67.079810529 \ & \ 67.0 \cr
\noalign{\hrule} \ 17 \ & \ 69.546401711 \ & \ 69.5 \cr
\noalign{\hrule} \ 18 \ & \ 72.067157674 \ & \ 72.0 \cr
\noalign{\hrule} \ 19 \ & \ 75.704690699 \ & \ 75.7 \cr
\noalign{\hrule} \ 20 \ & \ 77.144840069 \ & \ 77.1 \cr
\noalign{\hrule} \ 21 \ & \ 79.337375020 \ & \ 79.3 \cr
\noalign{\hrule} \ 22 \ & \ 82.910380854 \ & \ 82.9 \cr
\noalign{\hrule} \ 23 \ & \ 84.735492981 \ & \ 84.7 \cr
\noalign{\hrule} \ 24 \ & \ 87.425274613 \ & \ 87.4 \cr
\noalign{\hrule} \ 25 \ & \ 88.809111208 \ & \ 88.8 \cr
\noalign{\hrule} \ 26 \ & \ 92.491899271 \ & \ 92.4 \cr
\noalign{\hrule} \ 27 \ & \ 94.651344041 \ & \ 94.6 \cr
\noalign{\hrule} \ 28 \ & \ 95.870634228 \ & \ 95.8 \cr
\noalign{\hrule} \ 29 \ & \ 98.831194218 \ & \ 98.8 \cr
\noalign{\hrule} } }
$$

{\centerline {\textbf{ Zeros versus group averages.}} }

\medskip

It is also interesting to study the structure of the distribution 
of the deficit of
deltas near the zeros. Once properly scaled, we observe a
universal distribution for all zeros. This distribution is a
negative Fresnel distribution, i.e. the distribution generated by
the Fresnel integral (also named sine integral, or sine cardinal
function)
$$
{\hbox {\rm sinc}}_\pi (x)= \left \{ \begin{array}{l l}
    \frac{\sin(\pi x)}{\pi x} & \ \ {\hbox {\rm for }} \ x\not=0 \\
    \ \ 1 \ \ \ & \ \ {\hbox {\rm for }} \ x=0 \\
  \end{array} \right.
$$

We have
$$
\int_{-\infty}^{+\infty } {\hbox {\rm sinc}}_\pi (x) \ dx =1 \ .
$$
The Fresnel distribution is the Fourier transform
$$
{\hbox {\rm sinc}}_\pi (x)=\int_{\RR} e^{-2\pi i x t } \ \Pi (t)
\ dt \ ,
$$
of the box function,
$$
\Pi (x)=H(x+1/2)-H(x-1/2)=\left \{ \begin{array}{l l}
    0 & \ \ {\hbox {\rm for }} \ |x| >1/2 \ , \\
    1 & \ \ {\hbox {\rm for }} \ |x| <1/2 \ . \\
  \end{array} \right.
$$

We can appreciate this for the histogram plotted in Figure 8.a. Figure
8.b shows a more intensive computation with deltas of the first 10
million zeros.


\begin{center}
  \resizebox{6cm}{!}{\includegraphics{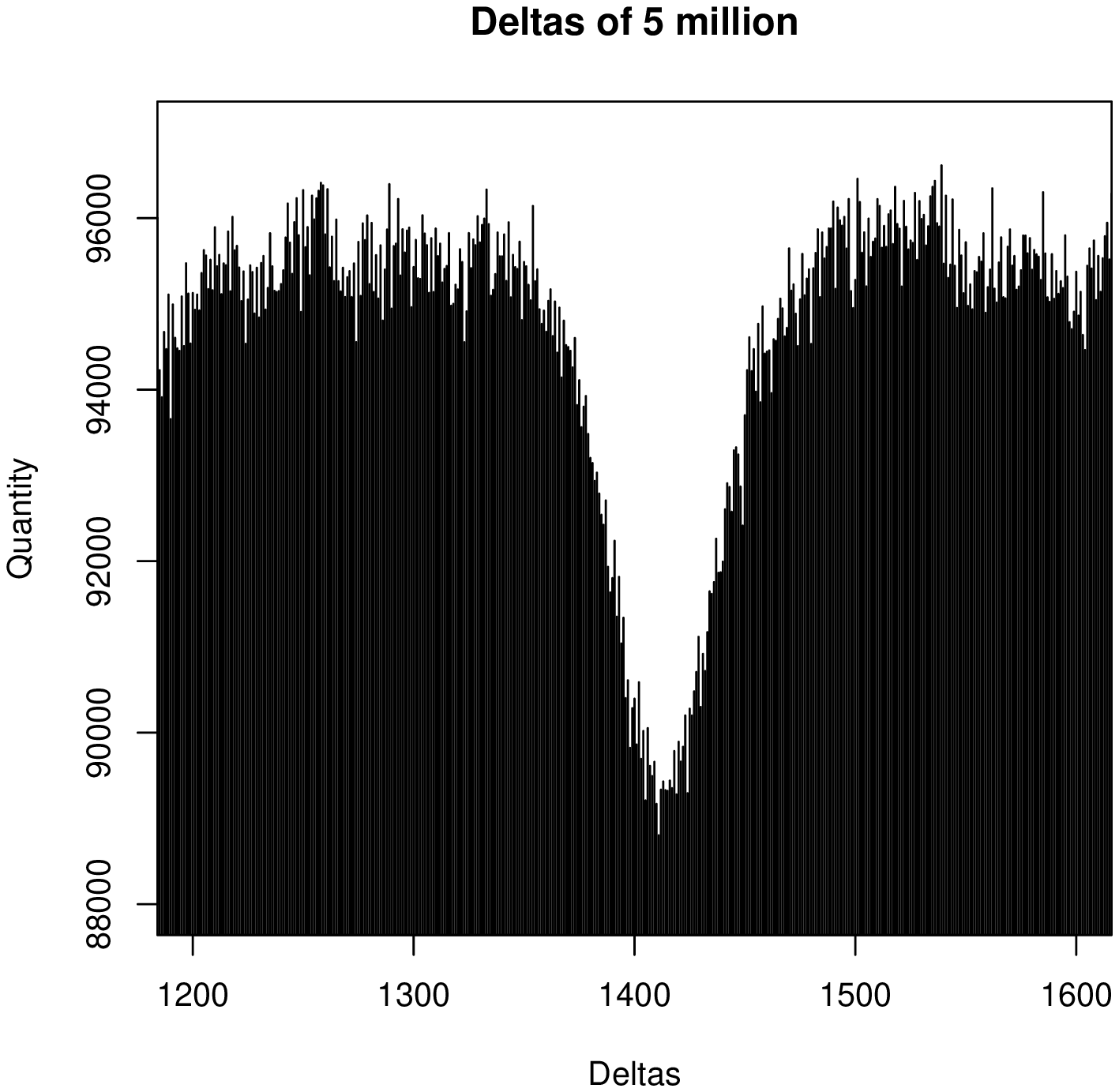}}    
  \resizebox{6cm}{!}{\includegraphics{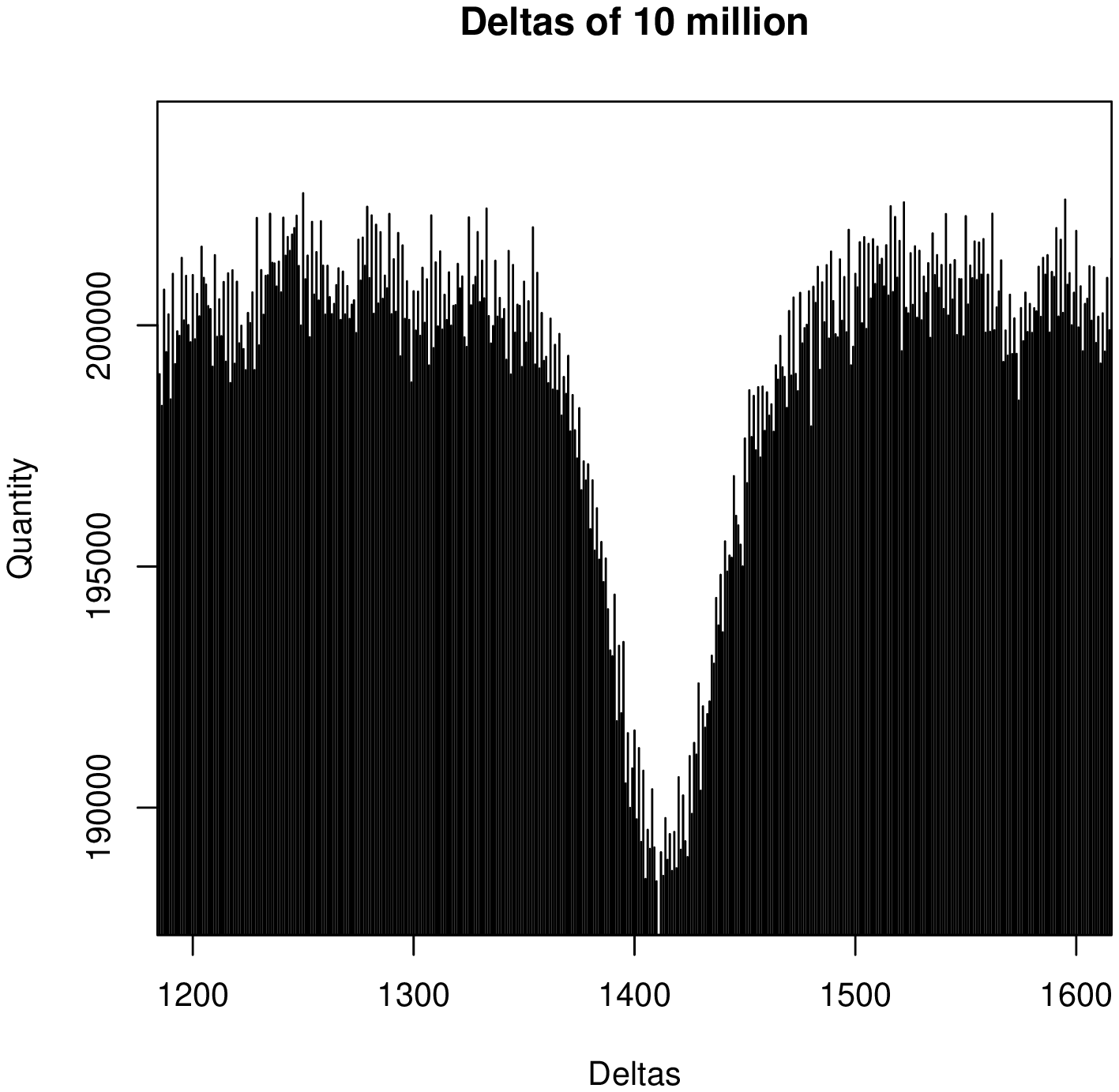}}
\end{center}

\centerline{Figures 8.a and 8.b.}

In the figures we appreciate a higher frequency noise that blurs
the picture. We can filter the noise out by standard filtering procedures. The
simplest one would be to replace (for example) the sequence
$(x_k)$ by the sequence $(fx_k)$ where

$$
fx_k={1\over \tau } \sum_{i=-\tau/2}^{\tau/2} x_k \ .
$$
The new figures 9.a and 9.b show the pictures with the noise
filtered.


\begin{center}
  \resizebox{6cm}{!}{\includegraphics{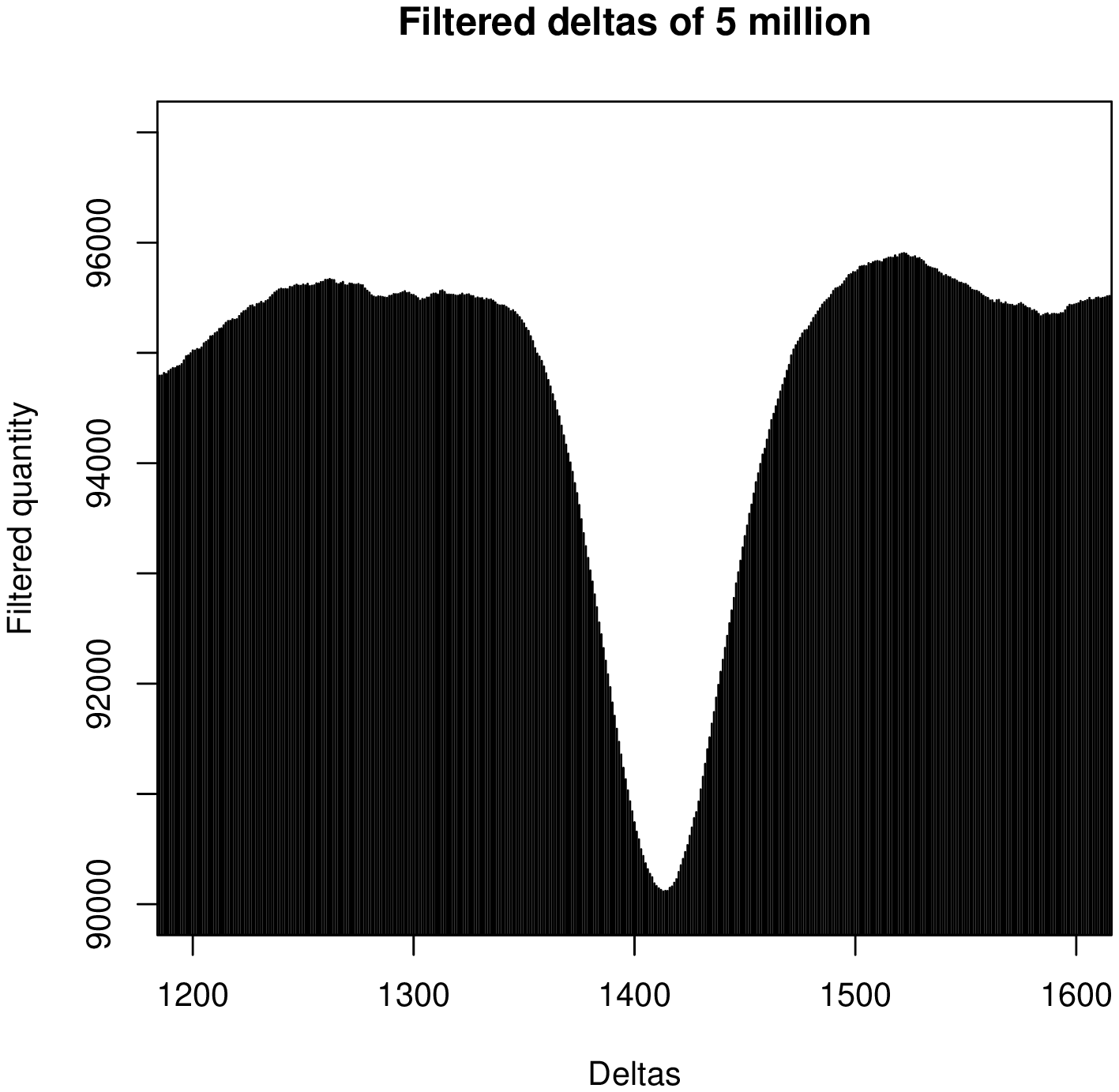}}    
  \resizebox{6cm}{!}{\includegraphics{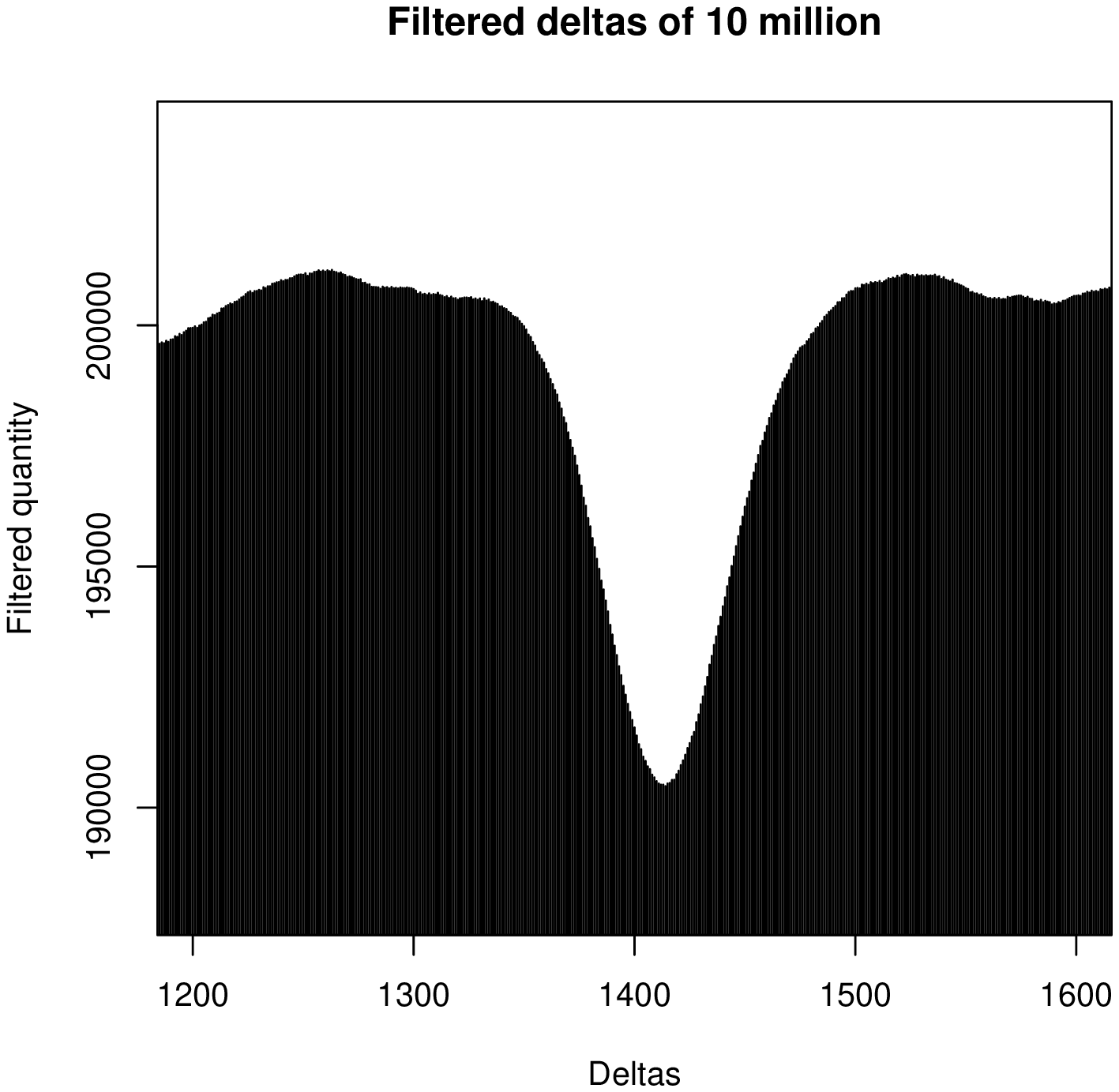}}
\end{center}

\centerline{Figures 9.a and 9.b.}

\medskip

\textbf{E\~ne product computation.}

\medskip

These numerical observations come from the analytic divisor interpretation of the e\~ne product, noted $\bar \star$, 
developped in \cite{PM1} and \cite{PM2}. The space of Dirichlet $L$-functions endowed with usual the multiplication 
and the e\~ne product is a commutative ring, having a proper normalization of the Riemann zeta function as the
e\~ne-multiplicative unit. The e\~ne product is associative, not only with respect to multiplication, but also to 
infinite arithmetic Euler products. Given two Euler products with polynomials $(F_p)$ and $(G_p)$ (with $F_p(0)=G_p(0)=1$),

\begin{align*}
F(s) &= \prod_p F_p(p^{-s}) \\
G(s) &= \prod_p G_p(p^{-s}) 
\end{align*}
then
\begin{equation*}
F\ \bar \star \ G (s) = \prod_p F_p \star G_p (p^{-s}) \ ,
\end{equation*}
where $F_p\star G_p$ is the plain e\~ne product in $\CC$ of polynomial whose zeros are the product of the zeros of $F_p$ with the zeros 
of $G_p$, i.e. if
\begin{align*}
F_p(X) &=\prod_\alpha \left (1-\frac{X}{\alpha } \right ) \\ 
G_p(X) &=\prod_\beta \left (1-\frac{X}{\beta } \right ) 
\end{align*}
then
\begin{equation*}
 F_p\star G_p=\prod_{\alpha , \beta} \left (1-\frac{X}{\alpha \beta} \right ) \ ,
\end{equation*}
in particular
\begin{equation*}
(1-ap^{-s}) \ \bar \star \ (1-bp^{-s}) = 1-ab p^{-s} \ .
\end{equation*}
The main arithmetic property is that for $p\not= q$, we have $\log p$ and $\log q$ $\QQ$-independent, and 
\begin{equation*}
F_p(p^{-s}) \ \bar \star \ G_q(q^{-s}) = 1 \ .
\end{equation*}

Now we denote that for a real analytic function $F$,
$$
\bar F(s)={\overline {F(\bar s)}}=F(s) \ .
$$

 The main statistics in this section have its origin in the
following computation

\begin{align*}  \zeta (s) \ \bar \star \ {\overline {\zeta}}(s)&= \zeta
(s) \ \bar \star \ \zeta(s) \\ 
&= \left ( \prod_p (1- p^{-s})^{-1} \right ) \bar \star \left ( \prod_q (1- q^{-s} )^{-1}
\right ) \\ 
&=\prod_p (1- p^{-s} )^{-1} \ \bar \star \ (1- p^{-s}
)^{-1} \\ &=\prod_p (1- p^{-1/2} p^{-s} ) \\ 
&=\zeta(s+1/2)^{-1} \ . 
\end{align*}

\textbf{R script.}

\medskip

The following script can be directly feeded into R in order to
plot the histogram for the deltas for 5 million zeros with
precision $10^{-2}$ (statistics (b)). The zeros are read from the
file "zero.data"(one zero per line in increasing order). The
reader can consult the R tutorial for more elaborate plotting
commands. The scripts for the other statistics are simple
modifications from this one.

\medskip

{\tt

zero1<-scan("zero.data", nlines=5000000)

zero2<-zeta1

x=rep(0,10000)

N=5000000

k=0

for (i in 1:N)

$\{$

while (((zero1[i]-zero2[i+k])<100.01) \& (k+i>1))

$\{$

 k<-k-1

$\}$

k=k+1

j=k

while ( (zero1[i]-zero2[i+j]>0) \& (zero1[i]-zero2[i+j]<100.01) )

$\{$

d=100*(zero1[i]-zero2[i+j]) x[as.integer(d)]=x[as.integer(d)]+1

j=j+1

$\}$

$\}$

barplot(x)

 }

\section{Large Riemann zeros know about all zeros.} \label{sec:large}

In this section we perform the same statistics as in section 1,
but only using deltas of large zeros. The convergence is
slower, but the results are the same. This indicates that zeros
with large imaginary part contain full information on the
location of all zeros. Indeed a density $1$
proportion of zeros with large imaginary part contains the
information about the location of all zeros.

We perform the statistics with two sets of data. The first one,
using Rubinstein's file for the first 35 million
zeros and selecting one million after the 30th million zero. The
second and the third are
performed with a much larger set of zeros using Odlyzko's file
containing $10 \ 000$ zeros after $10^{12}$ and after $10^{21}$
respectively. The number of zeros in  these last two statistics is
insufficient. These sets of $10 \ 000$ zeros are small and the
distribution of deltas is not even close to the uniform
distribution. In the first Odlyzko's file all deltas are smaller
than $2568$ and in the second smaller than $1409$. Therefore we
observe a linear deficit of deltas even for small values of delta
when delta increases. In order to pinpoint the deficit at the
zeros we filter the cumulative data on deltas by removing a moving
average. Although statistically not as significant as the other
statistics, the deficit phenomena is still clearly visible at the
location of the zeros. It is also less visible for large zeros.
This indicates a slower convergence.

The following figures 10a, 10b and 10c are from the first statistics with
10 million zeros $(\g_i)$ with $20.10^6<i\leq 30.10^6$. Figures
10b and 10c show the details near the top of the uniform distribution. 
We can appreciate the similarity of these pictures with the
previous ones. Figure 10c is centered around the segment $[10,30]$
and is almost identical to figures 3a and 3b.


\begin{center}
  \resizebox{6cm}{!}{\includegraphics{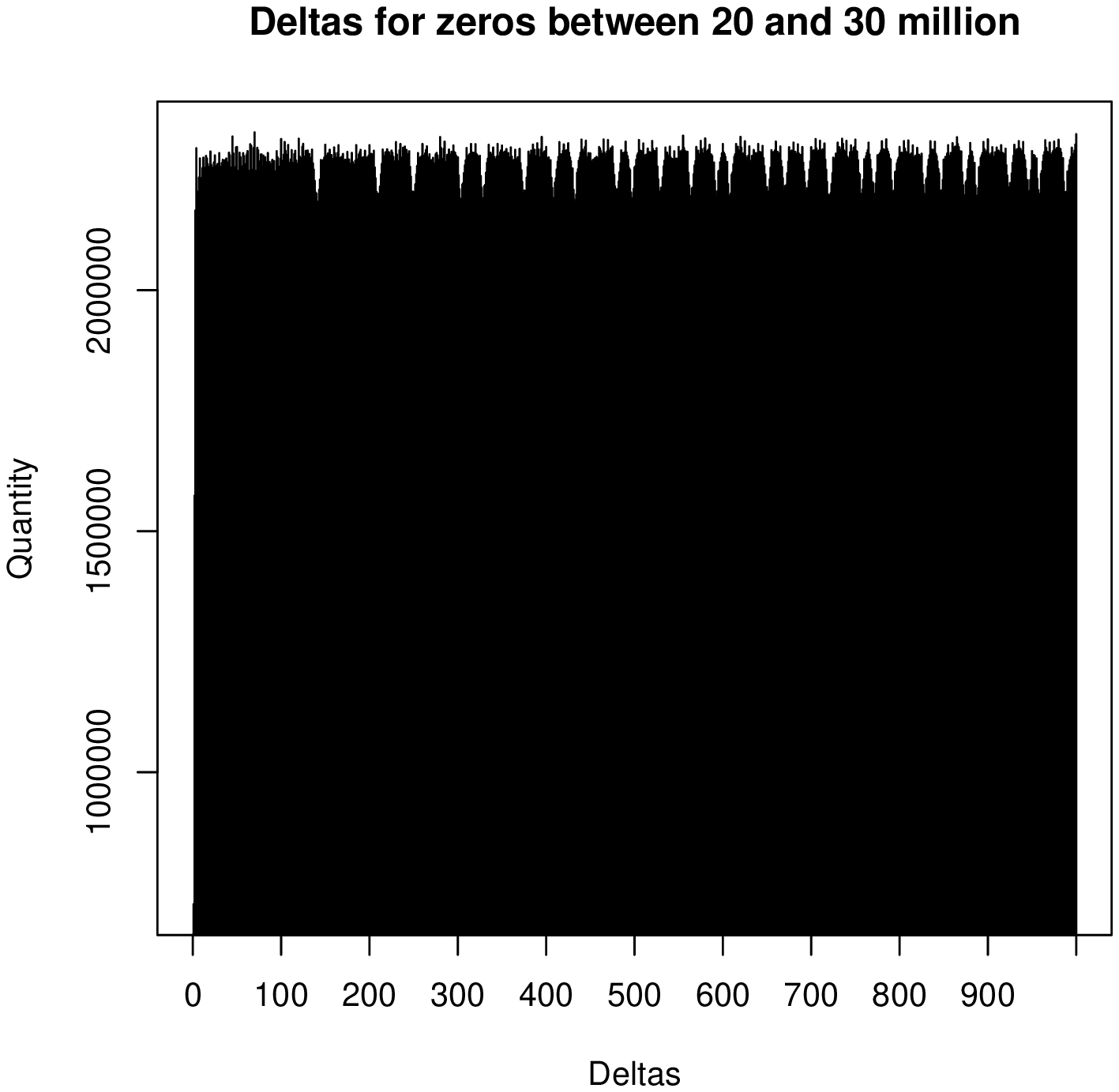}}    
  \resizebox{6cm}{!}{\includegraphics{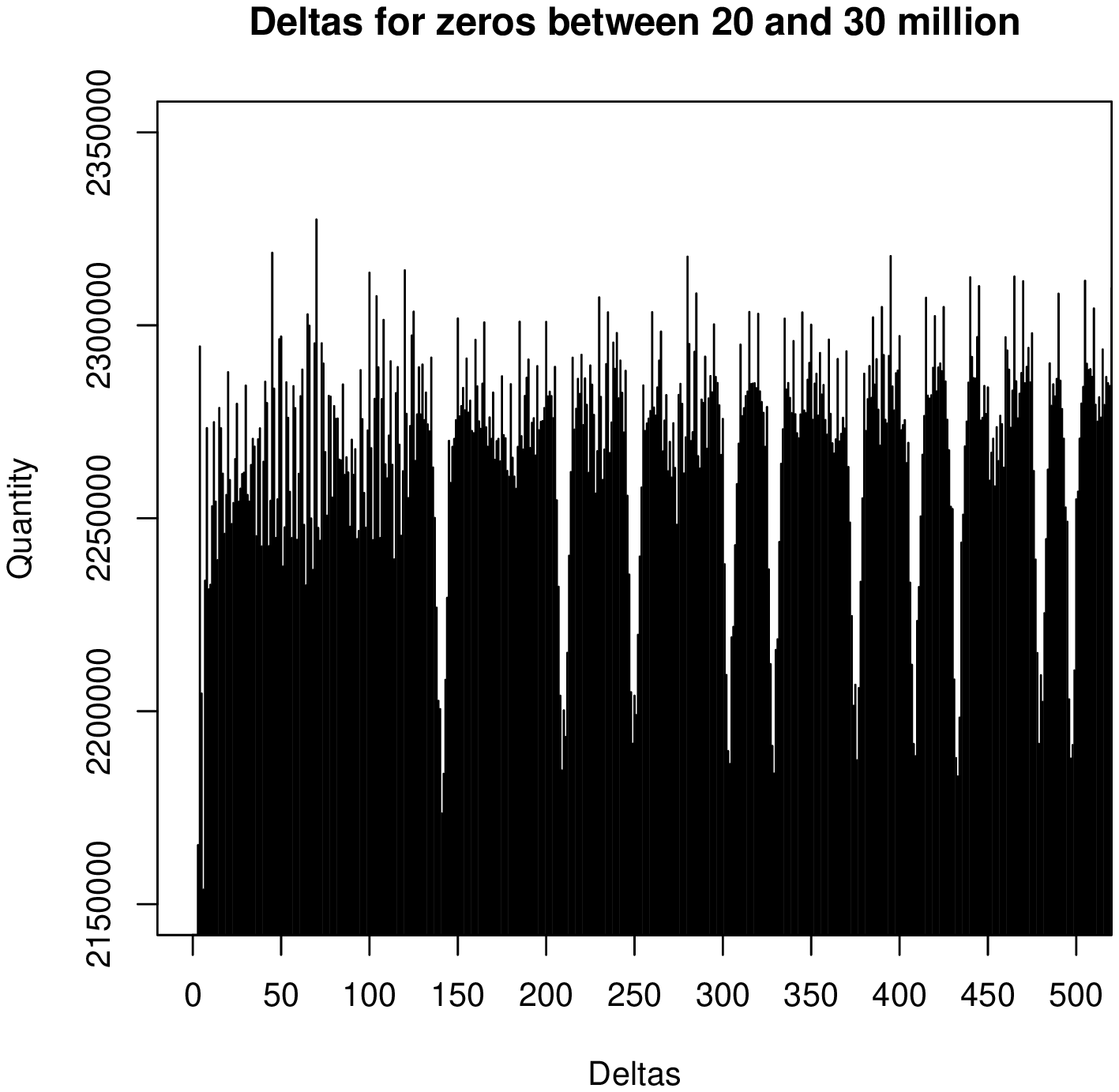}}
  \resizebox{6cm}{!}{\includegraphics{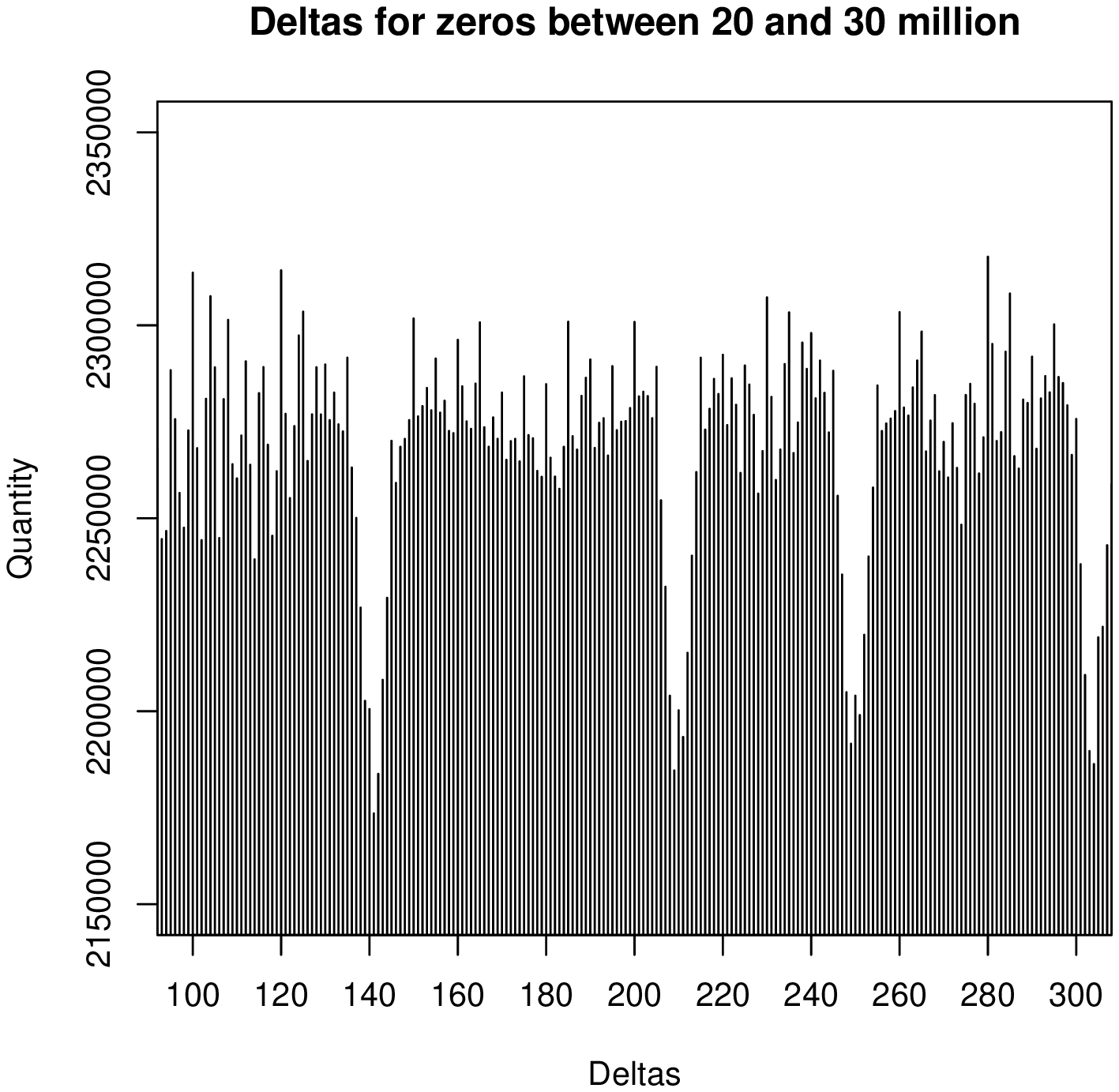}}
\end{center}

\centerline{Figures 10.a, 10.b and 10.c.}

\bigskip

In the following figures we illustrate the results for Odlyzko's
large zeros near $10^{12}$ (statistics (a)) and near $10^{21}$
(statistics (b)). We can observe the linear decreasing of the amount of deltas due
to the small number of zeros used. We worked with Odlyzko's files
containing only $10 \ 000$ zeros. Paying close attention we can
discern the deficit of deltas at the location of the zeros. This
can be better seen by filtering the data by removing a moving
average. Figures 12.a and 12.b show that. In Figures 13.a and 13.b we
have the details for deltas smaller than $50$.

%
%

\begin{center}
  \resizebox{6cm}{!}{\includegraphics{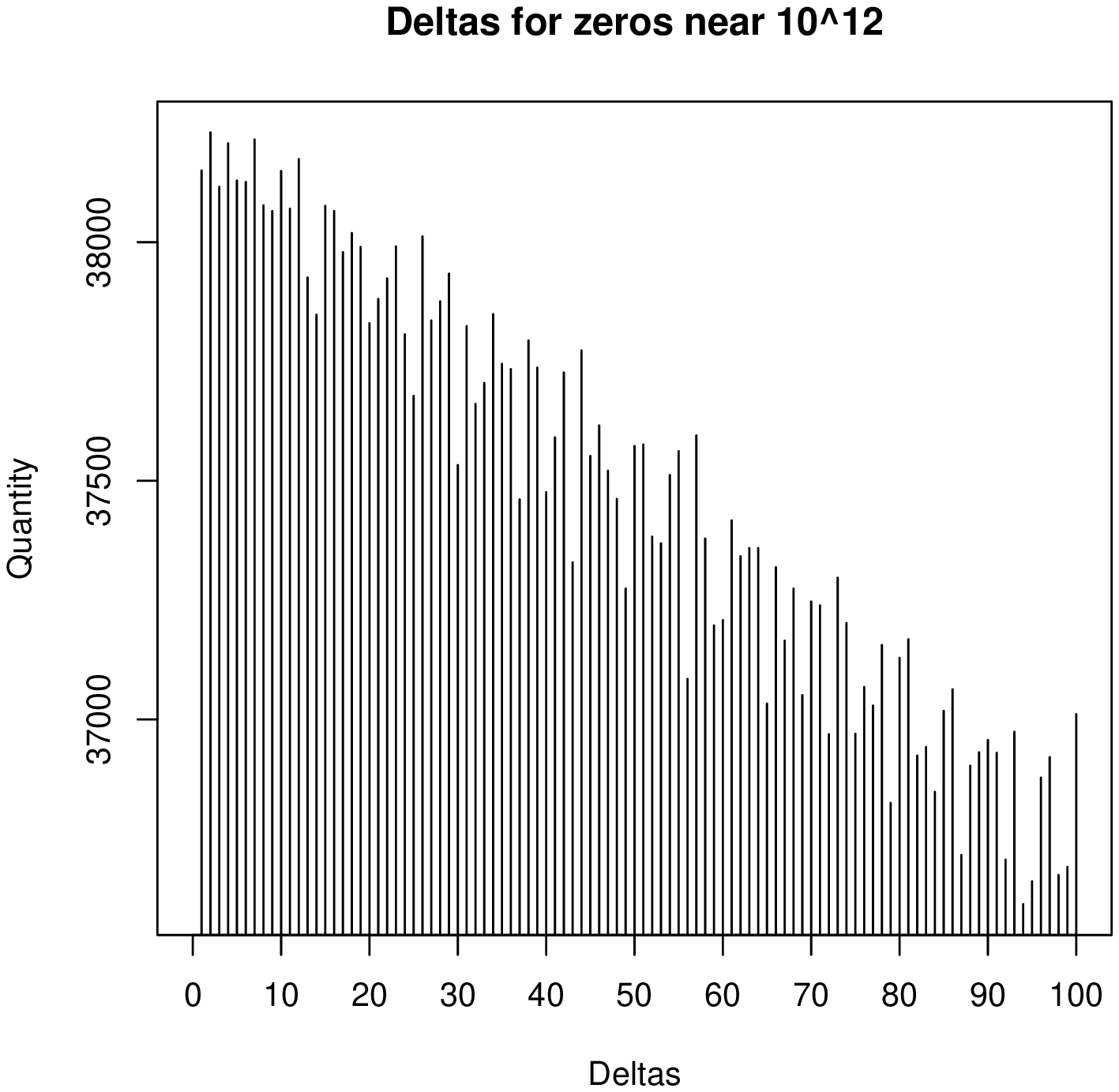}}    
\resizebox{6cm}{!}{\includegraphics{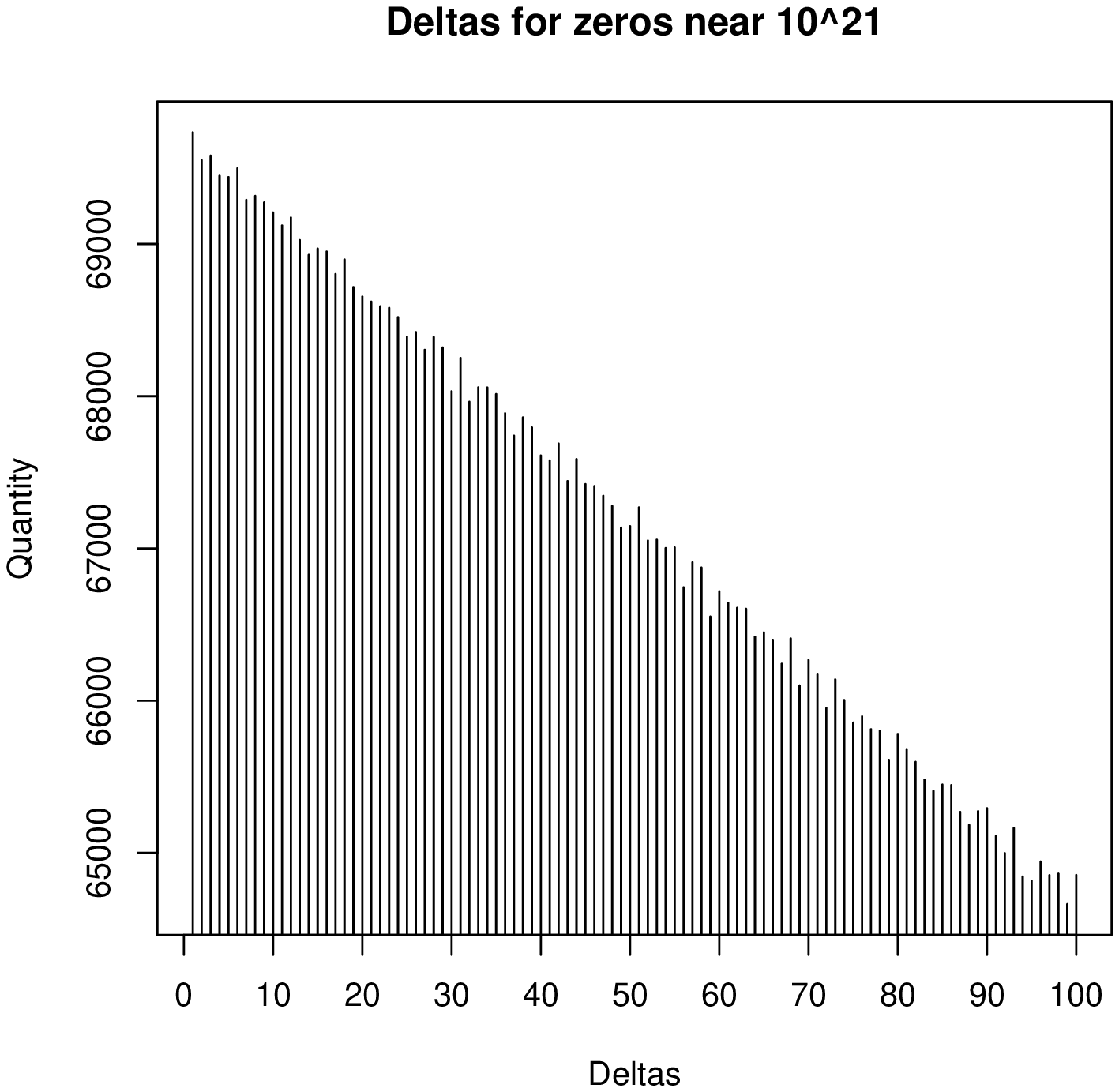}}
\end{center}

\centerline{Figures 11.a and 11.b.}


\begin{center}
  \resizebox{6cm}{!}{\includegraphics{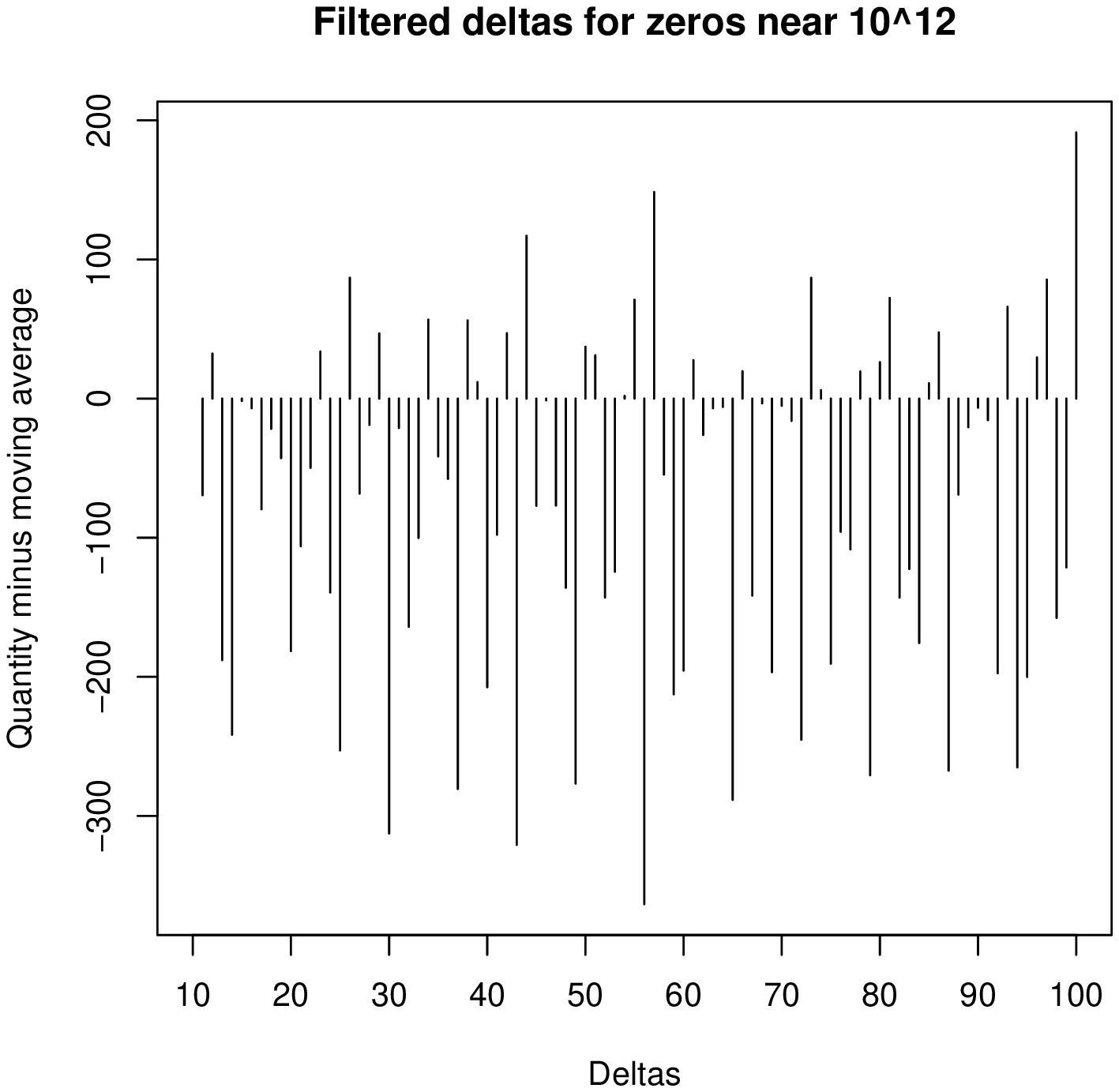}}    
\resizebox{6cm}{!}{\includegraphics{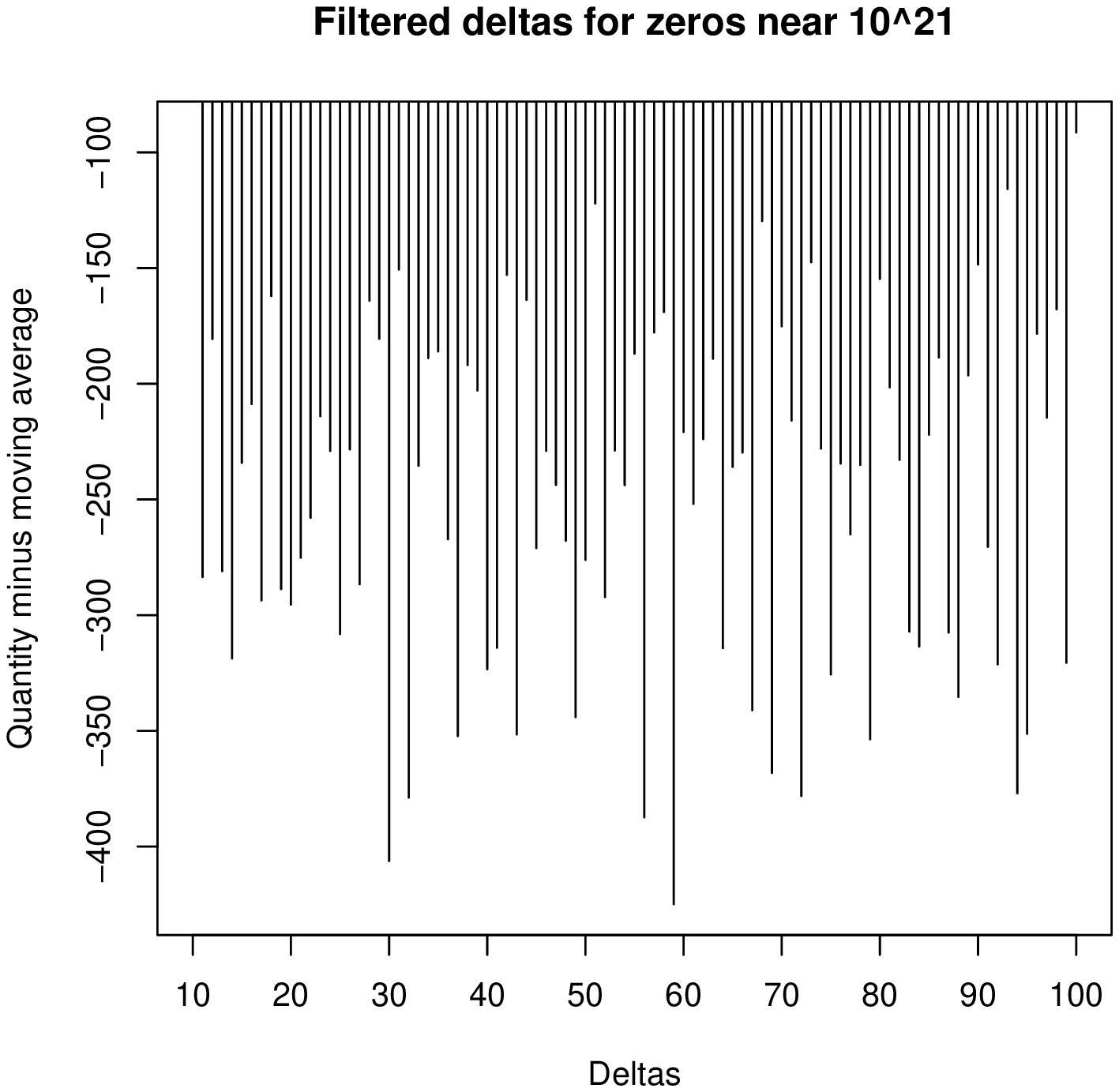}}
\end{center}

\centerline{Figures 12.a and 12.b.}


\begin{center}
  \resizebox{6cm}{!}{\includegraphics{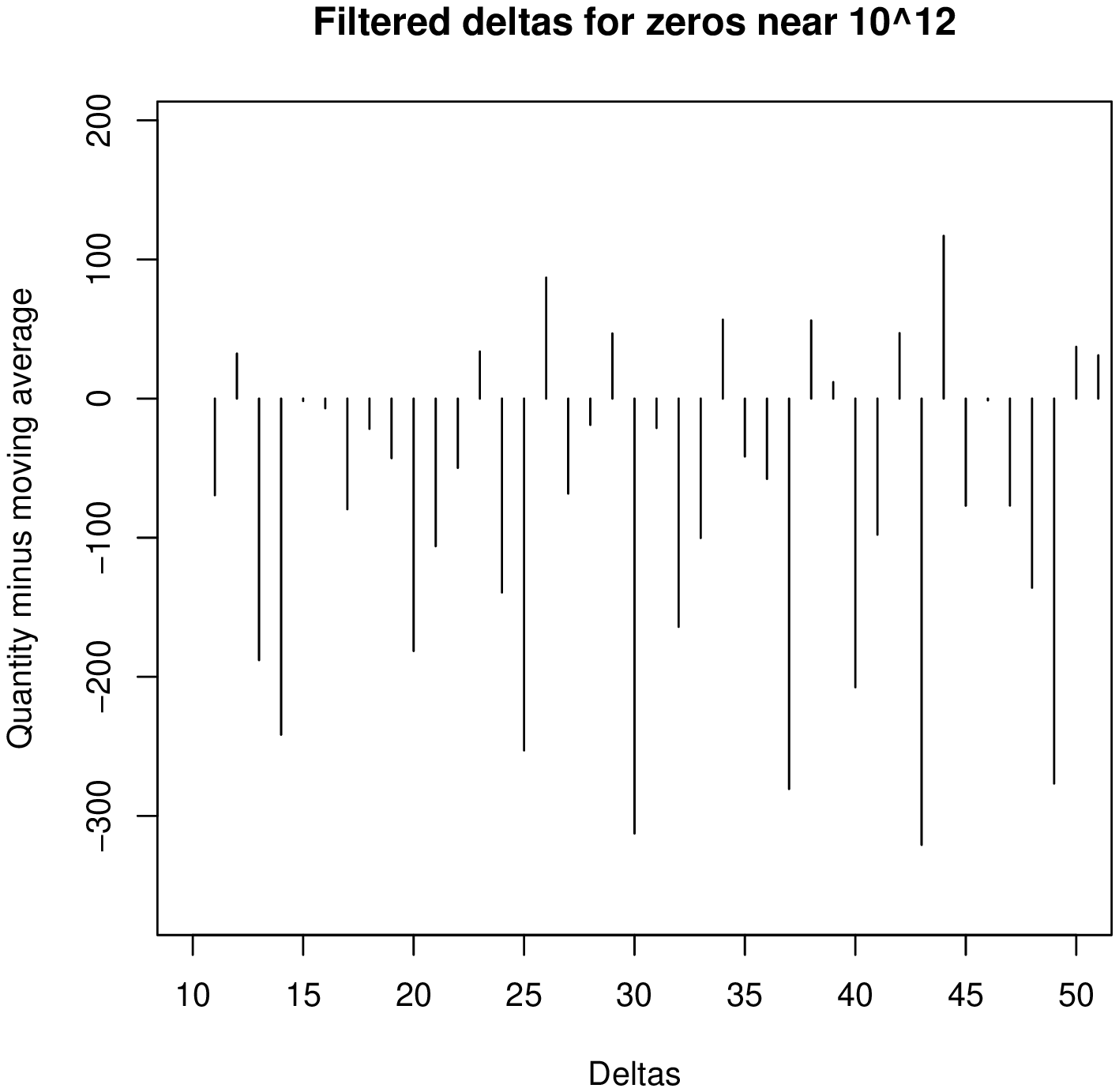}}    
\resizebox{6cm}{!}{\includegraphics{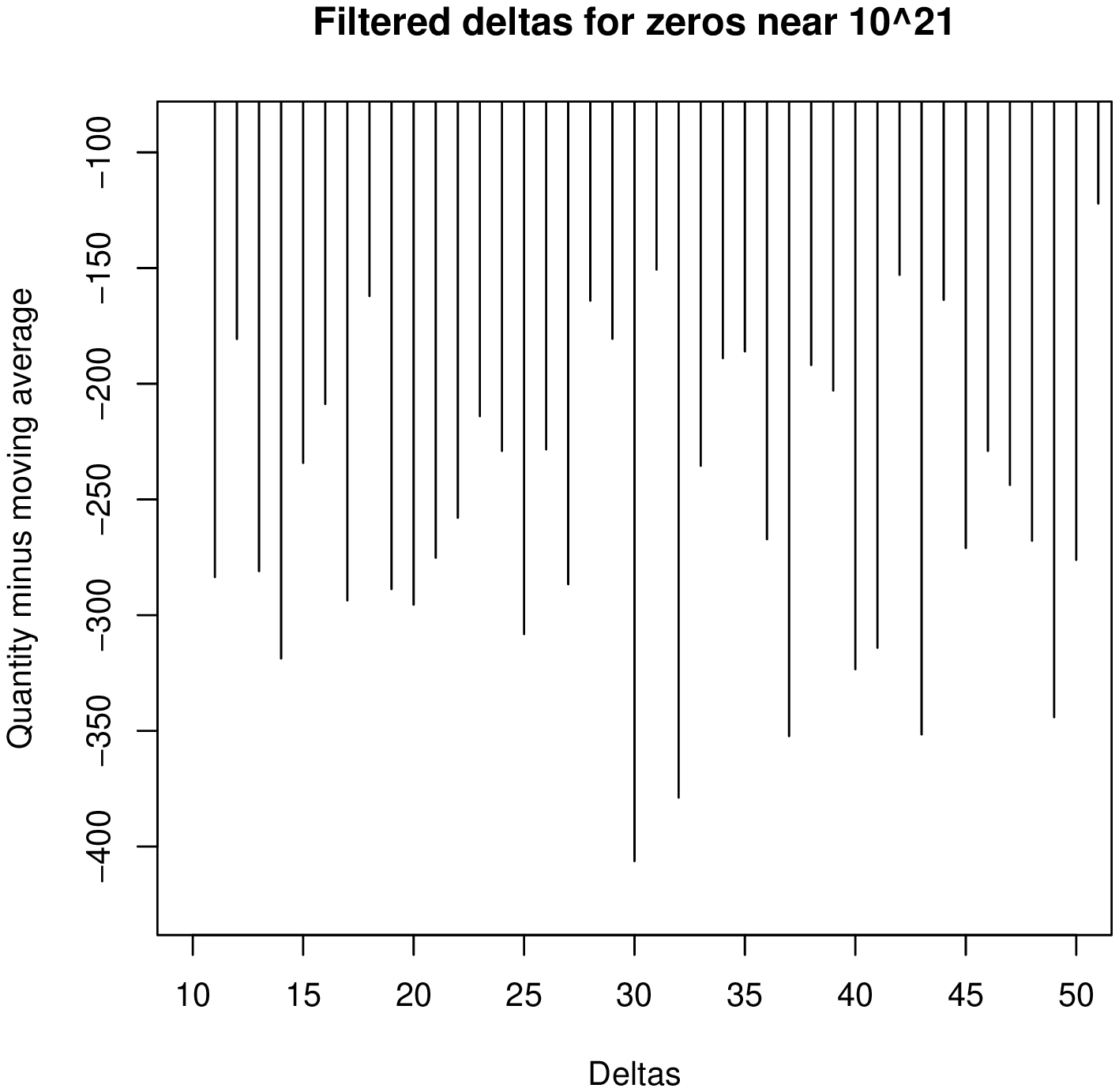}}
\end{center}

\centerline{Figures 13.a and 13.b.}


\section{Zeros of $L$-functions do know Riemann zeros.} \label{sec:L-functions_know}

In the survey article of B. Conrey on the Riemann Hypothesis
(\cite{Co}) we can read in the section entitled "The conspiracy of
$L$-functions",


{\it There is a growing body of evidence that there is a
conspiracy among $L$-functions (...) The first clue that zeta- and
L-functions even know about each other appears perhaps in works of
Deuring and Heilbronn (...)These results together (...) gave the
first indication of a connection between the zeros of $\zeta(s)$
and those of $L(s, \chi_d)$. }


We confirm in this section that zeros of $L$-functions do know about all Riemann
zeros. Indeed we provide a simple algorithm that builds the sequence of
Riemann zeros from the sequence of zeros of any Dirichlet
$L$-functions. Our first example is for the simplest non-trivial
$L$-function: We show how to recover Riemann zeros from the zeros of
$L_{\chi_3}$, where $\chi_3$ is the only character of conductor
$3$.

We perform the  statistics for the deltas of the zeros of
$L_{\chi_3}$ as done in section \ref{sec:self} for Riemann zeros. This time
we observe that the deficit values for the deltas of zeros of
$L_{\chi_3}$ is located precisely at Riemann zeros. As in
section \ref{sec:self} we perform one statistic with $100\ 000$ zeros of
$L_{\chi_3 }$ and precision $0.1$ for the deltas, and another,
more intensive, with $5$ million zeros of $L_{\chi_3}$. Figures 14
show the histogram of deltas for both statistics. Figures 15 show
the details in the interval $[10,30]$, figures 16 for $[30,50]$,
and figures 17 for $[80,100]$. We observe in figures 18 the
deficit of deltas near $0$ verifying Montgomery's prediction.

The similarity of these figures with those in section \ref{sec:self} is
clear. Recall though that they are generated from a very different
set of data.


\begin{center}
  \resizebox{6cm}{!}{\includegraphics{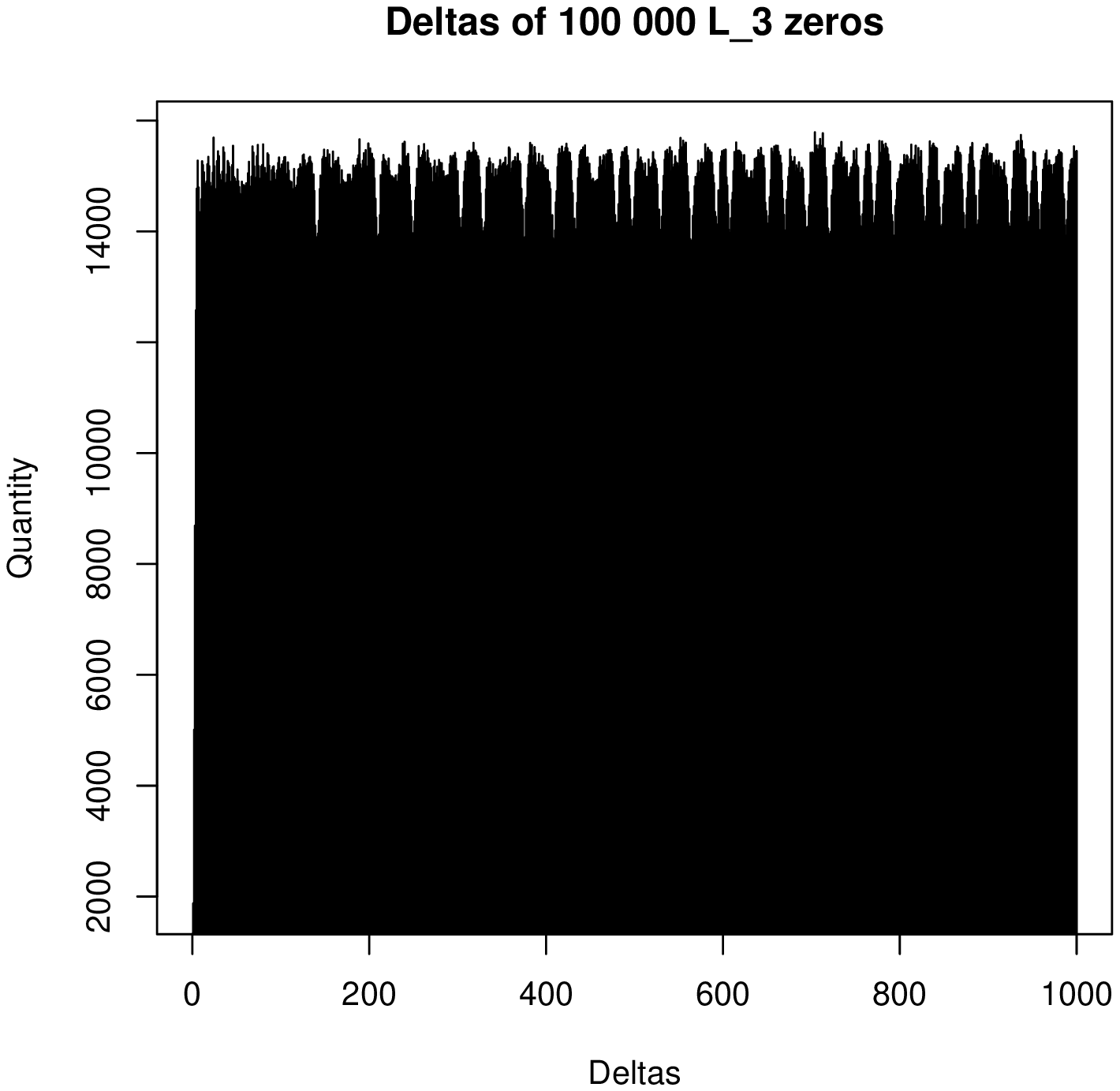}}    
\resizebox{6cm}{!}{\includegraphics{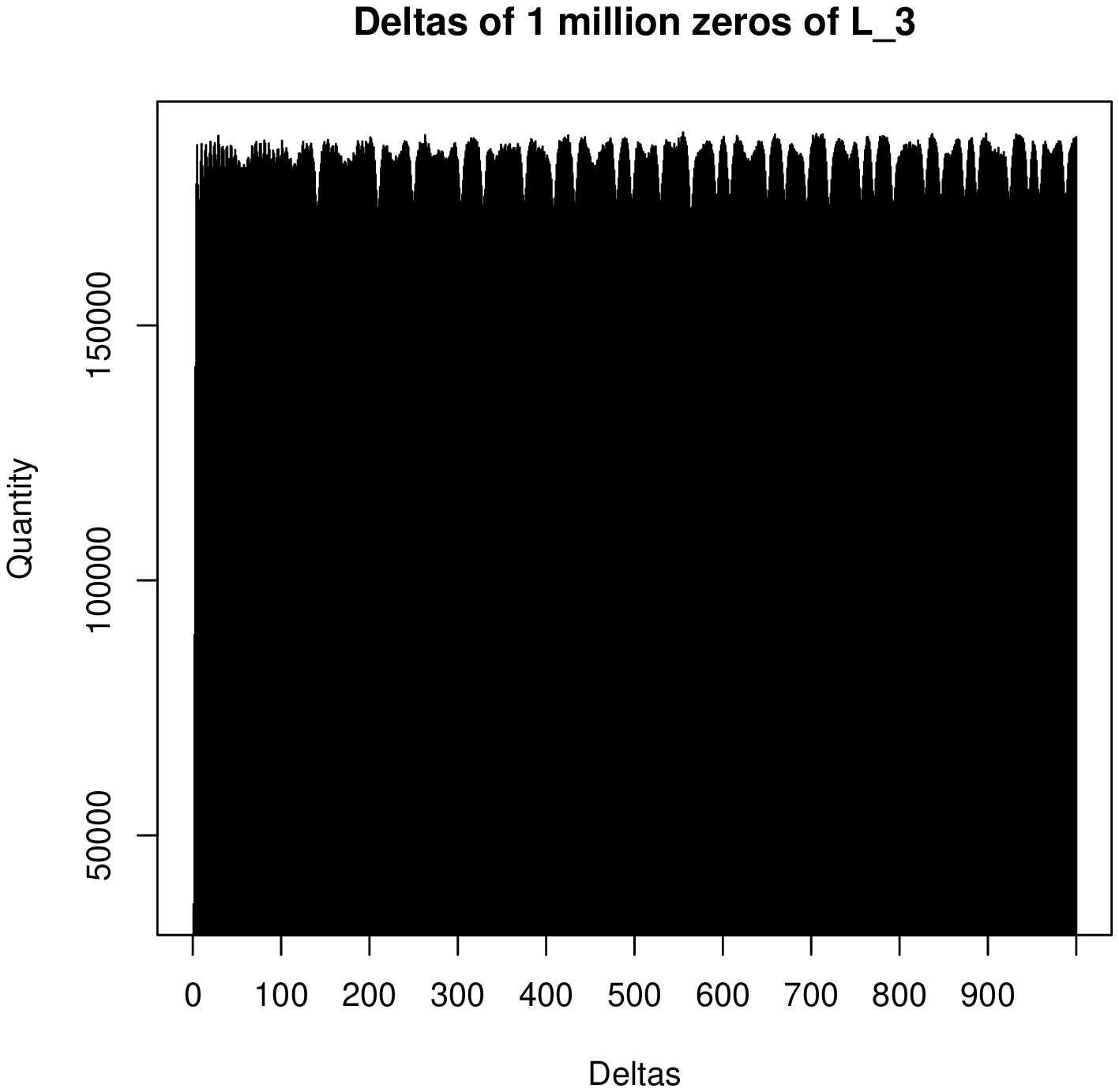}}
\end{center}

\centerline{Figures 14.a and 14.b.}


\begin{center}
  \resizebox{6cm}{!}{\includegraphics{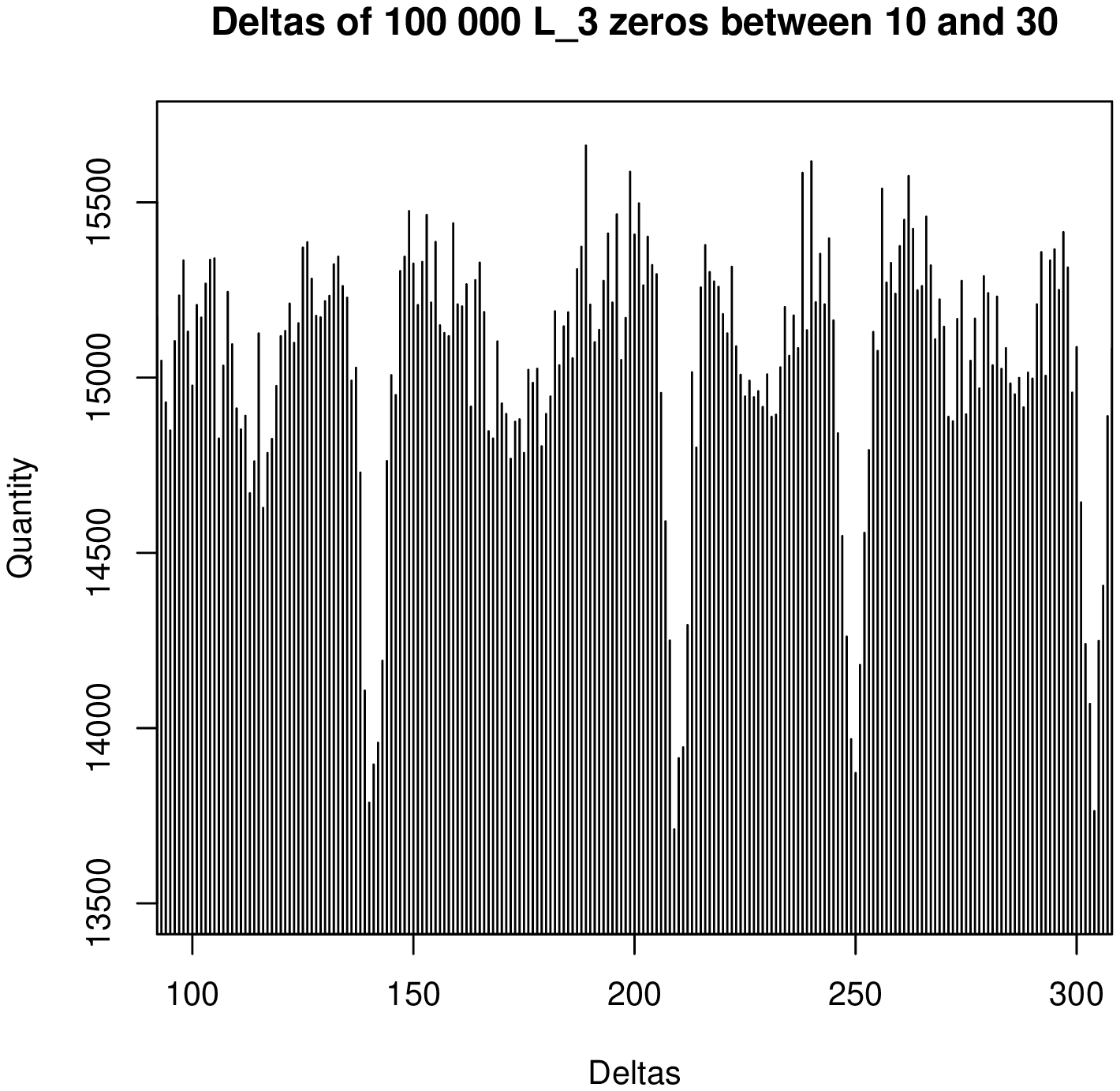}}    
\resizebox{6cm}{!}{\includegraphics{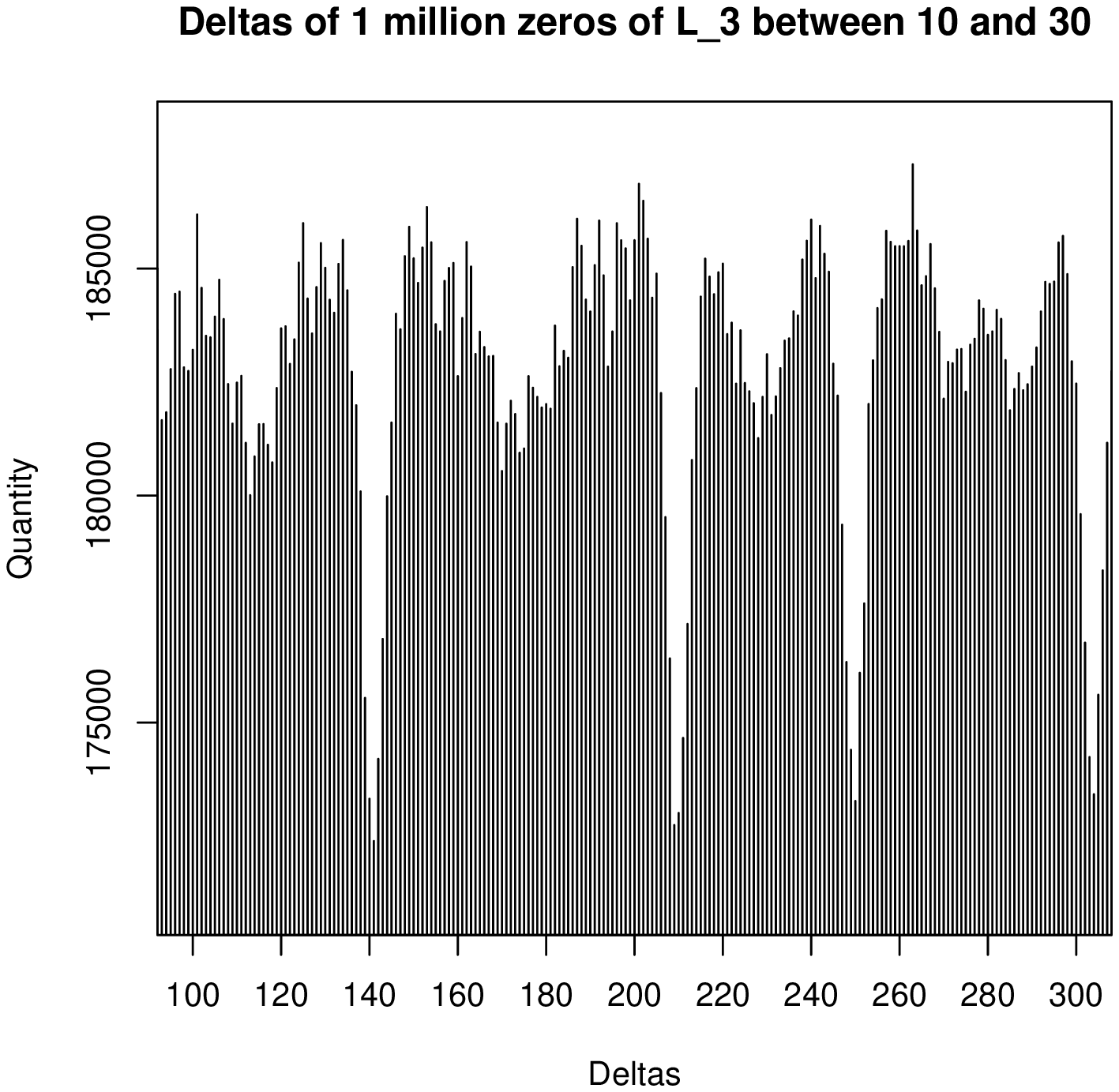}}
\end{center}

\centerline{Figures 15.a and 15.b.}


\begin{center}
  \resizebox{6cm}{!}{\includegraphics{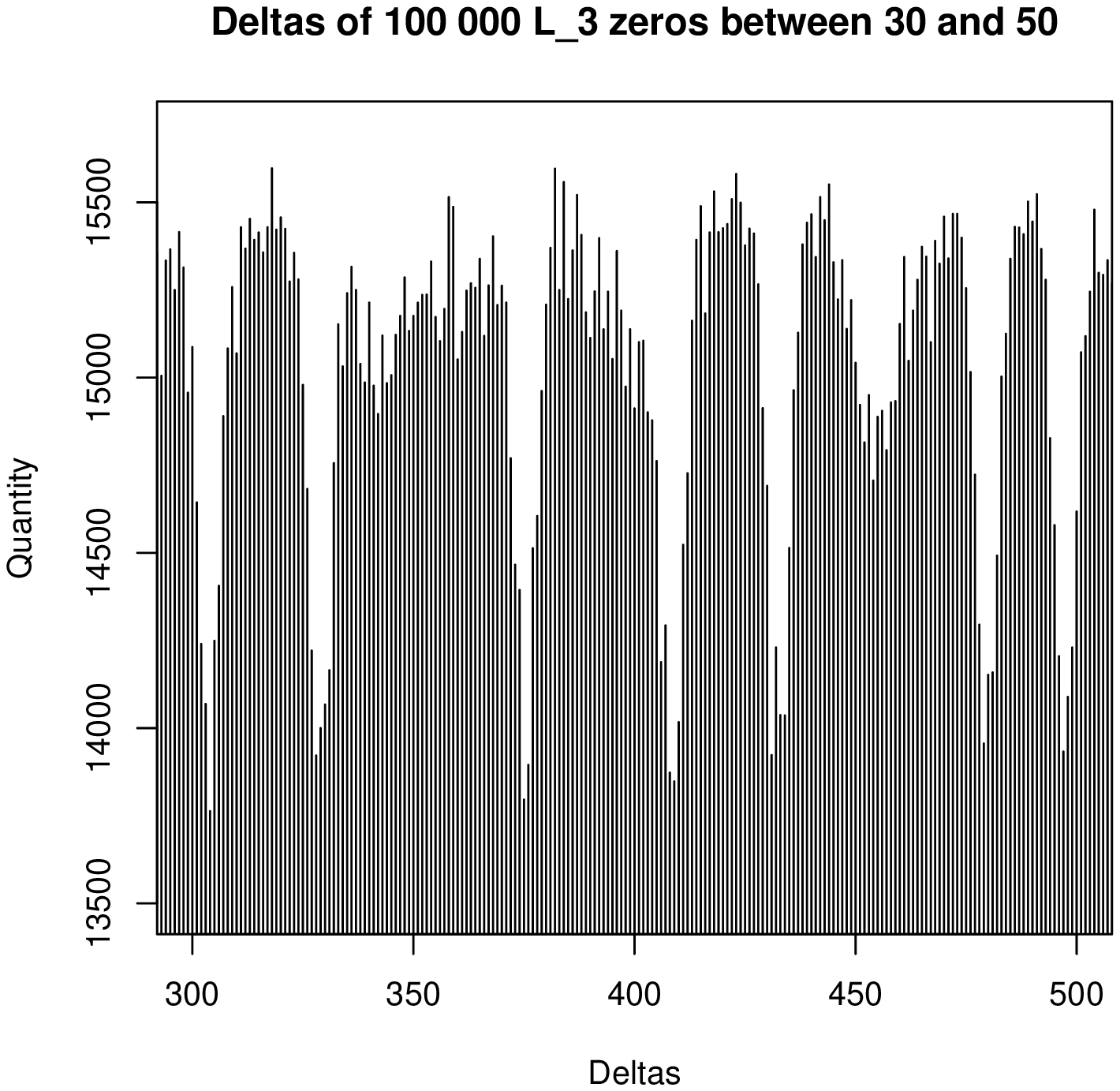}}    
\resizebox{6cm}{!}{\includegraphics{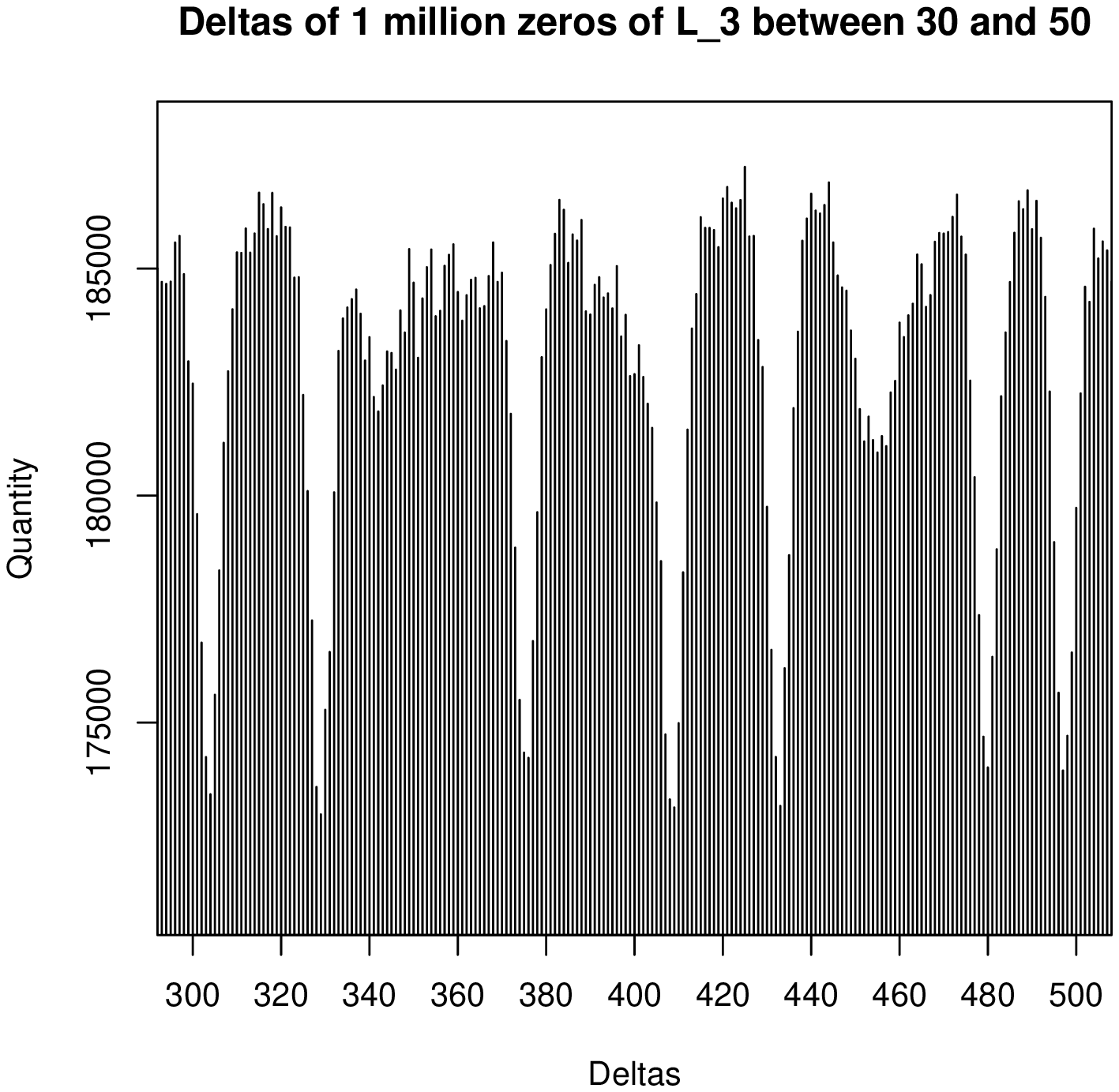}}
\end{center}

\centerline{Figures 16.a and 16.b.}


\begin{center}
  \resizebox{6cm}{!}{\includegraphics{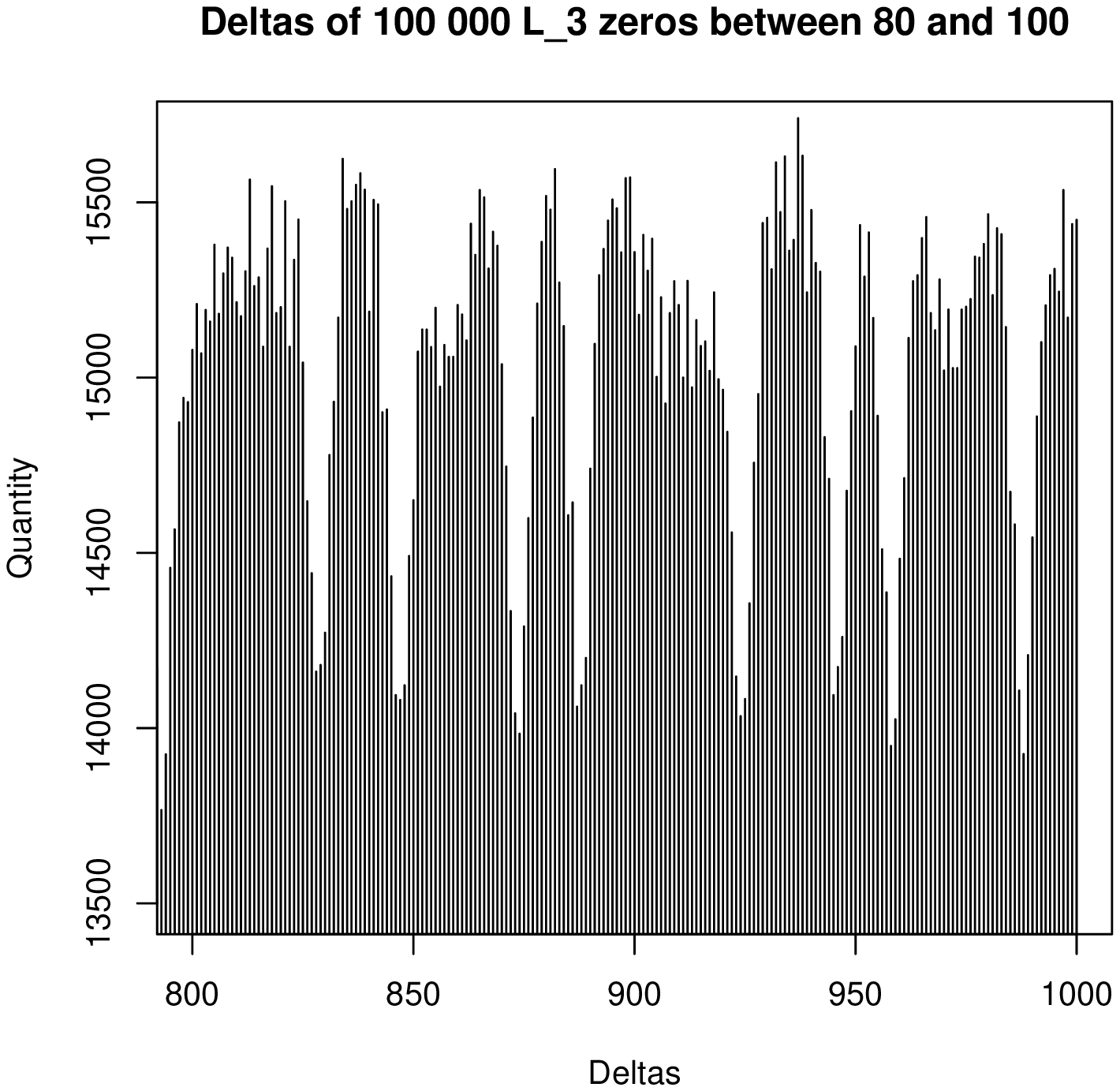}}    
\resizebox{6cm}{!}{\includegraphics{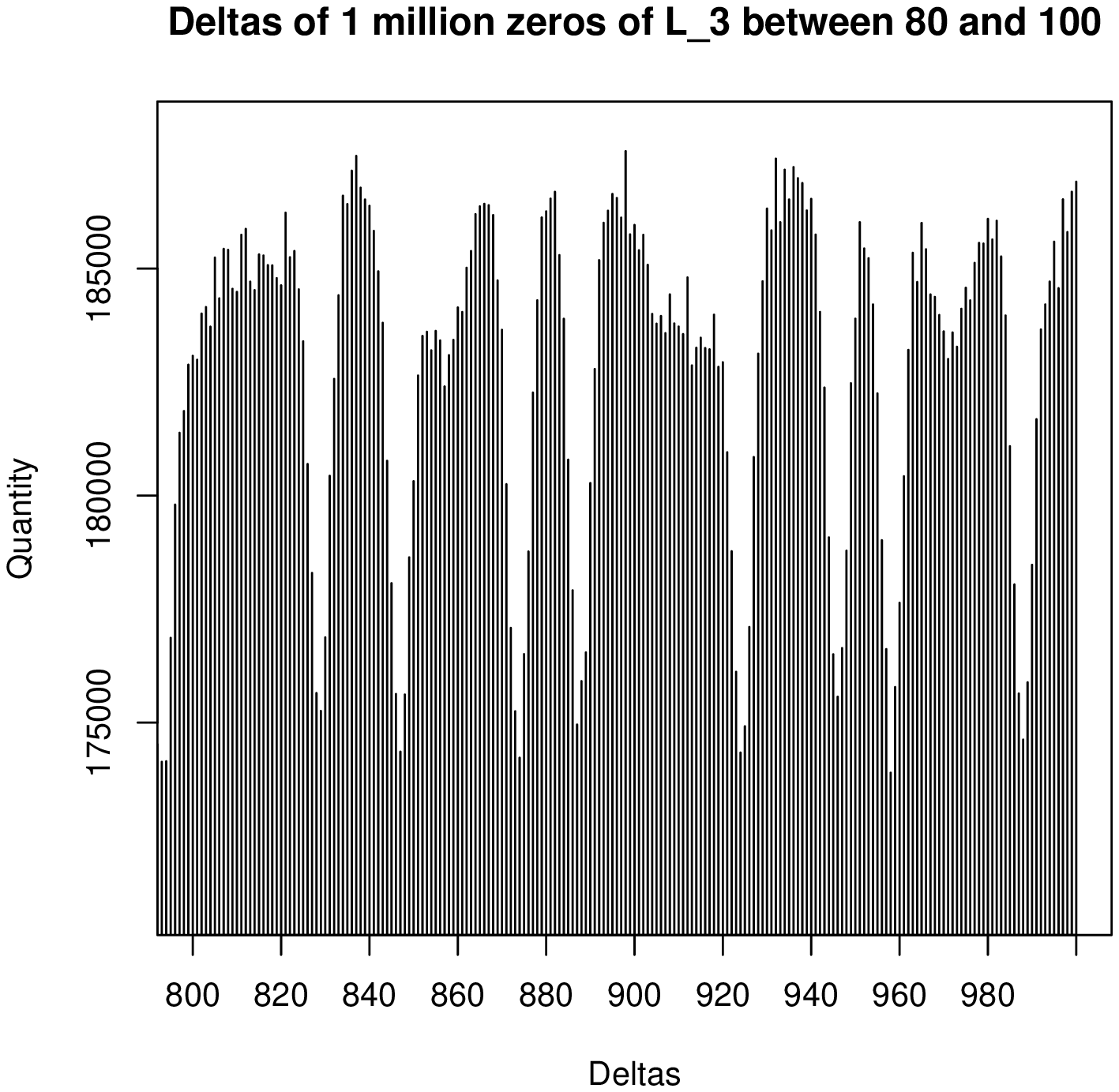}}
\end{center}

\centerline{Figures 17.a and 17.b.}


\begin{center}
  \resizebox{6cm}{!}{\includegraphics{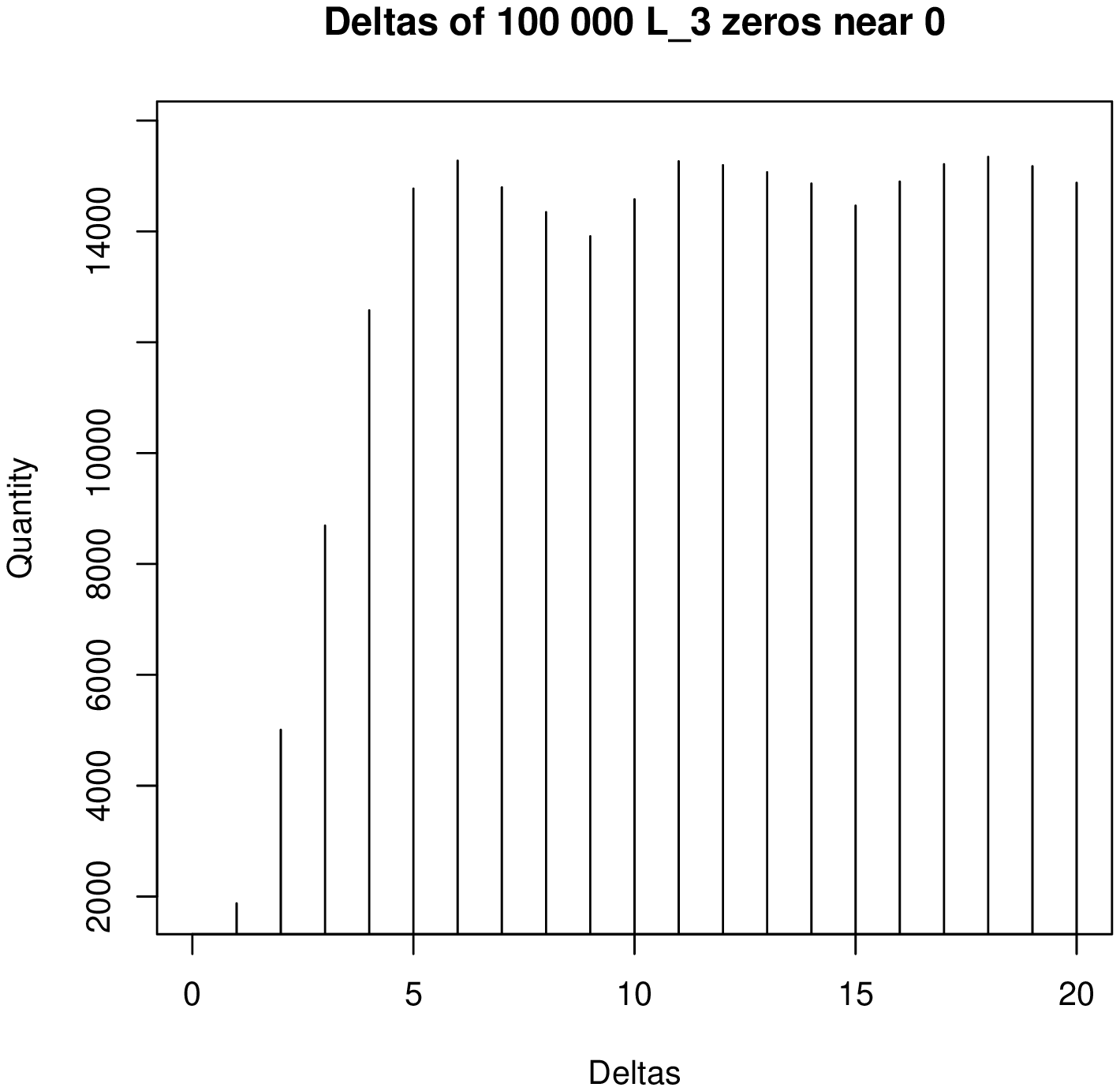}}    
\resizebox{6cm}{!}{\includegraphics{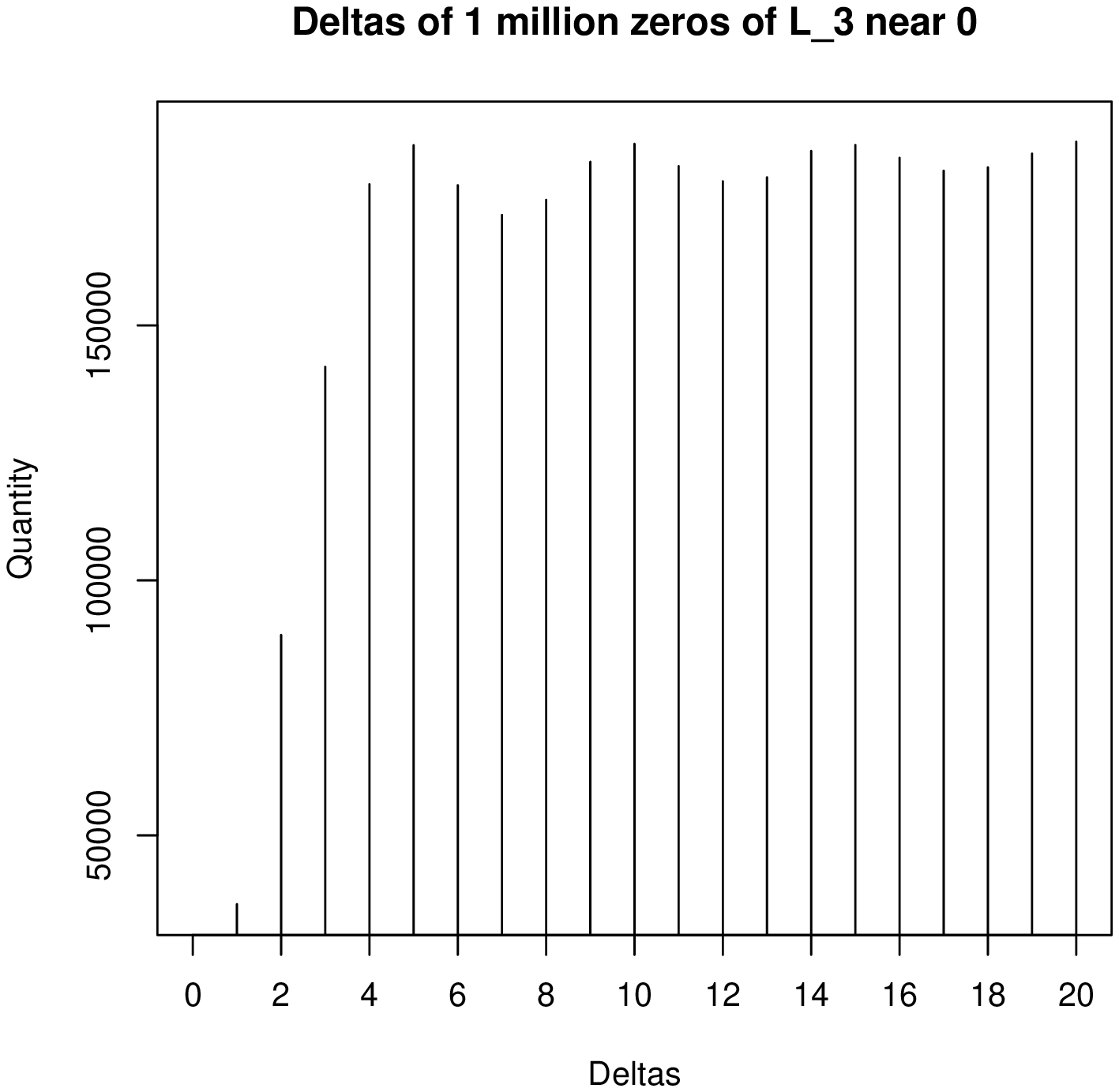}}
\end{center}

\centerline{Figures 18.a and 18.b.}

\bigskip

Now we perform the same statistics for the zeros of other
Dirichlet $L$-functions $L_{\chi }$. We perform the statistics for
the deltas of $1$ million zeros, for deltas in $[0,100]$, and with
precision $0.1$. This time we consider a real and a complex
non-real character. For a complex non-real character, the
associated $L$-function is not real analytic, and the zeros are no
longer symmetric with respect to the real axes. Therefore we
compute the deltas of those with positive imaginary part and the
deltas for those with negative imaginary part, and we compute the
cumulative result. Since the sequence of zeros is not symmetric
with respect to $0$, we take the first million zeros in the
following sense: We order the zeros by absolute value and we
consider the first million of them for the statistics.

The first statistics is for $\chi=\chi_{4}$, the only primitive
character of conductor $4$. The character $\chi_4$ is real and the
associated Dirichlet function real-analytic. The second statistics is for
$\chi=\chi_{7,3}$, one of the primitive complex characters of
conductor $7$. Figures 19 show the histograms of deltas in
$[10,30]$. Figures 20 show the histograms of deltas in $[30,50]$.
Again we find that the deficit locations coincide with
Riemann zeros.


\begin{center}
  \resizebox{6cm}{!}{\includegraphics{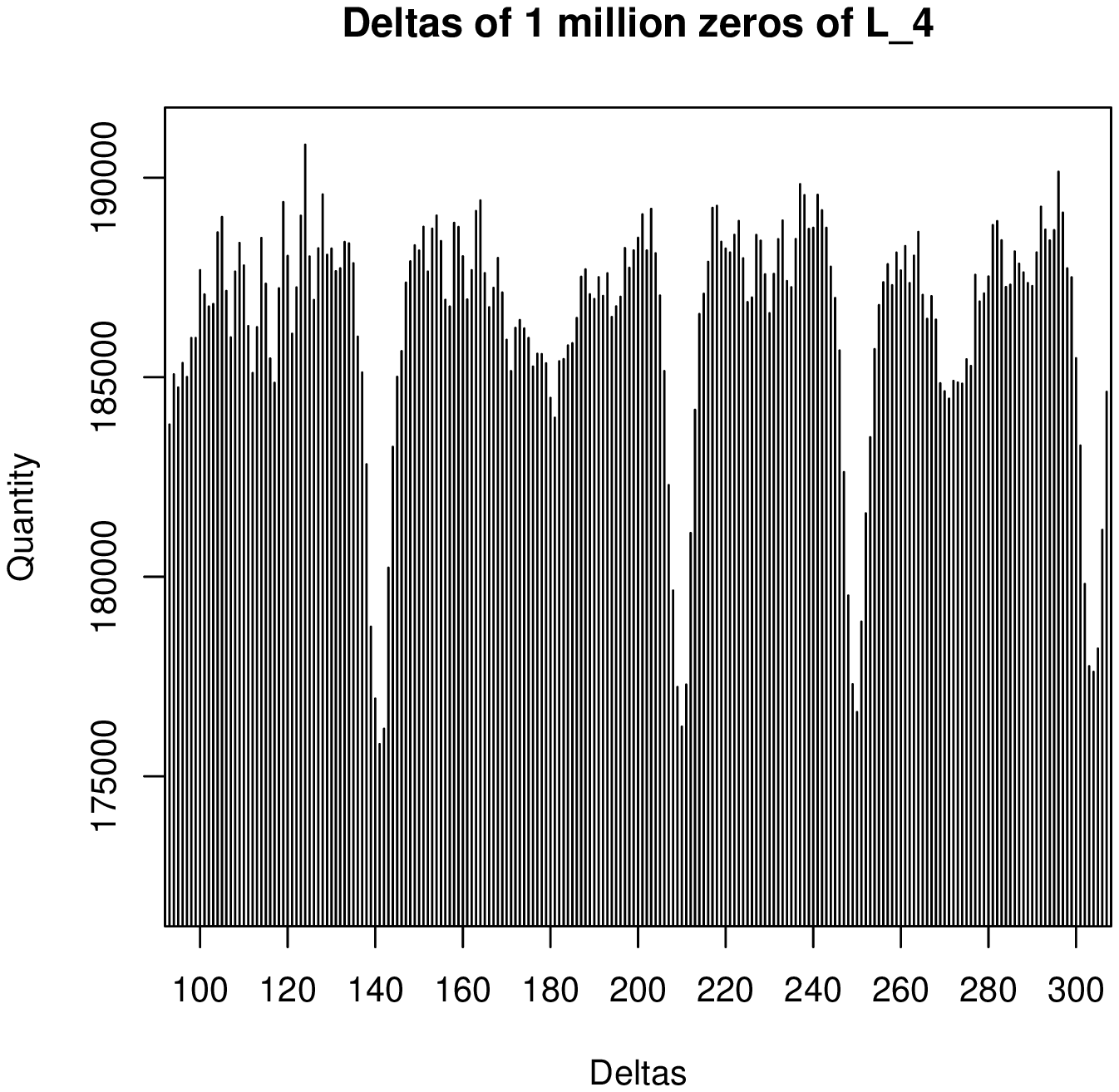}}    
\resizebox{6cm}{!}{\includegraphics{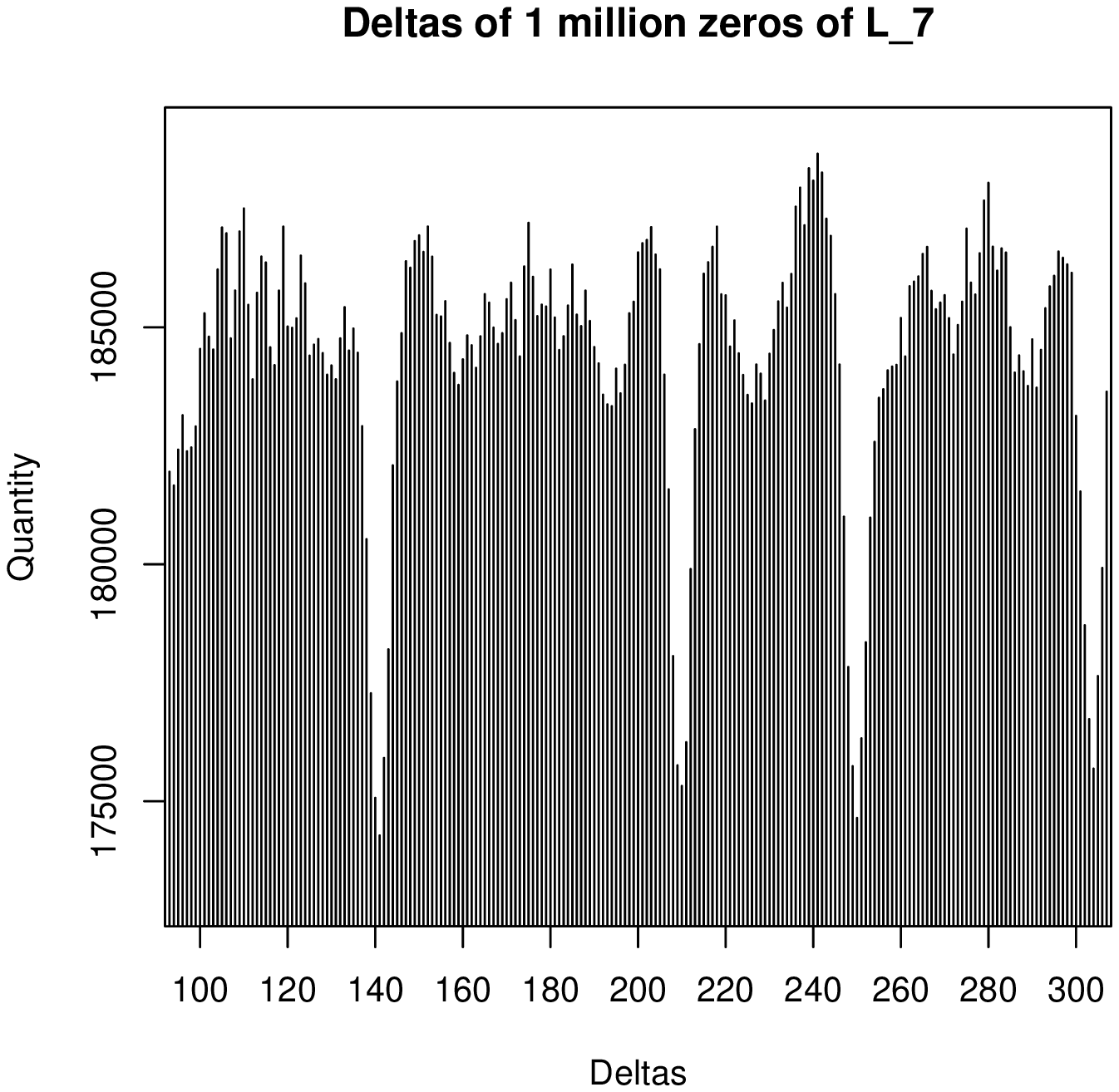}}
\end{center}

\centerline{Figures 19.a and 19.b.}


\begin{center}
  \resizebox{6cm}{!}{\includegraphics{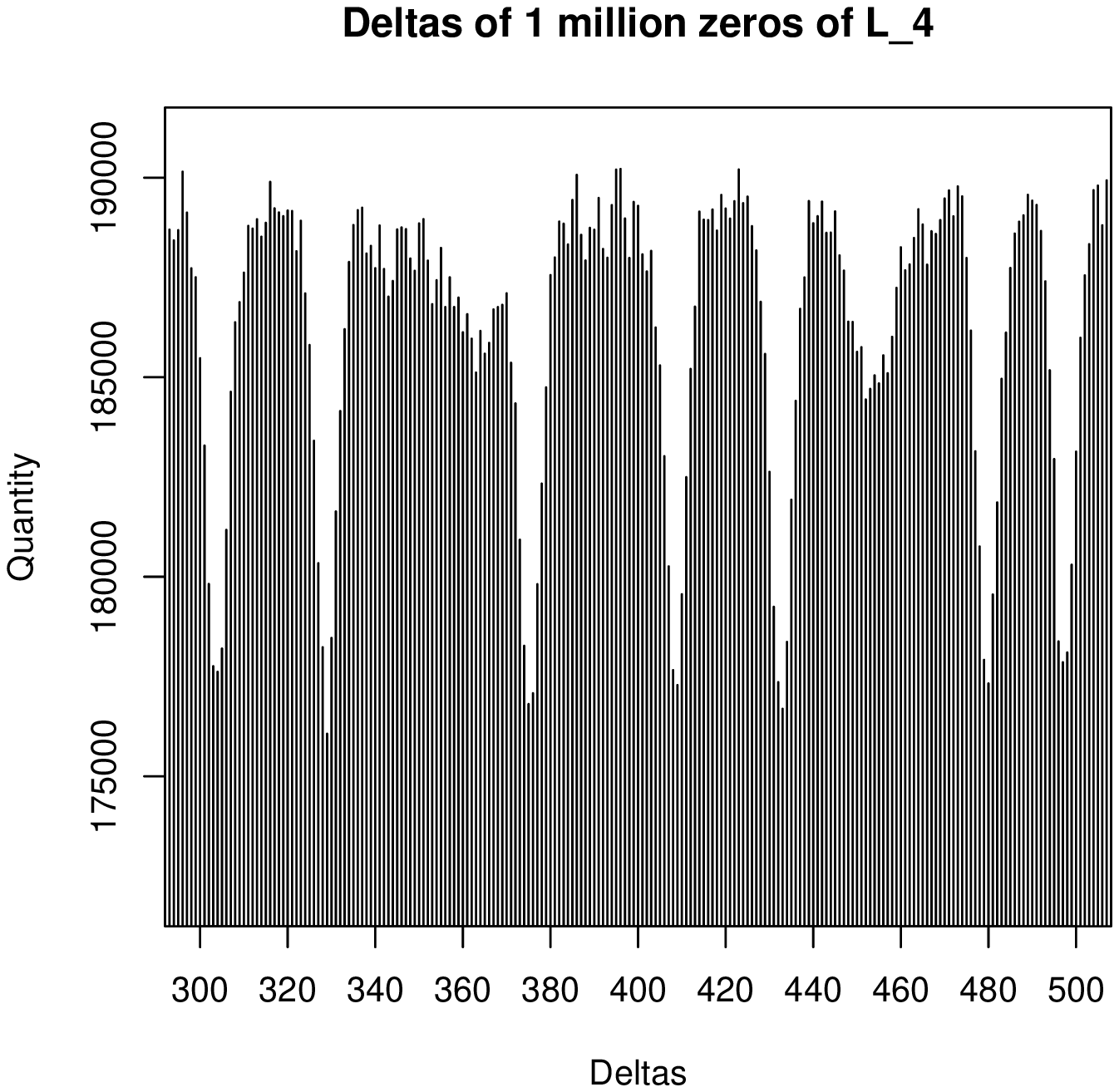}}    
\resizebox{6cm}{!}{\includegraphics{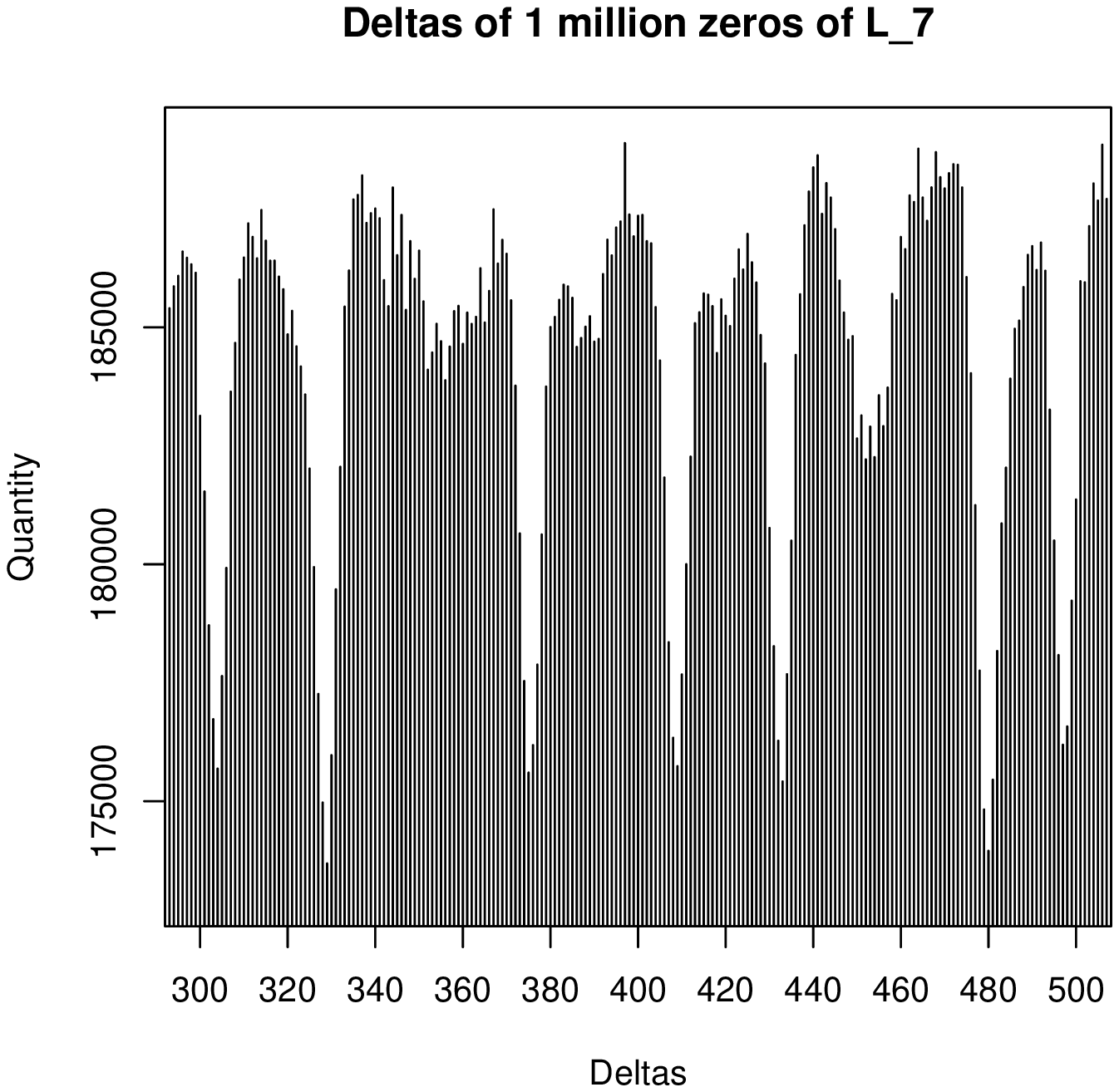}}
\end{center}

\centerline{Figures 20.a and 20.b.}

\bigskip

We analyze next the distribution of deltas near $0$. We plot the
histogram near $0$ of the deltas of $1$ million zeros with
precision $0.01$. We observe the predicted GUE pair correlation
distribution as pictures 21 show.


\begin{center}
  \resizebox{6cm}{!}{\includegraphics{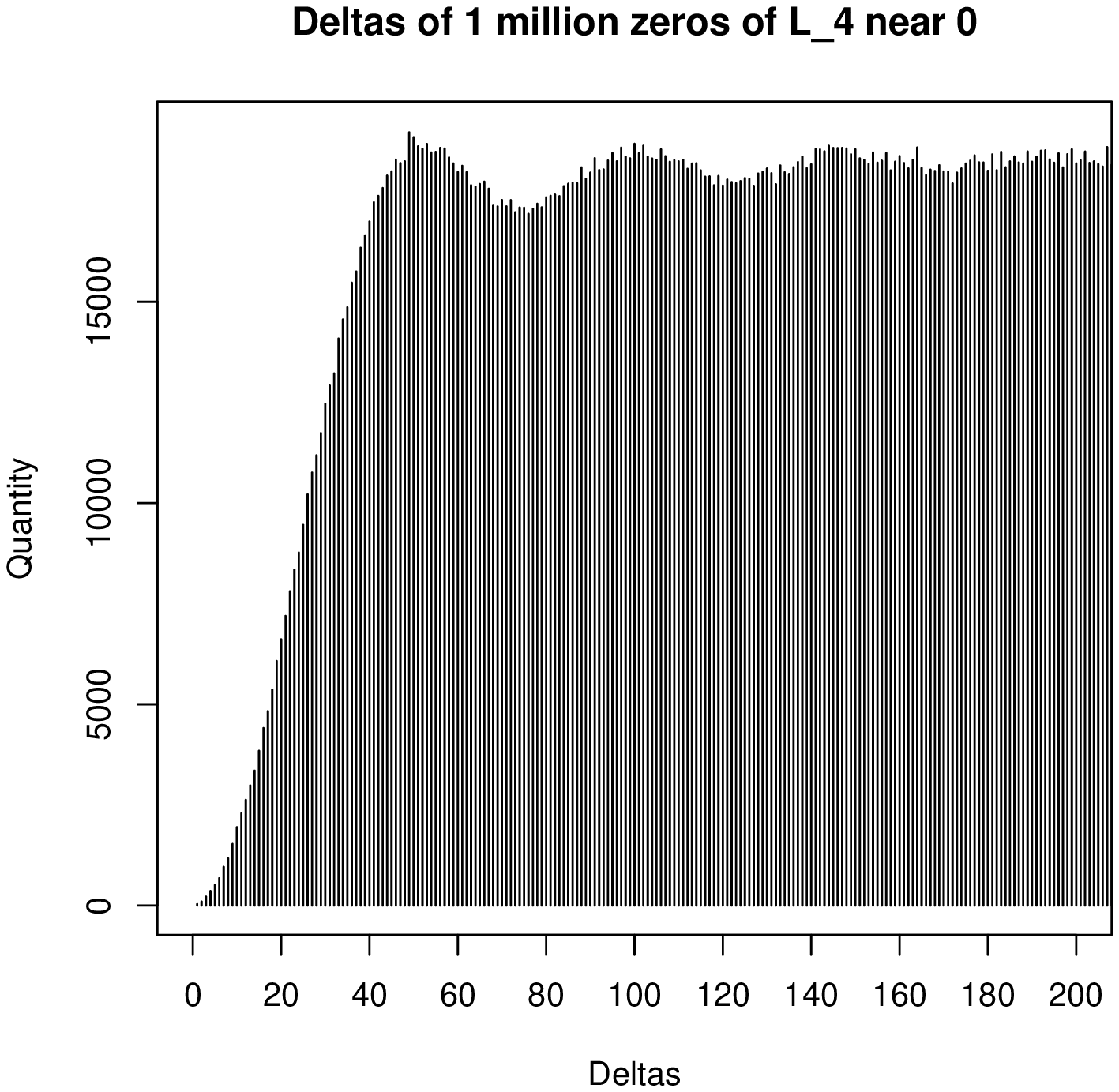}}    
\resizebox{6cm}{!}{\includegraphics{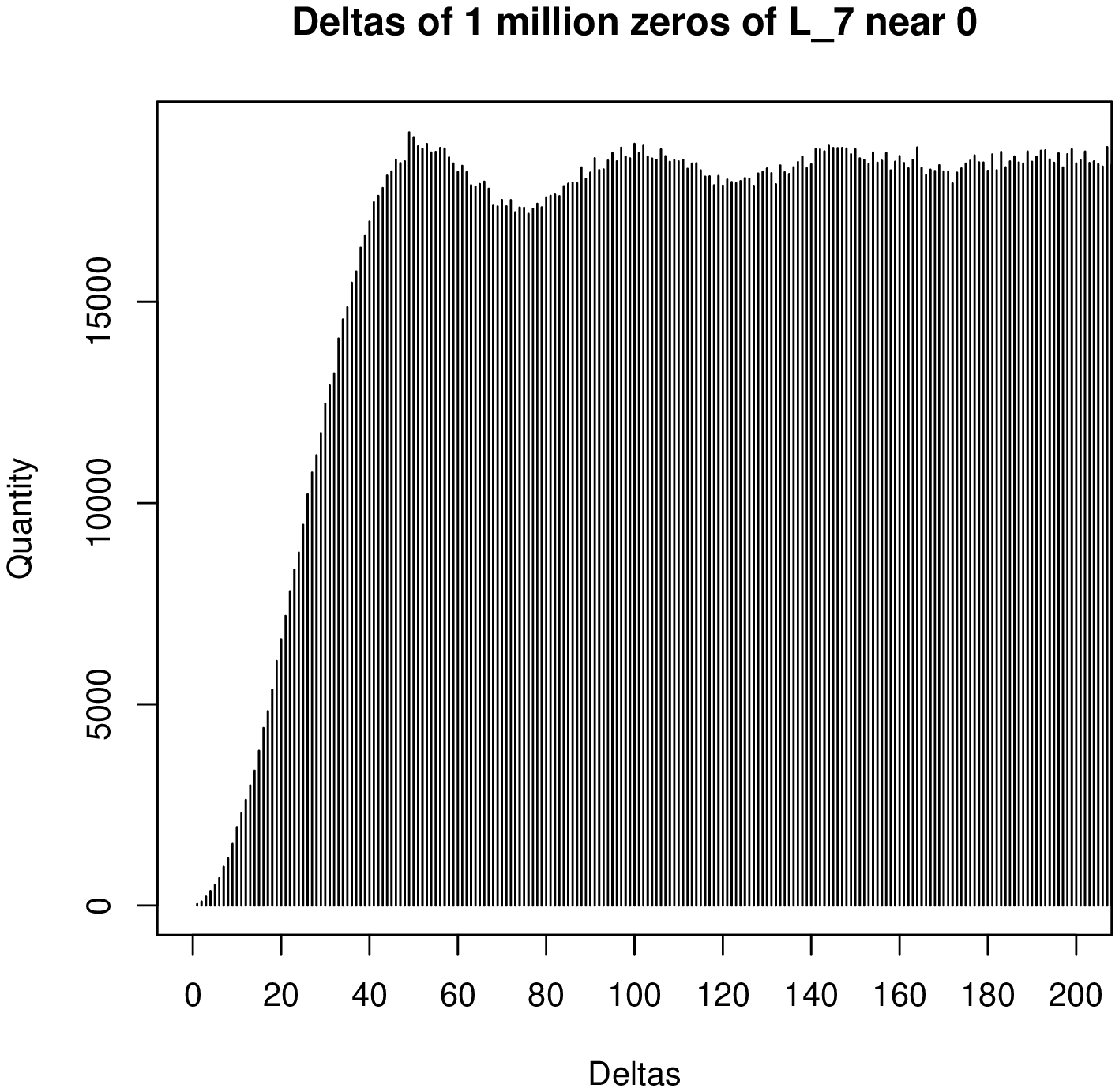}}
\end{center}

\centerline{Figures 21.a and 21.b.}

\bigskip

\textbf{E\~ne product explanation.}

\medskip

 We note that for a real character $\chi $,
 $$
 \bar L_{\chi }(s)={\overline {L_{\chi }(\bar s)}}=L_{\chi }(s)\ .
 $$
 The e\~ne product explanation of the first numerical result for $\chi_3$ is based on the
 following computation (where we denote by $\chi_0$ the principal character modulo $3$)

\begin{align*}
L_{\chi_3} \ \bar \star \ \bar L_{\chi_3}&= L_{\chi_3} \ \bar \star \  L_{\chi_3} \\ 
&= \left ( \prod_p (1-\chi_3(p)
p^{-s} )^{-1} \right ) \bar \star \left ( \prod_q (1-\chi_3(q)
q^{-s} )^{-1} \right ) \\ 
&=\prod_p (1-\chi_3(p) p^{-s} )^{-1} \
\bar \star  \ (1-\chi_3(p) p^{-s} )^{-1} \\ 
&=\prod_p (1-\chi_3(p)^2 p^{-1/2}  p^{-s} ) \\ &=\prod_p (1-\chi_0(p)
p^{-1/2-s} ) \\ 
&=(1-3^{-1/2-s})^{-1} \prod_p (1- p^{-1/2-s} ) \\ &= (1-3^{-1/2-s})^{-1} \zeta(s+1/2)^{-1} \ . 
\end{align*}

In general, for an arbitrary character $\chi$ modulo $n$, we
recognize the distribution of the deltas of the zeros of
$L_{\chi}$ in the result of the e\~ne product of $L_{\chi}$ with
$\bar L_{\chi}$. We have
$$
{\bar L_{\chi }}=L_{\overline {\chi }} \ .
$$
Also observe that
$$
\chi .{\overline  {\chi }} =| \chi |^2=\chi_0 \ ,
$$
where $\chi_0$ is the principal character modulo $n$.

Therefore we have

\begin{align*}
L_{\chi} \ \bar \star \ \bar L_{\chi} &= L_{\chi } \
\bar \star \  L_{\bar \chi} \\ &=\left ( \prod_{p \not| n }
(1-\chi (p) p^{-s} )^{-1} \right ) \ \bar \star \ \left ( \prod_{q
\not| n} (1-{\overline {\chi }}  (q) q^{-s} )^{-1} \right ) \\
&=\prod_p (1-\chi (p) p^{-s} )^{-1} \ \bar \star  \ (1-{\overline
{\chi }}(p) p^{-s} )^{-1} \\ &=\prod_p (1-|\chi(p)|^2 p^{-1/2}
p^{-s} ) \\ &=\prod_p (1-\chi_0(p)  p^{-1/2-s} ) \\
&=\prod_{p|n} (1-p^{-1/2-s})^{-1} . \prod_p (1-  p^{-1/2-s} ) \\
&= \zeta(s+1/2)^{-1} \prod_{p|n} (1-p^{-1/2-s})^{-1} \ . 
\end{align*}

Observe that the zeros of each Euler factor
$$
f_p(s)=(1-p^{-1/2-s})^{-1}
$$
are for $k\in \ZZ$,
$$
s_k=-{1\over 2} +i {2\pi \over \log p} k \ .
$$
According to the explanation with the e\~ne product, we should
observe a deficit of deltas (with a lower order amplitude) near
locations multiples of the fundamental harmonic
$$
{2\pi \over \log p} k
$$
for $k\in \ZZ$.

For $p=2$, $p=3$ and $p=7$ we have,

\begin{align*}
{2\pi \over \log 2} &=9.0647\ldots \\ {2\pi \over \log
3} &=5.7192\ldots \\ {2\pi \over \log 7} &=3.2289\ldots\ldots
\end{align*}

A good eye can spot a trace of these deficit locations in the
figures 16, 17, 19 and 20. In particular comparing these figures
with figures 3 and 4.


\bigskip

\textbf{Script.}

\bigskip

We provide the script for computing the deltas of the zeros of a
complex non-real L-function since it is slightly different from
the previous ones. Here we feed the program by reading into
Rubinstein's table "zeros-0007-2000000" which contains the first
zeros for each of the three primitive characters of conductor $7$.
The zeros for the complex character that we are considering are
those after row $2\ 000\ 000$.

\medskip

{\tt

zeros<-read.table("zeros-0007-2000000",skip=2000000,nrows=1000000)

z<-zeros[,3]

z.plus<-z[z>0]

z.minus<--z[z<0]

x=rep(0,10000)

zeta1<-z.plus

zeta2<-z.plus

N=length(z.plus)

k=0

for (i in 1:N)

$\{$

while( ((zeta1[i]-zeta2[i+k])<100.01) \& (k+i >1) )

$\{$

k<-k-1

$\}$

k=k+1

j=k

while ( (zeta1[i]-zeta2[i+j]>0) \& (zeta1[i]-zeta2[i+j]<5.01) )

$\{$

d=100*(zeta1[i]-zeta2[i+j])

x[as.integer(d)]=x[as.integer(d)]+1

j=j+1

$\}$

$\}$

zeta1=numeric()

zeta2=numeric()

zeta1<-z.minus

zeta2<-z.minus

N=length(z.minus)

k=0

for (i in 1:N)

$\{$

while( ((zeta1[i]-zeta2[i+k])<100.01) \& (k+i >1) )

$\{$

k<-k-1

$\}$

k=k+1

j=k

while ( (zeta1[i]-zeta2[i+j]>0) \& (zeta1[i]-zeta2[i+j]<5.01) )

$\{$

d=100*(zeta[i]-zeta3[i+j])

x[as.integer(d)]=x[as.integer(d)]+1

j=j+1

$\}$

$\}$

}

\section{Zeros of $L$-functions replicate from their mating with Riemann zeros.} \label{sec:L-functions_replicate}

We present in this section and the next one a new type of statistics.
We do study the statistics of differences of zeros of
an $L$-function $L_{\chi_1}$ with the zeros of another $L$-function $L_{\chi_2}$.
We name this operation the "mating" of zeros of $L_{\chi_1}$ and
$L_{\chi_2}$. As predicted by the e\~ne product theory, it appears that the sequence of Riemann zeros plays
the role of the unit for this mating operation. More precisely,
the statistics of this section verify that the mating of Riemann
zeros with those of another $L$-function $L$ yield as deficit values
the zeros of $L$ itself.

We perform the statistics mating the Riemann zeros with the zeros
of $L_{\chi_3}$ where $\chi_3$ is as before the only primitive
character of conductor $3$. The function $L_{\chi_3}$ is real
analytic and its zeros are symmetric with respect to the real
axes. We consider only the non-real (i.e. non-trivial) zeros with
positive imaginary part. We denote by $(1/2+i \g^{(3)}_i)_{i\geq
1}$, or simply $(\g^{(3)}_i)_{i\geq 1}$, the zeros of
$L_{\chi_3}$, with $i\mapsto \g^{(3)}_i$ increasing.. The first
$18$ ones, less than $51$, are the following

\begin{align*}
\g^{(3)}_1 &=8.039737156\ldots \\ \g^{(3)}_2
&=11.24920621\ldots \\  \g^{(3)}_3 &=15.70461918\ldots \\
\g^{(3)}_4 &=18.2619975\ldots \\ \g^{(3)}_5 &=20.45577081\ldots
\\ \g^{(3)}_6 &=24.05941486\ldots \\ \g^{(3)}_7
&=26.57786874\ldots \\ \g^{(3)}_8 &=28.21816451\ldots \\
\g^{(3)}_9 &=30.74504026\ldots \\ \g^{(3)}_{10}
&=33.89738893\ldots \\ \g^{(3)}_{11} &=35.60841265\ldots \\
\g^{(3)}_{12} &=37.55179656\ldots \\ \g^{(3)}_{13}
&=39.48520726\ldots \\ \g^{(3)}_{14} &=42.61637923\ldots \\
\g^{(3)}_{15} &=44.12057291\ldots \\ \g^{(3)}_{16}
&=46.27411802\ldots \\ \g^{(3)}_{17} &=47.51410451\ldots \\
\g^{(3)}_{18} &=50.37513865\ldots 
\end{align*}

This time the "deltas" are differences
$$
\d_{i,j}=\g_i-\g^{(3)}_j \ .
$$

We perform statistics (a) with $1\leq i,j\leq 100 \ 000$ and
statistics (b) with $1\leq i,j\leq 10^6$. We look at deltas in
$[0,50]$ with precision $0.1$. The results are presented in the
following figures.


\begin{center}
  \resizebox{6cm}{!}{\includegraphics{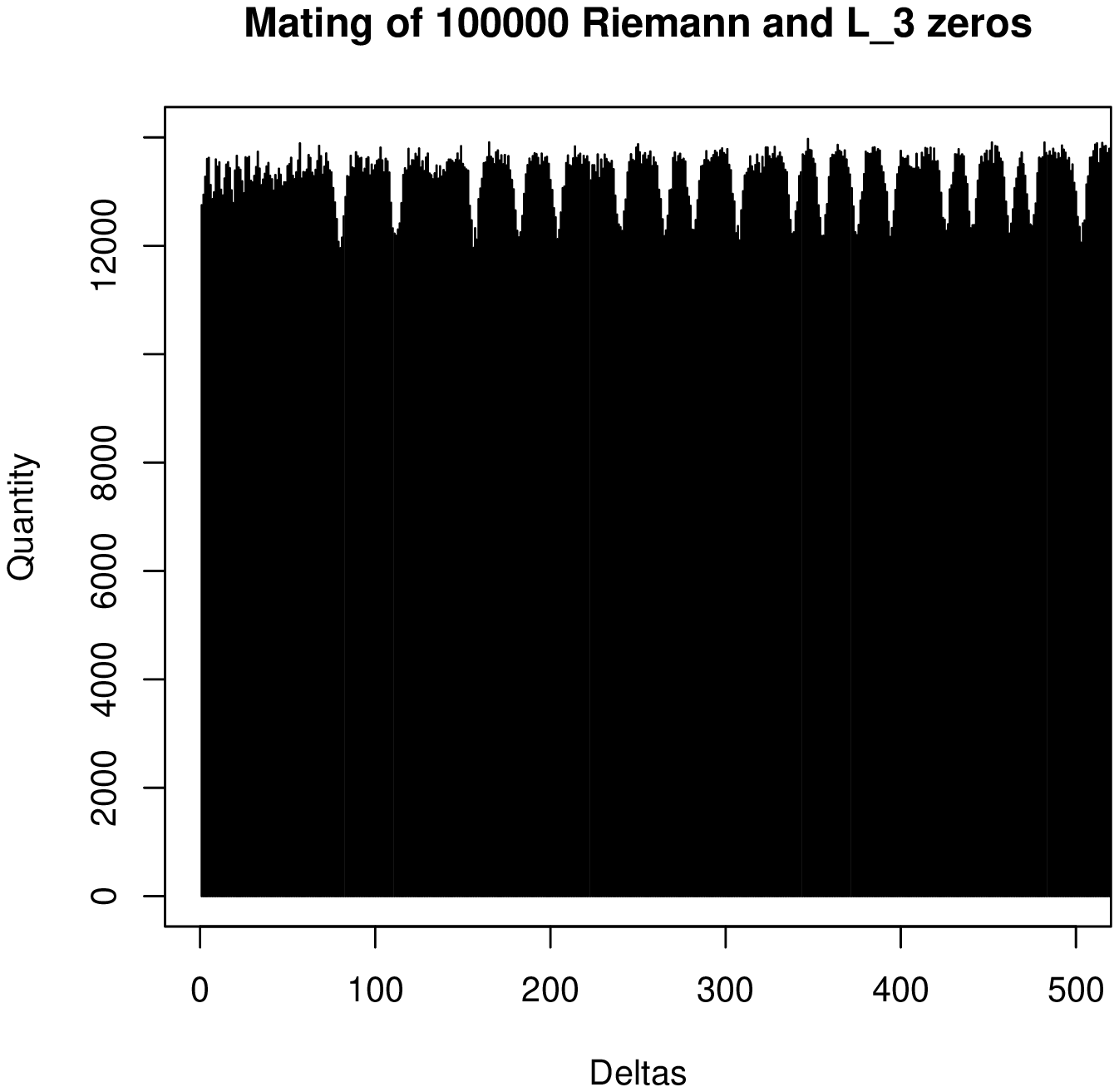}}    
\resizebox{6cm}{!}{\includegraphics{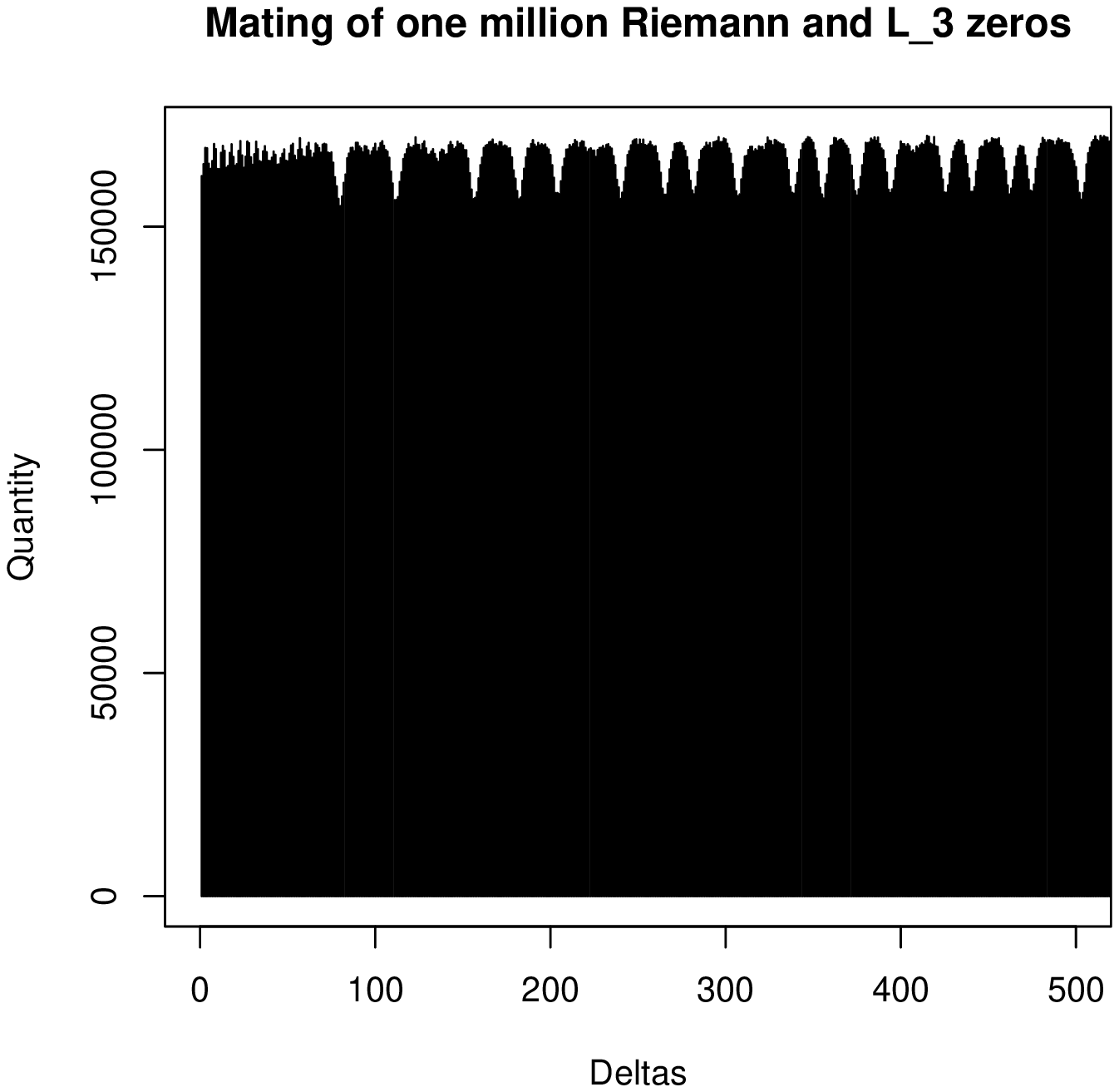}}
\end{center}

\centerline{Figures 22.a and 22.b.}


\begin{center}
  \resizebox{6cm}{!}{\includegraphics{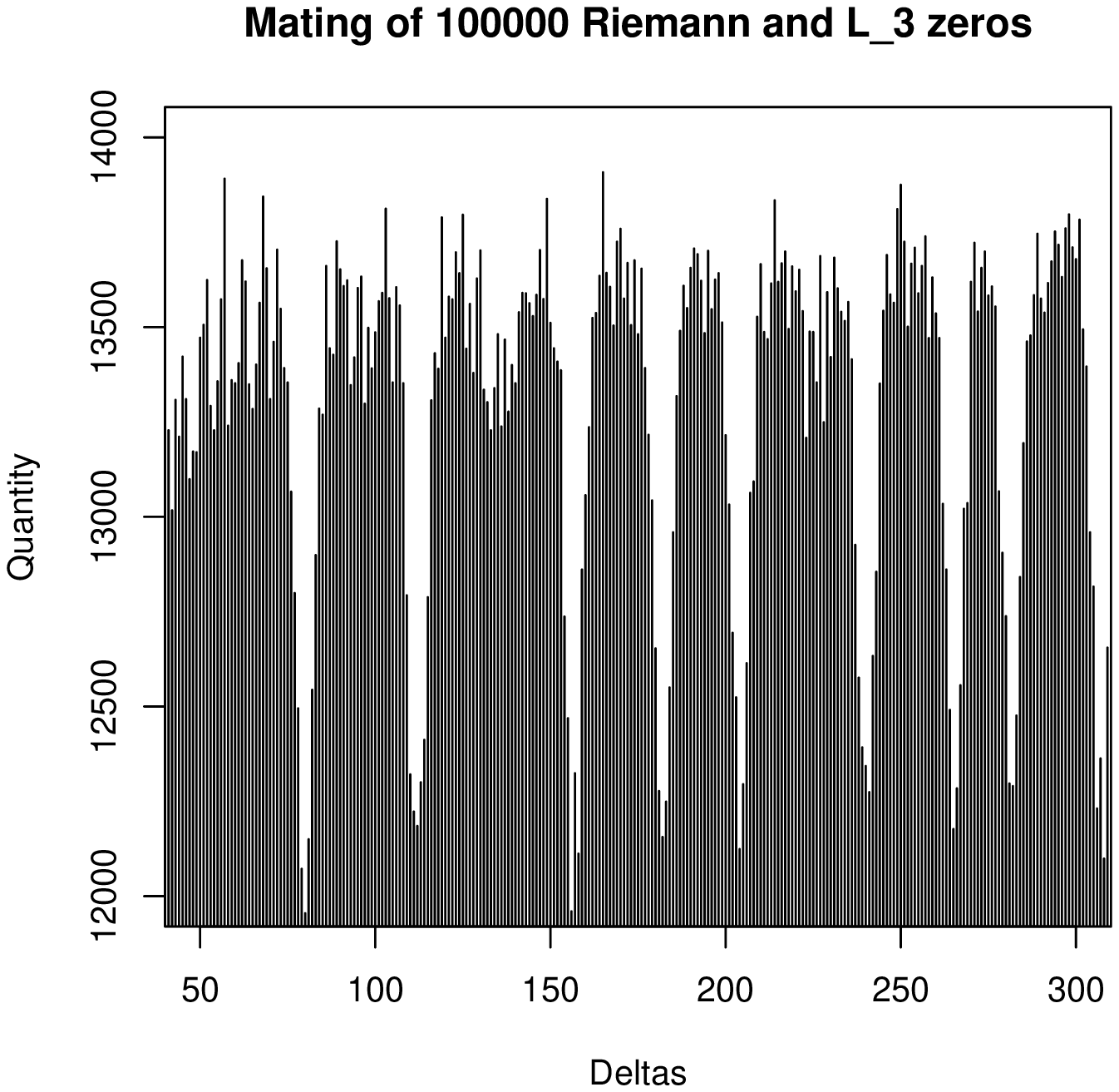}}    
\resizebox{6cm}{!}{\includegraphics{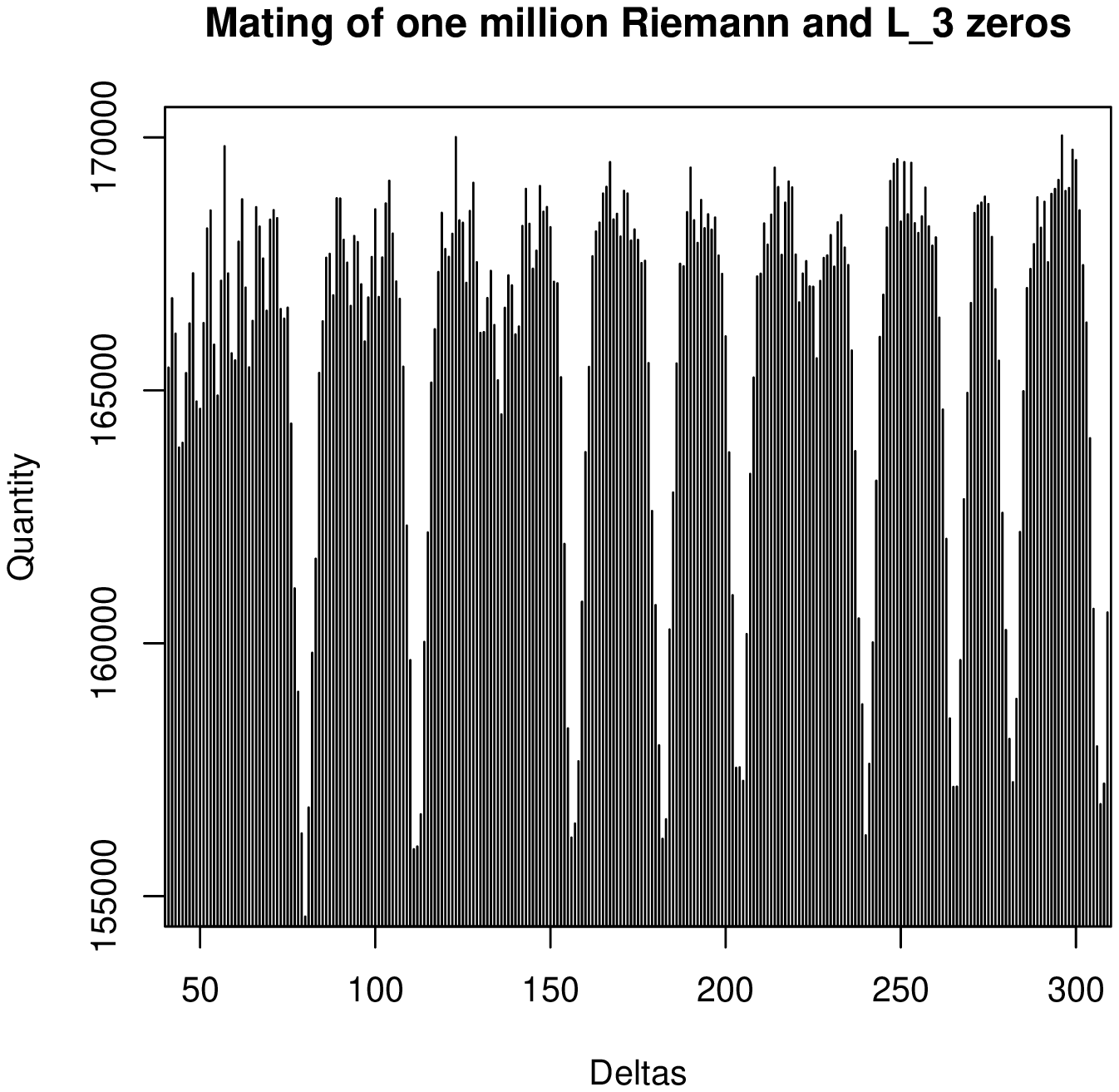}}
\end{center}

\centerline{Figures 23.a and 23.b.}


\begin{center}
  \resizebox{6cm}{!}{\includegraphics{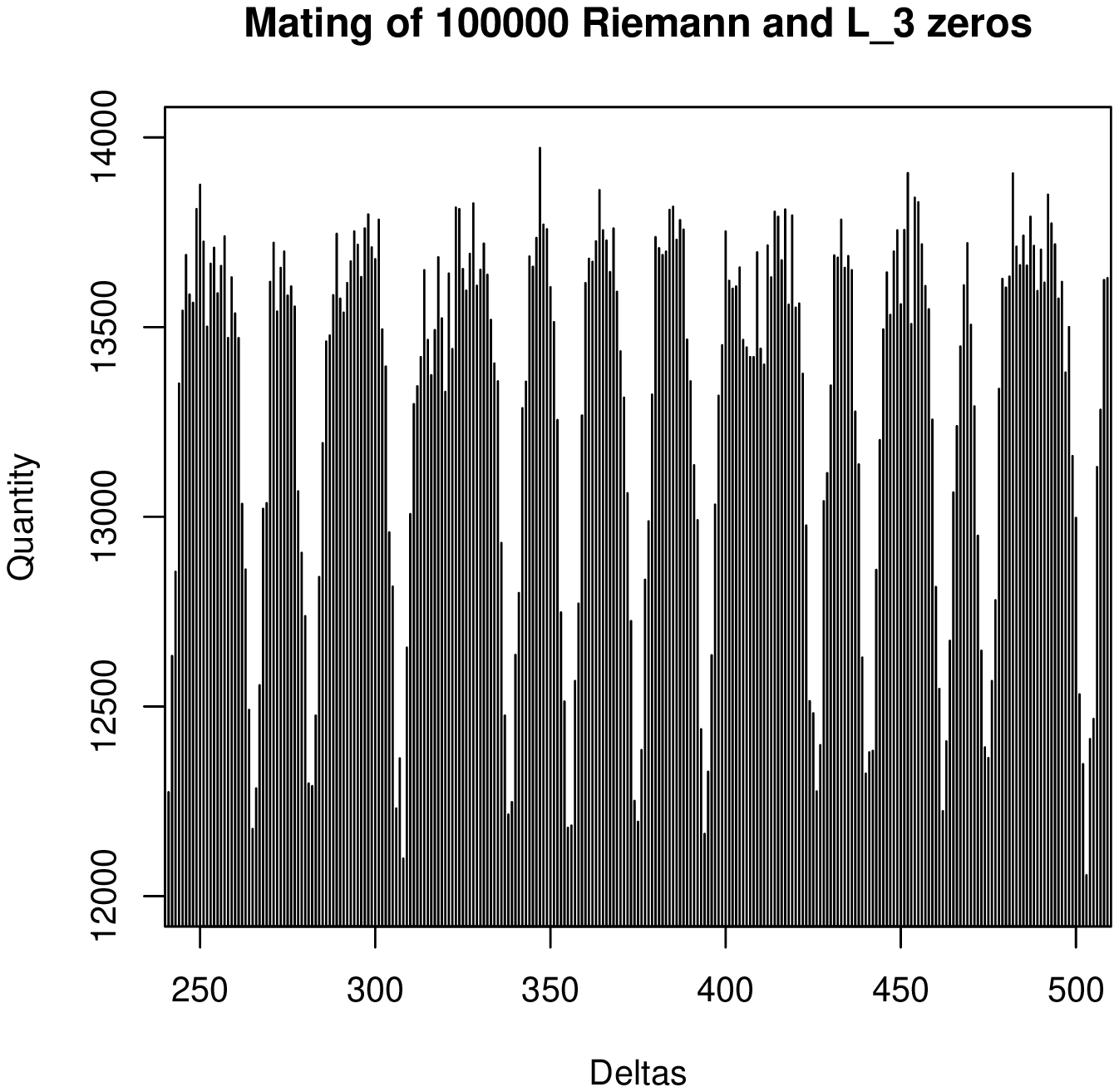}}    
\resizebox{6cm}{!}{\includegraphics{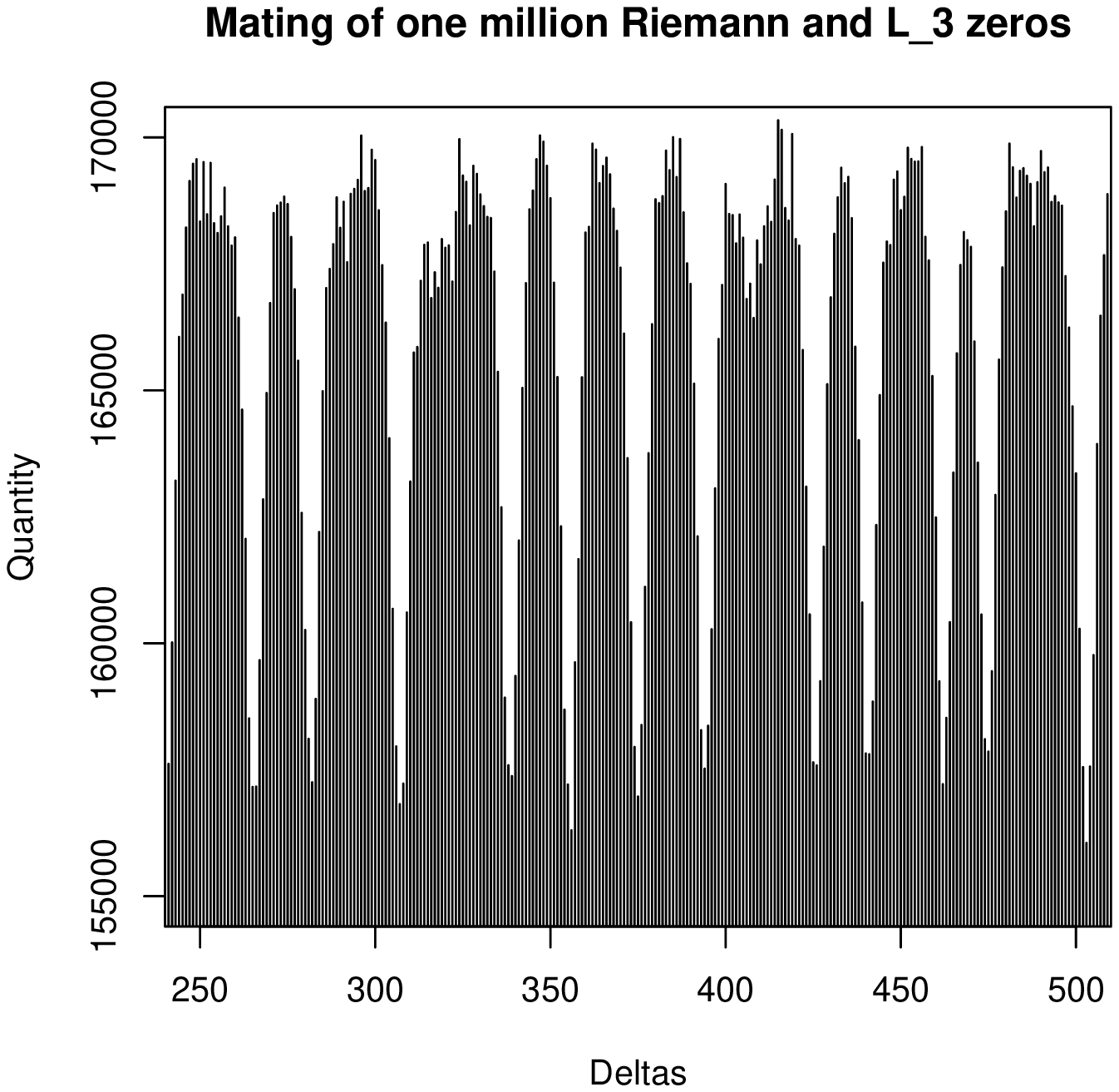}}
\end{center}

\centerline{Figures 24.a and 24.b.}

\bigskip

We observe that this time the deficient locations for the deltas
happen exactly at the location of the zeros of $L_{\chi_3}$. We
easily recognize in figures 23.a and 23.b the location of the first
zeros of $L_{\chi_3}$. We can check the full list of zeros less
than $50$ by looking also at the figures 24.a and 24.b. We conclude
that the zeros of $L$-functions replicate mating them with Riemann
zeros.

A new feature is that near $0$ we no longer have a GUE
distribution for the deltas. As the theory of the e\~ne product
explains, the deficit at $0$ only occurs when 
we have symmetric zeros, i.e. we mate the zeros of $L_{\chi_1}$ with 
those of $L_{\chi_2}$ when
$$
\chi_1=\bar \chi_2 \ ,
$$
and we have an atomic mass at $0$
that comes from the sum of symmetric zeros of $L_{\chi_1}$ and
$L_{\bar \chi_1}$. Thus if the character is not real, then we don't have 
a GUE distribution, not even a deficit, but the Riemann Hypothesis is still
conjectured, thus there is no direct relation between the Riemann Hypothesis 
and Montgomery Conjecture. The author knows no reference in the literature for this 
observation.


\begin{center}
  \resizebox{6cm}{!}{\includegraphics{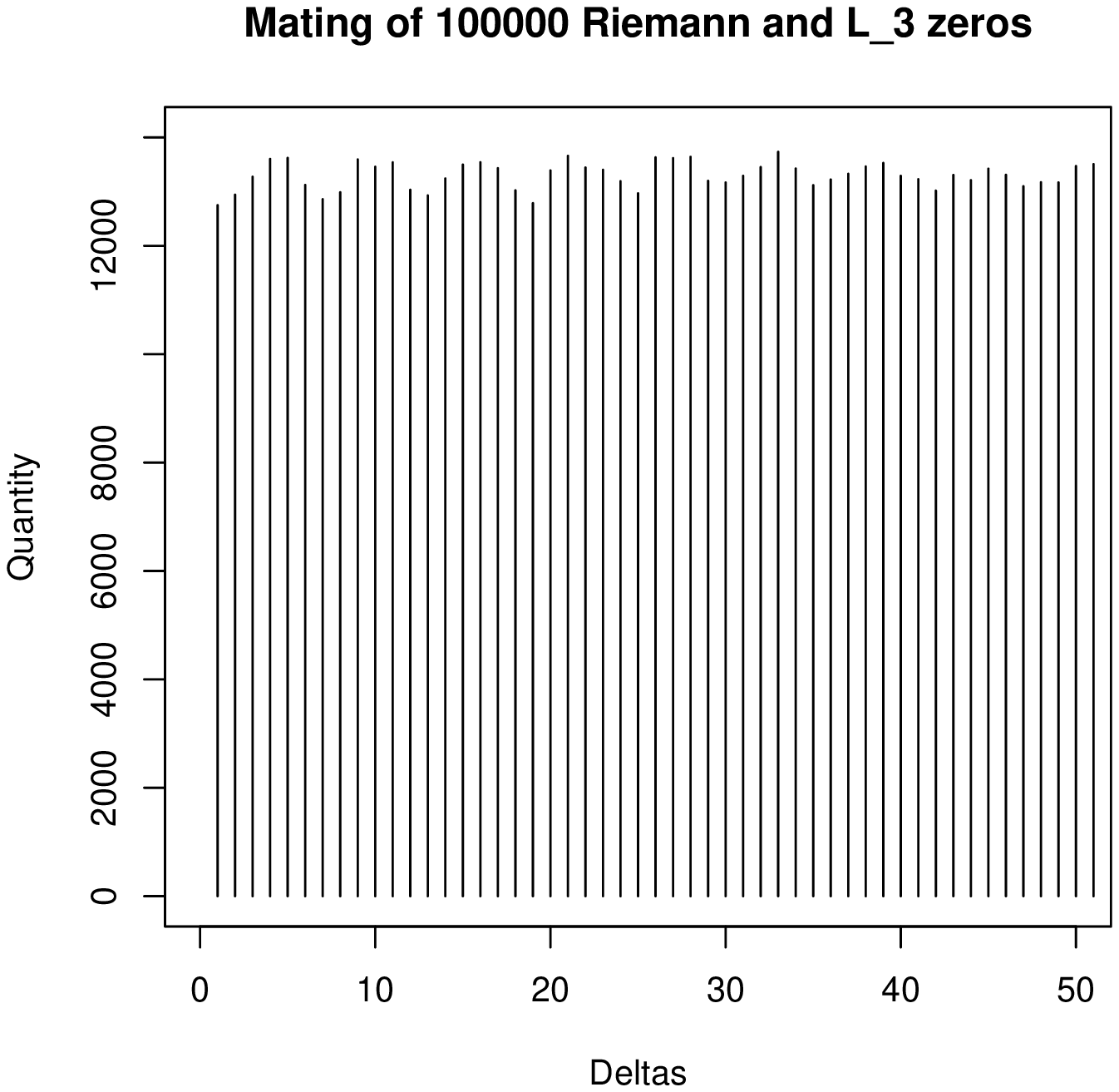}}    
\resizebox{6cm}{!}{\includegraphics{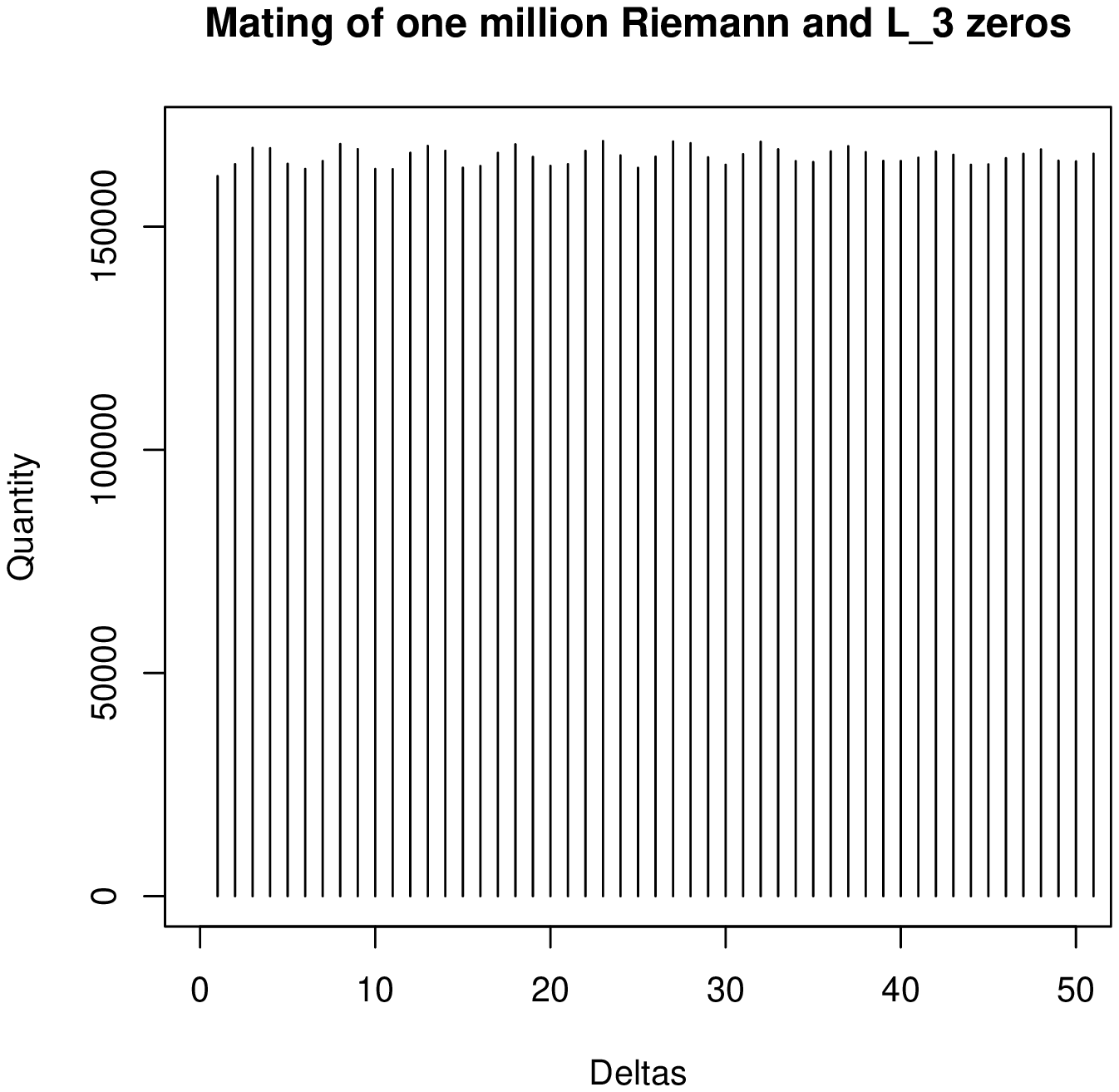}}
\end{center}

\centerline{Figures 25.a and 25.b.}

\medskip

Next we perform the same mating statistics of Riemann zeros with
zeros of   $L_{\chi_{7,3}}$.  Recall that the zeros of this
non-real analytic $L$-functions are not symmetric with respect to
$0$. We perform two statistics. We consider the first $100\ 000$
Riemann zeros and compute all deltas with positive (resp. negative
taking their negative value) zeros of $L_{\chi_{7,3}}$. The list
of the first positive zeros of $L_{\chi_{7,3}}$ less than $50$ are

\begin{align*}
\g^{(7+)}_1 &=4.356402\ldots \\ \g^{(7+)}_2
&=8.785555\ldots \\ \g^{(7+)}_3 &=10.736120\ldots \\ \g^{(7+)}_4
&=12.532548\ldots \\ \g^{(7+)}_5 &=15.937448\ldots \\
\g^{(7+)}_6 &=17.616053\ldots \\ \g^{(7+)}_7 &=20.030559\ldots
\\ \g^{(7+)}_8 &= 21.314647\ldots \\ \g^{(7+)}_9
&=23.203672\ldots \\ \g^{(7+)}_{10} &=26.169945\ldots \\
\g^{(7+)}_{11} &=27.873375\ldots \\ \g^{(7+)}_{12}
&=28.599794\ldots \\ \g^{(7+)}_{13} &=30.919561\ldots \\
\g^{(7+)}_{14} &=32.610089 \ldots \\ \g^{(7+)}_{15}
&=34.792503\ldots \\ \g^{(7+)}_{16} &=36.344756\ldots \\
\g^{(7+)}_{17} &=38.206755\ldots \\ \g^{(7+)}_{18}
&=39.338483\ldots \\ \g^{(7+)}_{19} &=40.476472\ldots
\\ \g^{(7+)}_{20} &= 43.539481\ldots \\ \g^{(7+)}_{21}
&=44.595772\ldots \\ \g^{(7+)}_{22} &=46.096099\ldots \\
\g^{(7+)}_{23} &=47.491559\ldots \\ \g^{(7+)}_{24}
&=49.126475\ldots \\ &\vdots 
\end{align*}

The list of the first negative zeros of $L_{\chi_{7,3}}$ less than
$51$ are

\begin{align*}
\g^{(7-)}_1 &=6.201230\ldots \\ \g^{(7-)}_2
&=7.927431\ldots \\ \g^{(7-)}_3 &=11.010445\ldots \\ \g^{(7-)}_4
&=13.829868\ldots \\ \g^{(7-)}_5 &=16.013727\ldots \\
\g^{(7-)}_6 &=18.044858\ldots \\ \g^{(7-)}_7 &=19.113886\ldots
\\ \g^{(7-)}_8 &= 22.756406\ldots \\ \g^{(7-)}_9
&=23.955938\ldots \\ \g^{(7-)}_{10} &=25.723104\ldots \\
\g^{(7-)}_{11} &=27.455596\ldots \\ \g^{(7-)}_{12}
&=29.338505\ldots \\ \g^{(7-)}_{13} &=31.284265\ldots \\
\g^{(7-)}_{14} &=33.672299 \ldots \\ \g^{(7-)}_{15}
&=34.774195\ldots \\ \g^{(7-)}_{16} &=35.973150\ldots \\
\g^{(7-)}_{17} &=37.786921\ldots \\ \g^{(7-)}_{18}
&=40.224566\ldots \\ \g^{(7-)}_{19} &=41.909138\ldots
\\ \g^{(7-)}_{20} &= 42.712631\ldots \\ \g^{(7-)}_{21} &=44.977200
\ldots \\ \g^{(7-)}_{22} &=46.086774\ldots \\ \g^{(7-)}_{23}
&=47.348801\ldots \\ \g^{(7-)}_{24} &=50.017326\ldots \\ &\vdots
\end{align*}

The following figures show the result of the numerical statistics.
Figures (a), resp. (b), are for the mating against positive, resp.
negative, zeros. We observe for statistics (a) that the deficient
locations do correspond to values of the positive zeros. For
statistics (b) we observe these locations at the values of the
negative zeros. In particular in figures 27 we appreciate the
location of the first zeros.


\begin{center}
  \resizebox{6cm}{!}{\includegraphics{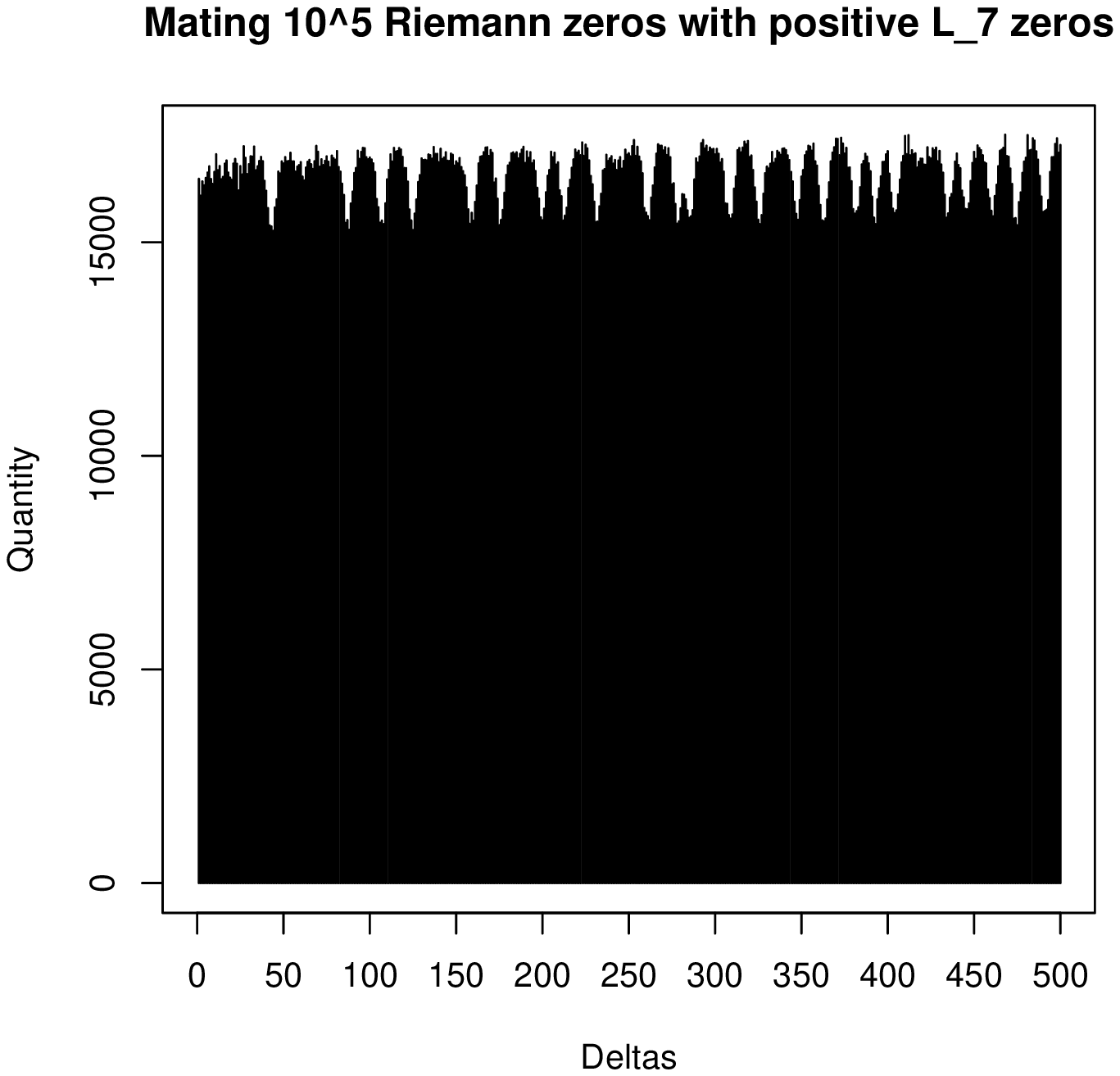}}    
\resizebox{6cm}{!}{\includegraphics{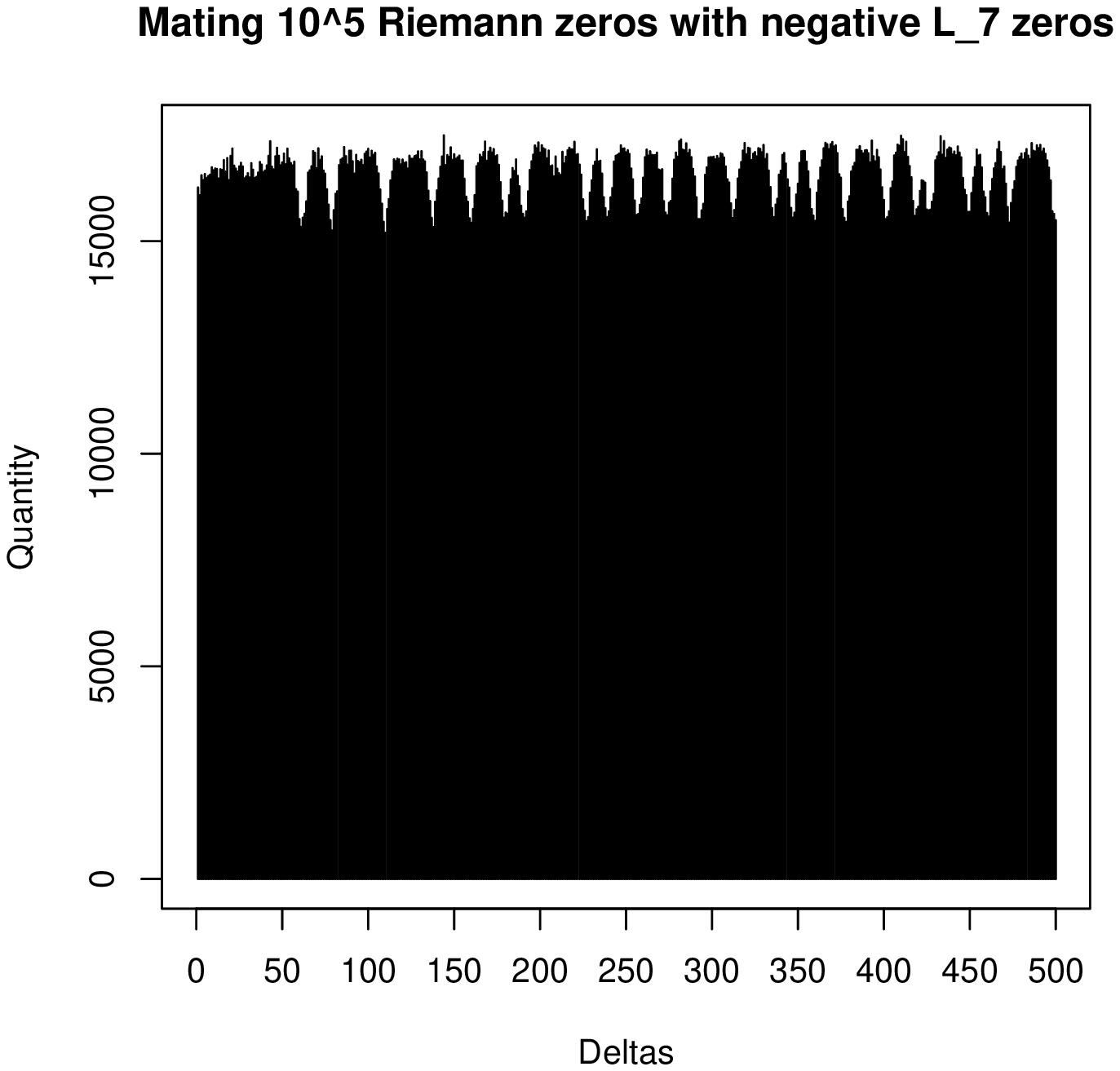}}
\end{center}

\centerline{Figures 26.a and 26.b.}


\begin{center}
  \resizebox{6cm}{!}{\includegraphics{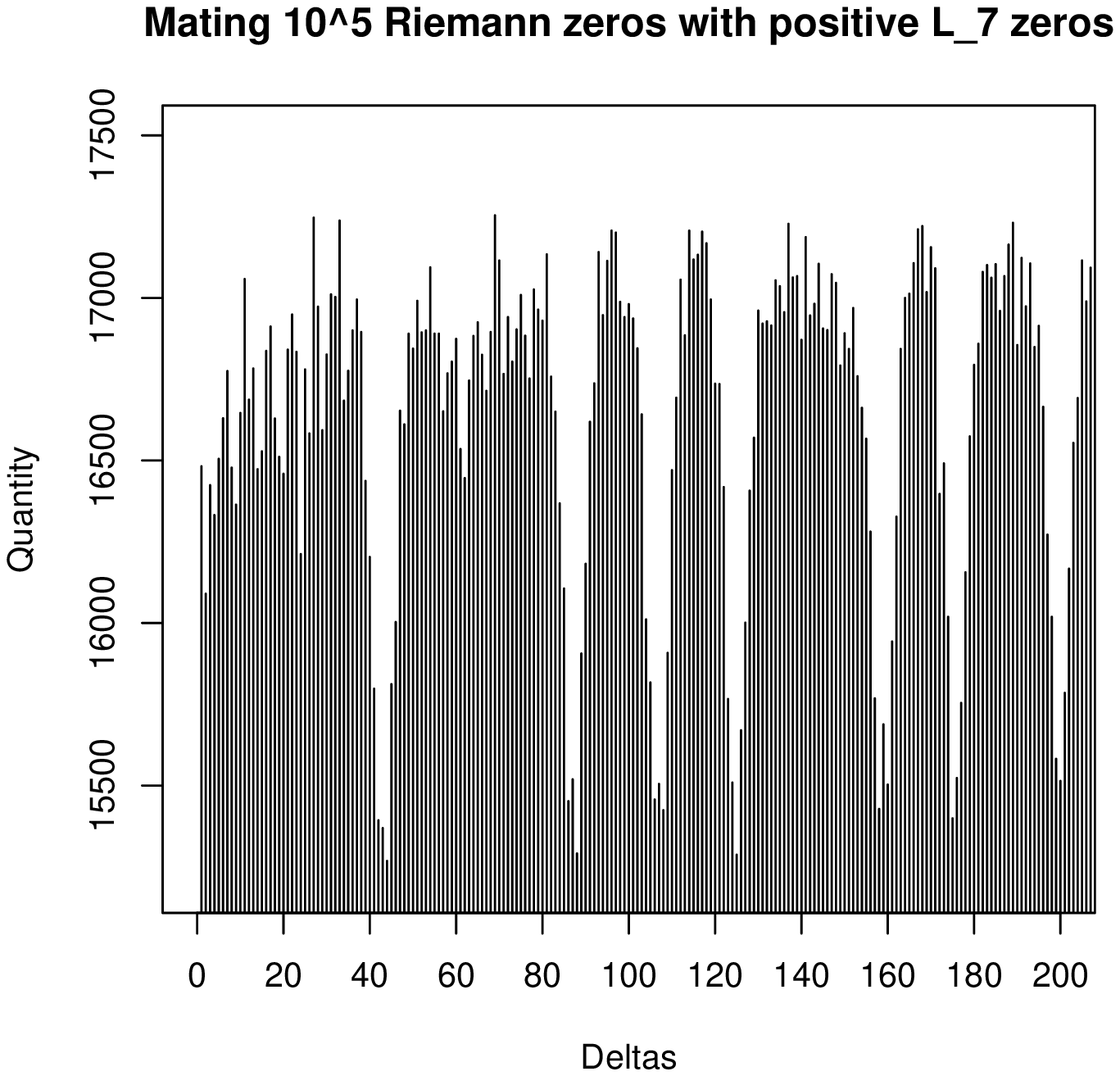}}    
\resizebox{6cm}{!}{\includegraphics{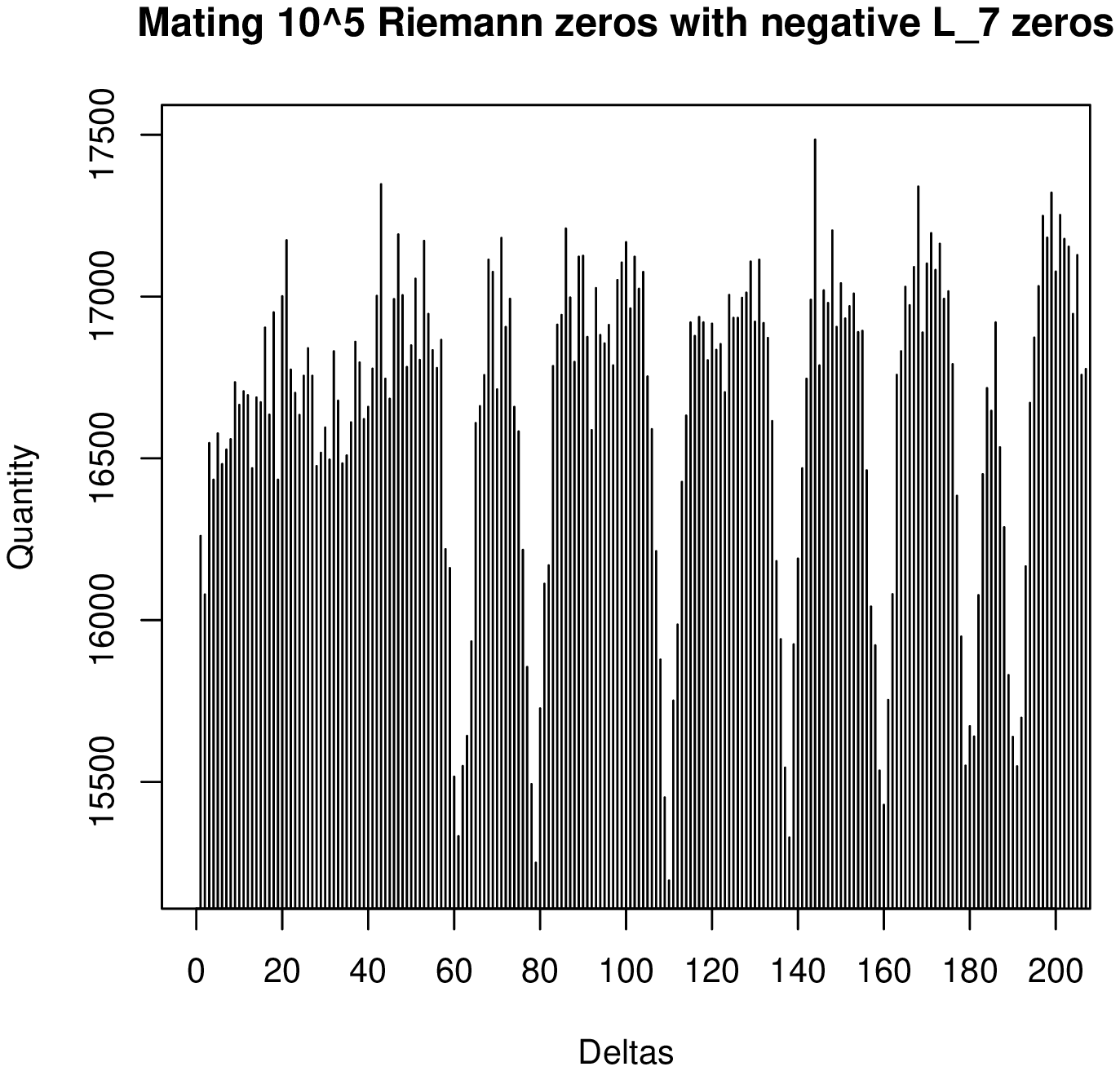}}
\end{center}

\centerline{Figures 27.a and 27.b.}

%

\begin{center}
  \resizebox{6cm}{!}{\includegraphics{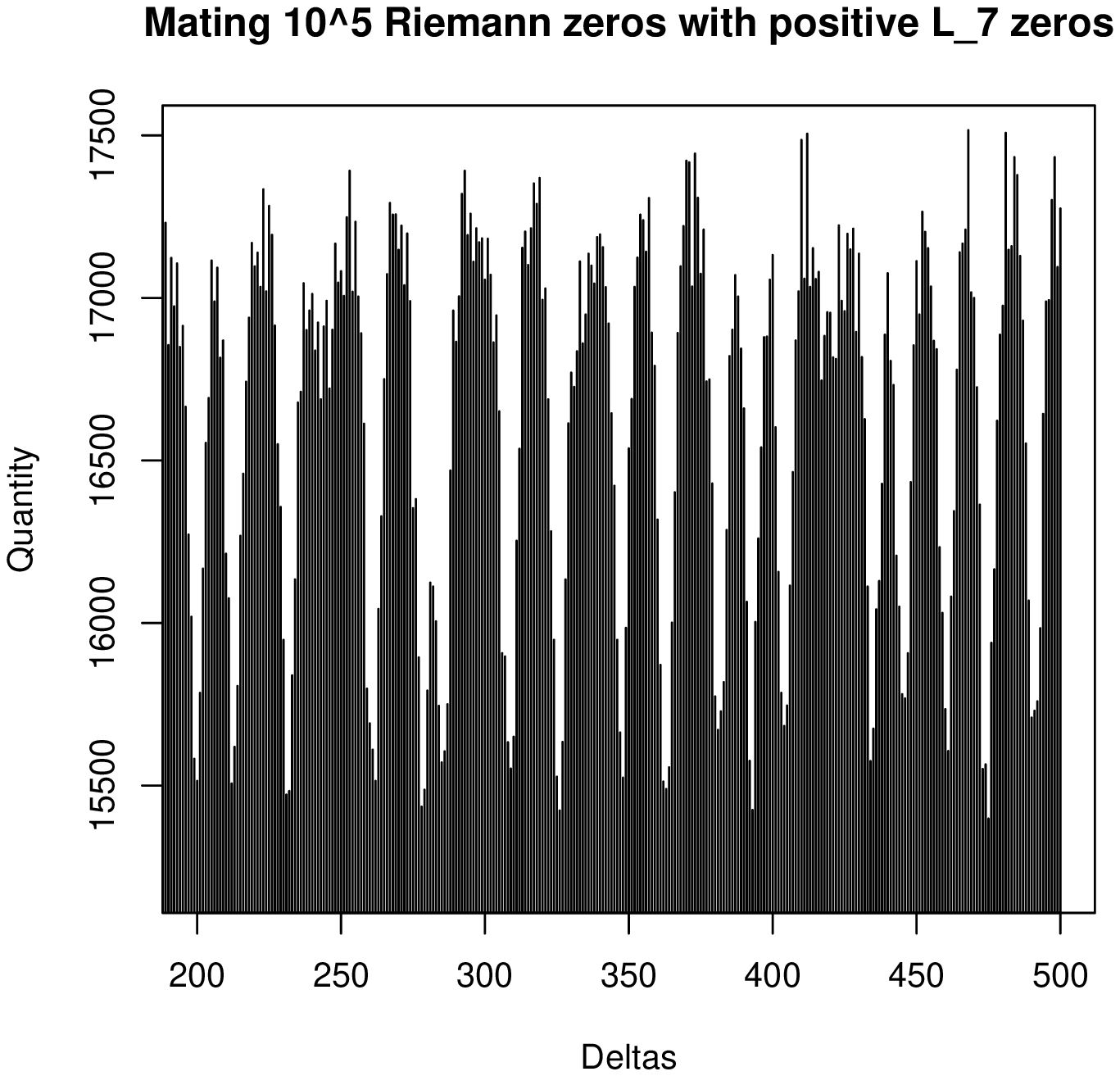}}    
\resizebox{6cm}{!}{\includegraphics{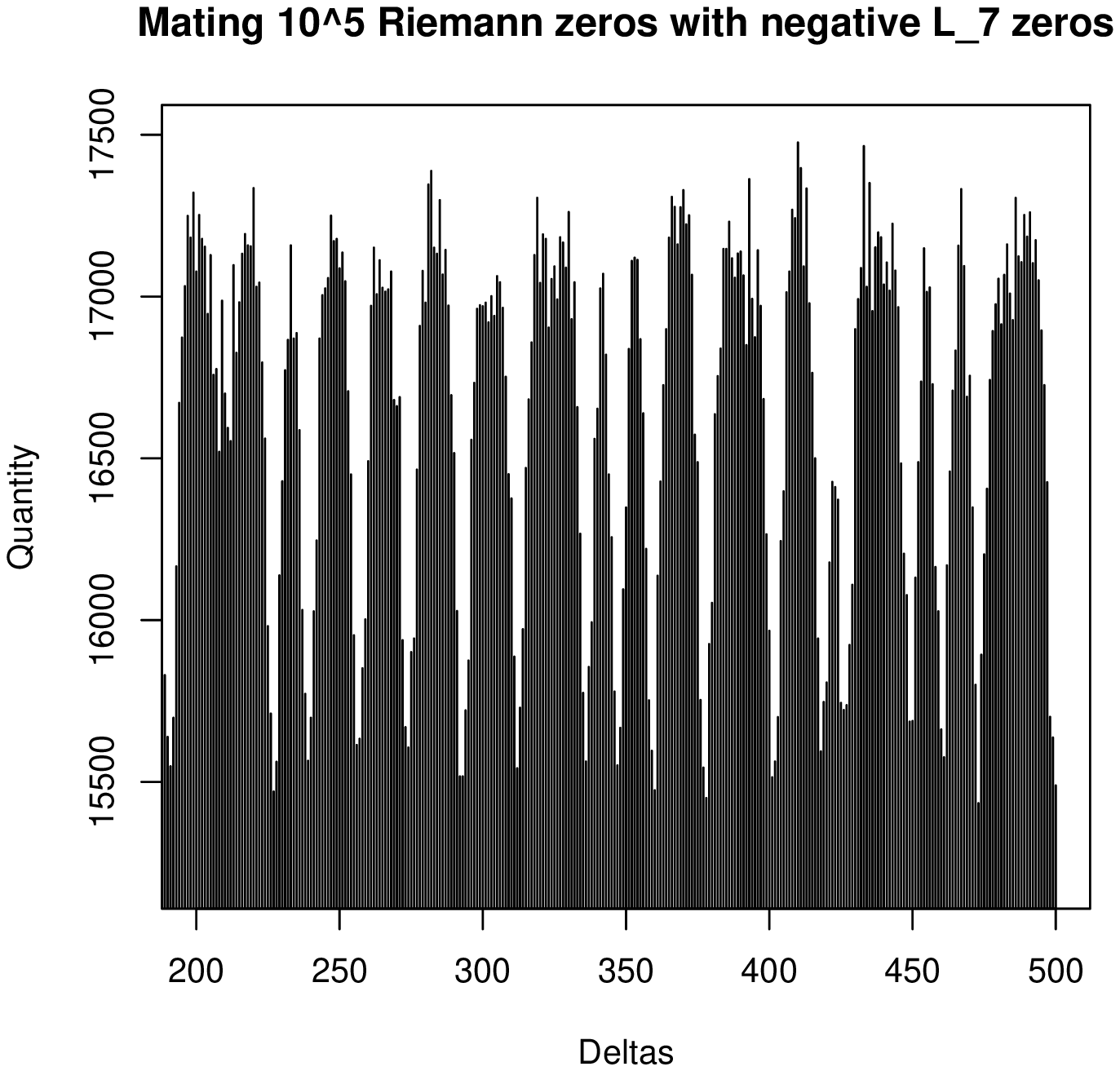}}
\end{center}

\centerline{Figures 28.a and 28.b.}

\bigskip

We perform a final statistic in order to check the nonexistence of 
the GUE distribution, and not even a deficit of deltas effect near $0$. 
We compute the deltas near
$0$ against half million Riemann zeros with precision $0.01$. The
results are shown in Figures 29, figure 29.a for positive deltas
and figure 29.b for negative ones. Figures 29 show the deltas with
double precision in the range $[0,2]$. The reader can compare
directly these figures  with Figures 21.c and 21.d. The conclusion
is clear: No GUE distribution near $0$.


\begin{center}
  \resizebox{6cm}{!}{\includegraphics{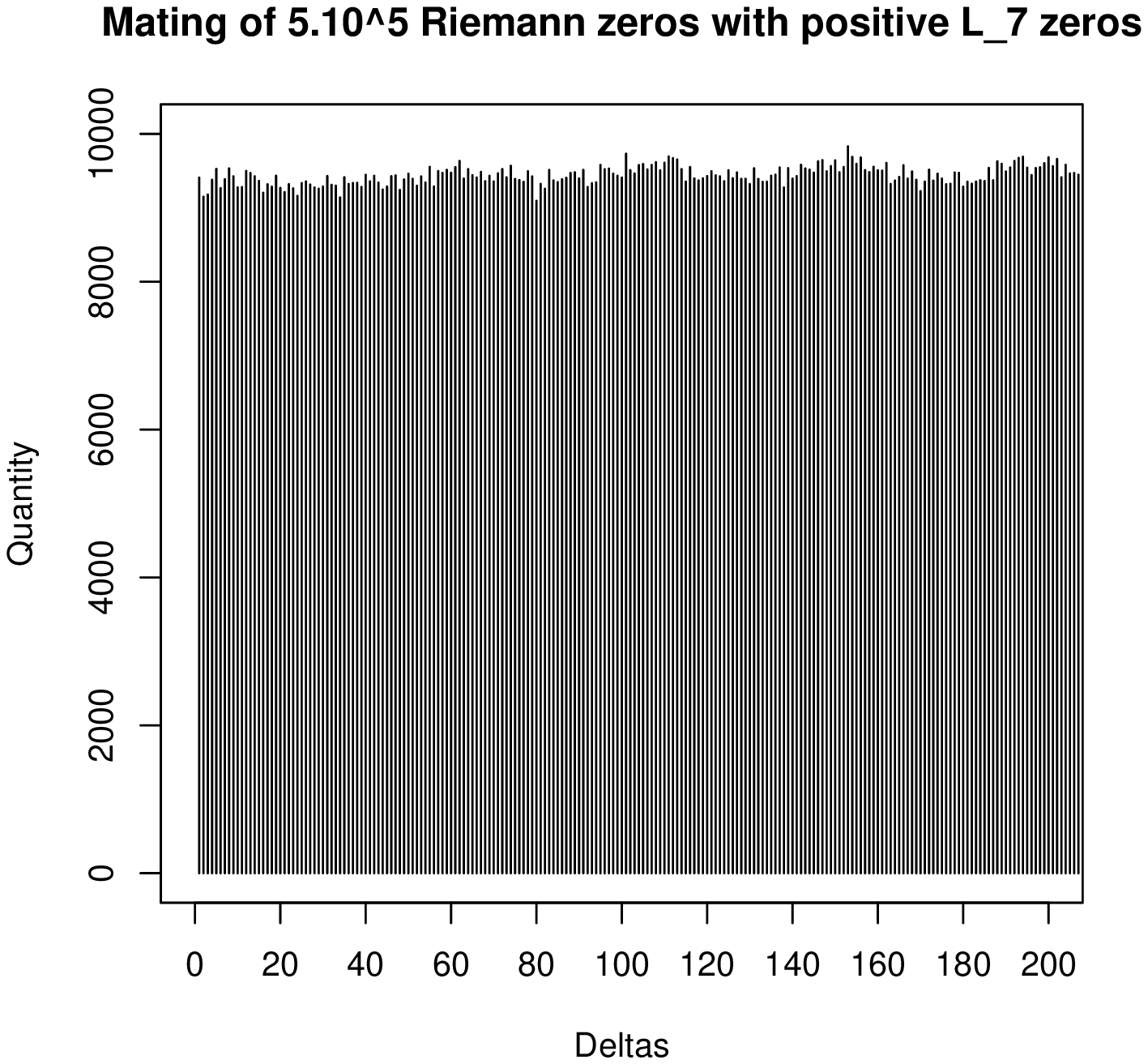}}    
\resizebox{6cm}{!}{\includegraphics{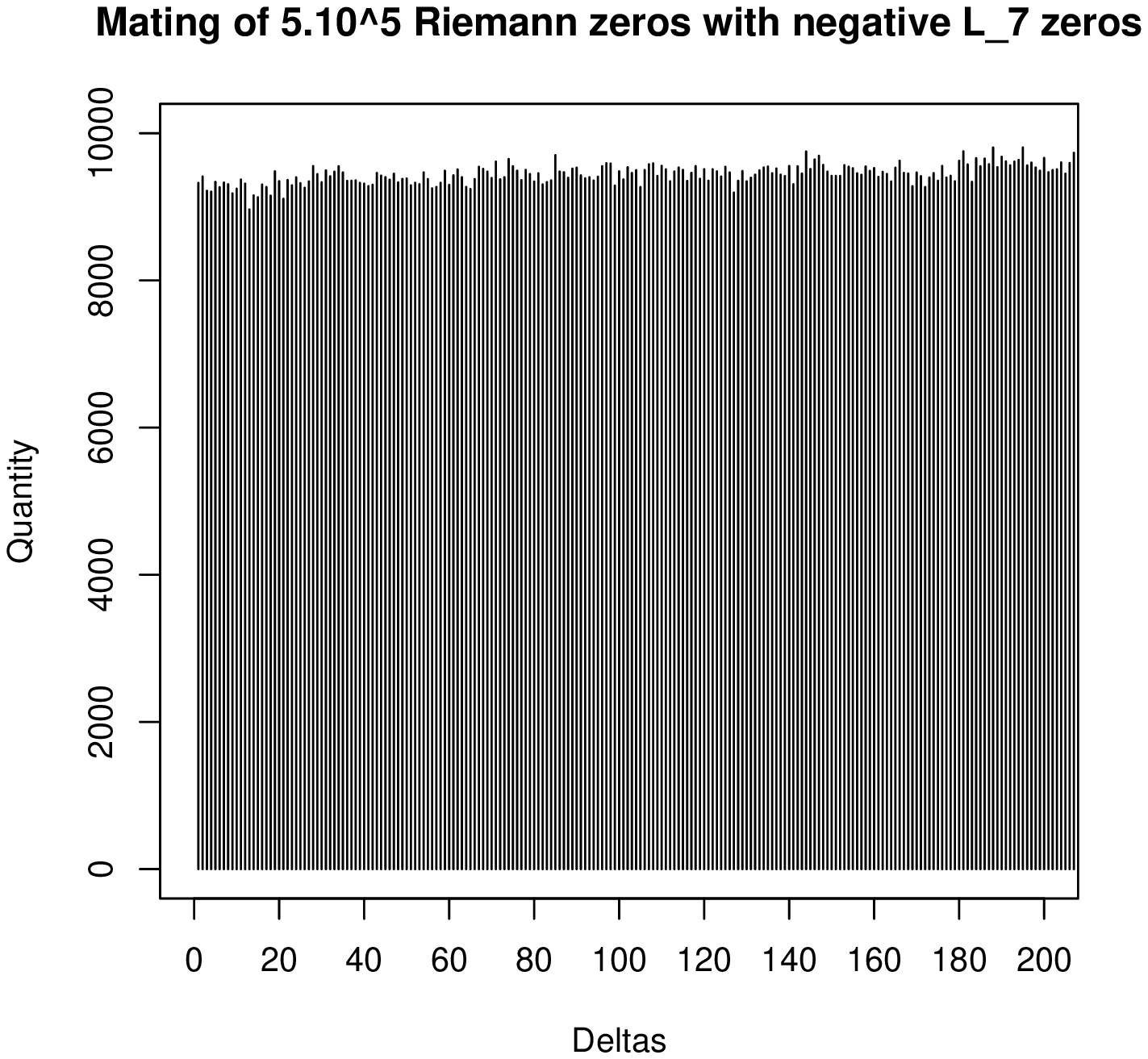}}
\end{center}

\centerline{Figures 29.a and 29.b.}

\bigskip

We check also from these statistics the location with double
precision the first positive and negative zero in figures 30.a and
30.b which show the deltas in the range $[0,10]$ with double
precision. We can appreciate distinctly with double precision in
figure 30.a both positive zeros less than $10$
$$
\g^{(7+)}_1 =4.356402\ldots \ \ \  \g^{(7+)}_2 =8.785555\ldots
$$
and in figure 30.b both "negative" zeros less than $10$
$$
\g^{(7-)}_1 =6.201230\ldots \ \ \ \g^{(7-)}_2 =7.927431\ldots
$$


\begin{center}
  \resizebox{6cm}{!}{\includegraphics{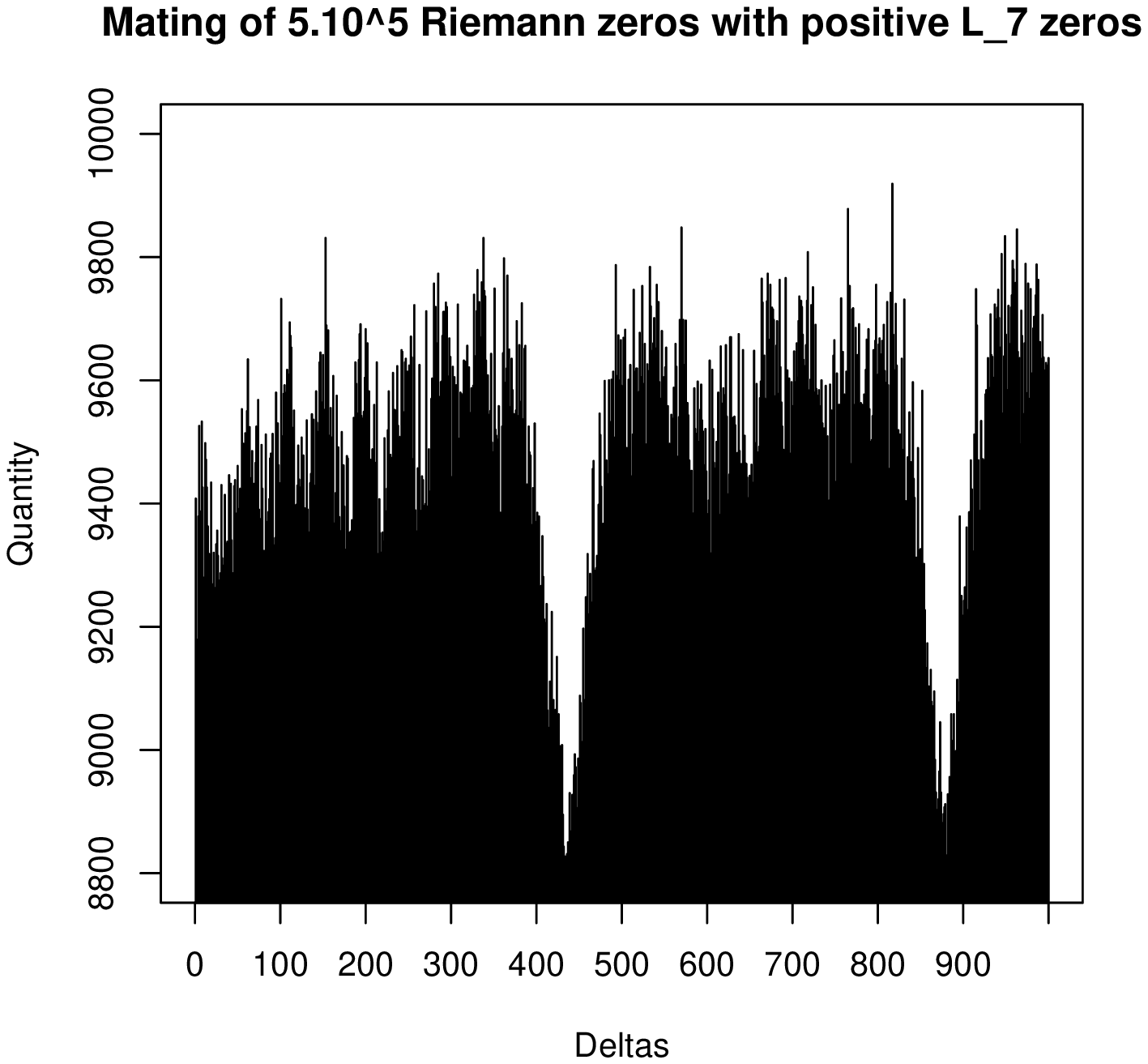}}    
\resizebox{6cm}{!}{\includegraphics{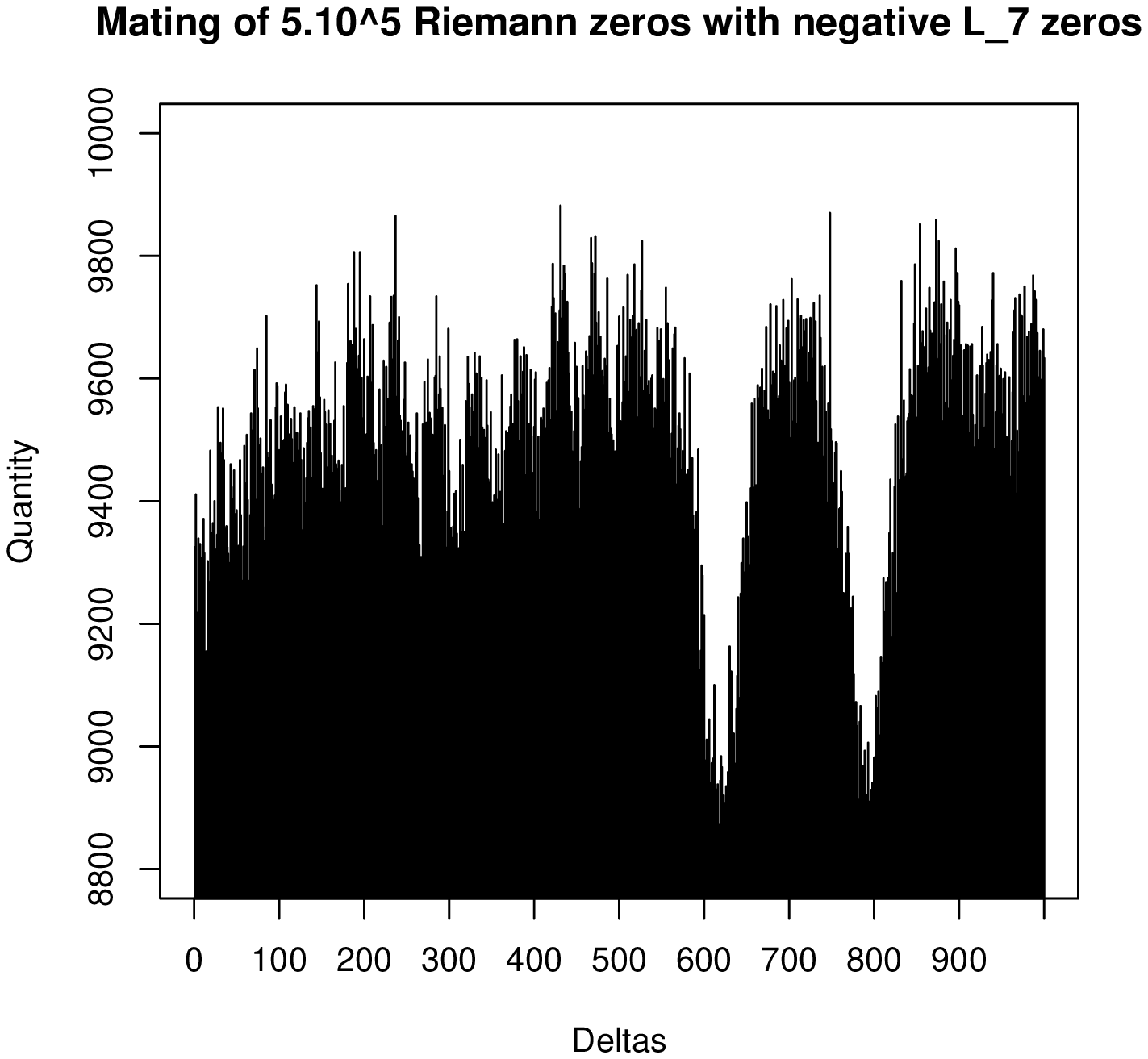}}
\end{center}

\centerline{Figures 30.a and 30.b.}

\medskip

\textbf{E\~ne product explanation.}

\medskip

The computation follows. We have for any Dirichlet $L$-function
$L_\chi$,

\begin{align*}
L_\chi \bar \star \bar \zeta &=L_\chi \bar \star \zeta
\\ &=\prod_p (1-\chi (p)p^{-s})^{-1} \ \bar \star \ \prod_q
(1-q^{-s})^{-1} \\ &=\prod_p (1-\chi (p)p^{-s})^{-1}\ \bar \star
\ (1-p^{-s})^{-1} \\ &=\prod_p (1-\chi (p)p^{-1/2} p^{-s}) \\
&=\prod_p (1-\chi (p)  p^{-(s+1/2)}) \\ &=L_\chi (s+1/2)^{-1} \ \
.
\end{align*}

Therefore we recognize that the mating of zeros of $L_\chi$ with
Riemann zeros has deficient deltas at the location corresponding
to the imaginary part of zeros of $L_\chi$.

\medskip

\textbf{Scripts.}

\medskip

Below is the script we used in order to produce the previous
figures. The script is slightly different from previous ones 
since we compute
separately positive and negative zeros. The cumulative positive
deltas are stored in the list "x" and the negative in the list
"y". The zeros of Dirichlet $L$-function of conductor $7$ are
stored in Rubinstein's file "zeros-0007-2000000", and Riemann
zeros are from Rubinstein's file "zeros-0001-35161820".

\medskip

{\tt

zerosL<-read.table("zeros-0007-2000000",skip=2000000,nrows=1000000)

z<-zerosL[,3]

z.plus<-z[z>0]

z.minus<--z[z<0]

zerosR<-scan("zeros-0001-35161820",skip=0,nlines=100000)

x=rep(0,500)

zeta1<-zerosR

zeta2<-z.minus

N=length(zeta) j=1

for (i in 1:N)

$\{$

while ( (zeta1[i]-zeta2[j])>50.1  )

$\{$

j=j+1

$\}$

l=0

while ( ((zeta1[i]-zeta2[j+l])>0) \& ((zeta1[i]-zeta2[j+l])<50.1)
)

$\{$

d=10*(zeta[i]-zeta3[j+l])

x[as.integer(d)]=x[as.integer(d)]+1

l=l+1

$\}$

$\}$

y=rep(0,500)

zeta1<-zerosR

zeta3<-z.plus

N=length(zeta)

j=1

for (i in 1:N)

$\{$

while ( (zeta1[i]-zeta3[j])>50.1  ) { j=j+1 }

l=0

while ( ((zeta1[i]-zeta3[j+l])>0) \& ((zeta1[i]-zeta3[j+l])<50.1)
)

$\{$

d=10*(zeta1[i]-zeta3[j+l])

y[as.integer(d)]=y[as.integer(d)]+1

l=l+1

$\}$

$\}$

}

\section{Mating of general $L$-functions.} \label{sec:L-functions_mating}

In this section we perform similar statistics to those in the
previous section but mating the zeros of two Dirichlet
$L$-functions $L_{\chi_1}$ and $L_{\chi_2}$. This time we observe
that the deficient locations for the statistics of deltas
correspond to the zeros of an arithmetically  well determined,
namely  $L_{\chi_1 \bar \chi_2}$.

For a character $\chi$ we denote by $f_\chi$ its conductor. We
have
$$
f_{\bar \chi}=f_\chi \ .
$$
All characters considered are primitive, i.e.
defined modulo its conductor. Let $\chi_1$
and $\chi_2$ be  two characters. If $f_{\chi_1}\wedge f_{\chi_2}=1$ then the
conductor of $\chi_1 \bar \chi_2$ is
$$
f_{\chi_1 \bar \chi_2}=f_{\chi_1}.f_{\chi_2} \ .
$$

The first complex non-real Dirichlet character has conductor $5$.
Therefore the mating of two Dirichlet $L$-functions of complex
non-real characters with distinct conductors has conductor at
least $35$. We have only access to Rubinstein's public data that
contains large files of zeros for  Dirichlet $L$-functions with
conductor $\leq 19$. Therefore we limit our numerical computation
to real characters for which we can check the result with the
available data. This is done only for checking purposes. Note that
we could indeed compute, with a rough precision, the zeros of
higher conductor Dirichlet $L$-functions (for example $35$) by
using Rubinstein's data of conductors $\leq 19$.

We choose to mate the zeros of Dirichlet $L$-functions
$L_{\chi_3}$ of conductor $3$, and $L_{\chi_4}$ of conductor $4$.
We should obtain the zeros of the only Dirichlet $L$-function of
conductor $12$, $L_{\chi_{12}}$. The list of the first zeros of
$L_{\chi_{12}}$ less than $50$ is

\begin{align*}
\g^{(12)}_1 &=3.8046276331\ldots \\ \g^{(12)}_2
&=6.6922233205\ldots \\ \g^{(12)}_3 &=8.8905929587\ldots \\
\g^{(12)}_4 &=11.188392745\ldots \\ \g^{(12)}_5
&=12.966178808\ldots \\ \g^{(12)}_6 &=15.181480876\ldots \\
\g^{(12)}_7 &=16.632633275\ldots \\ \g^{(12)}_8 &=
18.884369457\ldots \\ \g^{(12)}_9 &=20.103928191\ldots \\
\g^{(12)}_{10} &=22.285839107\ldots \\ \g^{(12)}_{11}
&=23.561319713\ldots \\ \g^{(12)}_{12} &=25.411633892\ldots \\
\g^{(12)}_{13} &=27.013943986\ldots \\ \g^{(12)}_{14}
&=28.442203258 \\
 \g^{(12)}_{15} &=30.204006556\ldots \\ \g^{(12)}_{16}
&=31.648077615\ldots \\ \g^{(12)}_{17} &=33.03713288\ldots \\
\g^{(12)}_{18} &=35.027378485\ldots \\ \g^{(12)}_{19}
&=35.778044577\ldots \\ \g^{(12)}_{20} &= 37.926816821\ldots
\\\g^{(12)}_{21} &=38.973998822 \ldots \\ \g^{(12)}_{22}
&=40.484154751\ldots \\ \g^{(12)}_{23} &=42.235143018\ldots \\
\g^{(12)}_{24} &=43.192847103\ldots \\ \g^{(12)}_{25}
&=44.948822502\ldots \\ \g^{(12)}_{26} &=46.243369979\ldots
\\ \g^{(12)}_{27} &=47.646400501\ldots \\ \g^{(12)}_{28}
&=48.943728012\ldots \\ & \vdots 
\end{align*}


\begin{center}
  \resizebox{6cm}{!}{\includegraphics{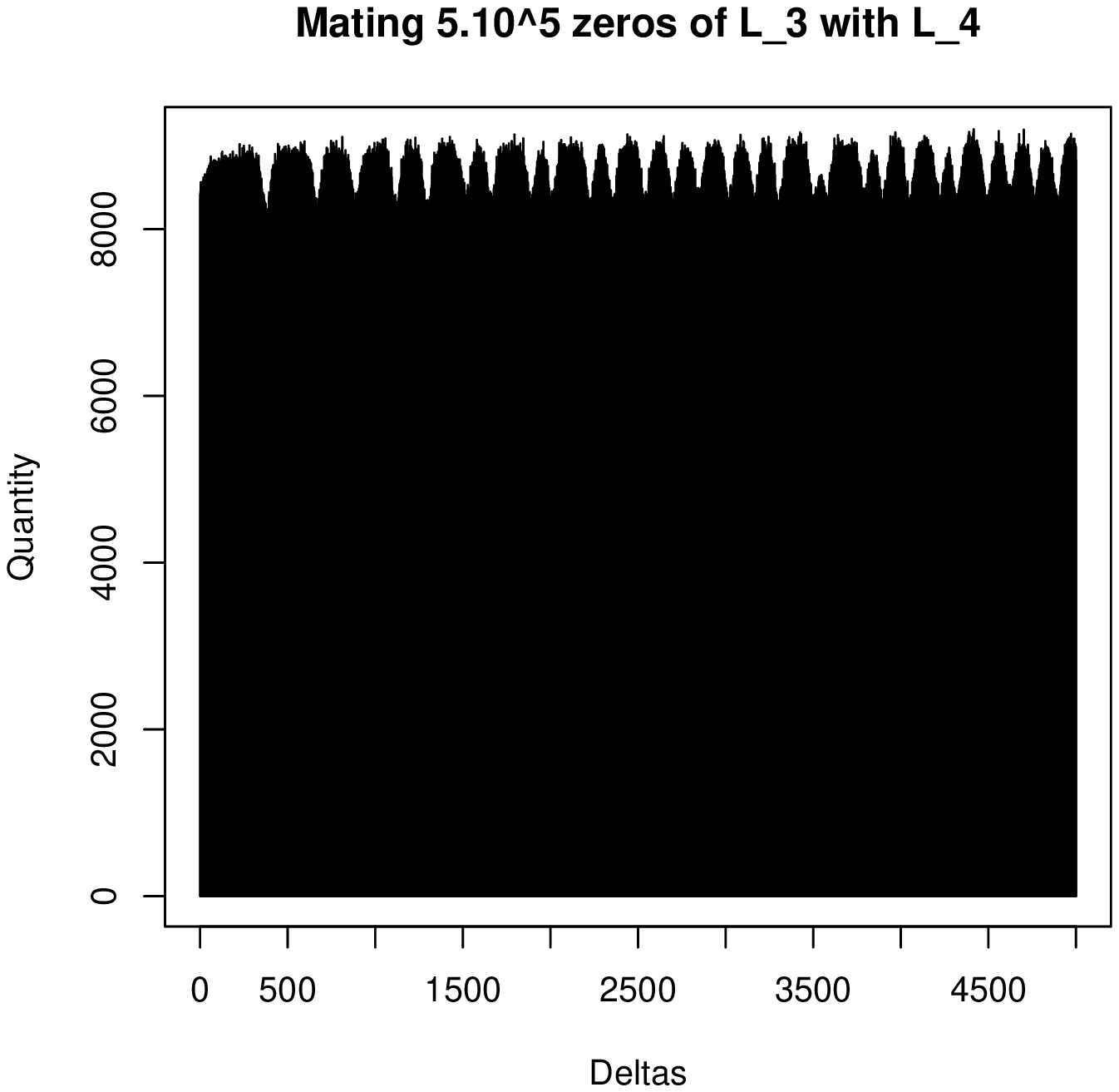}}
\end{center}

\centerline{Figure 31.}


\begin{center}
  \resizebox{6cm}{!}{\includegraphics{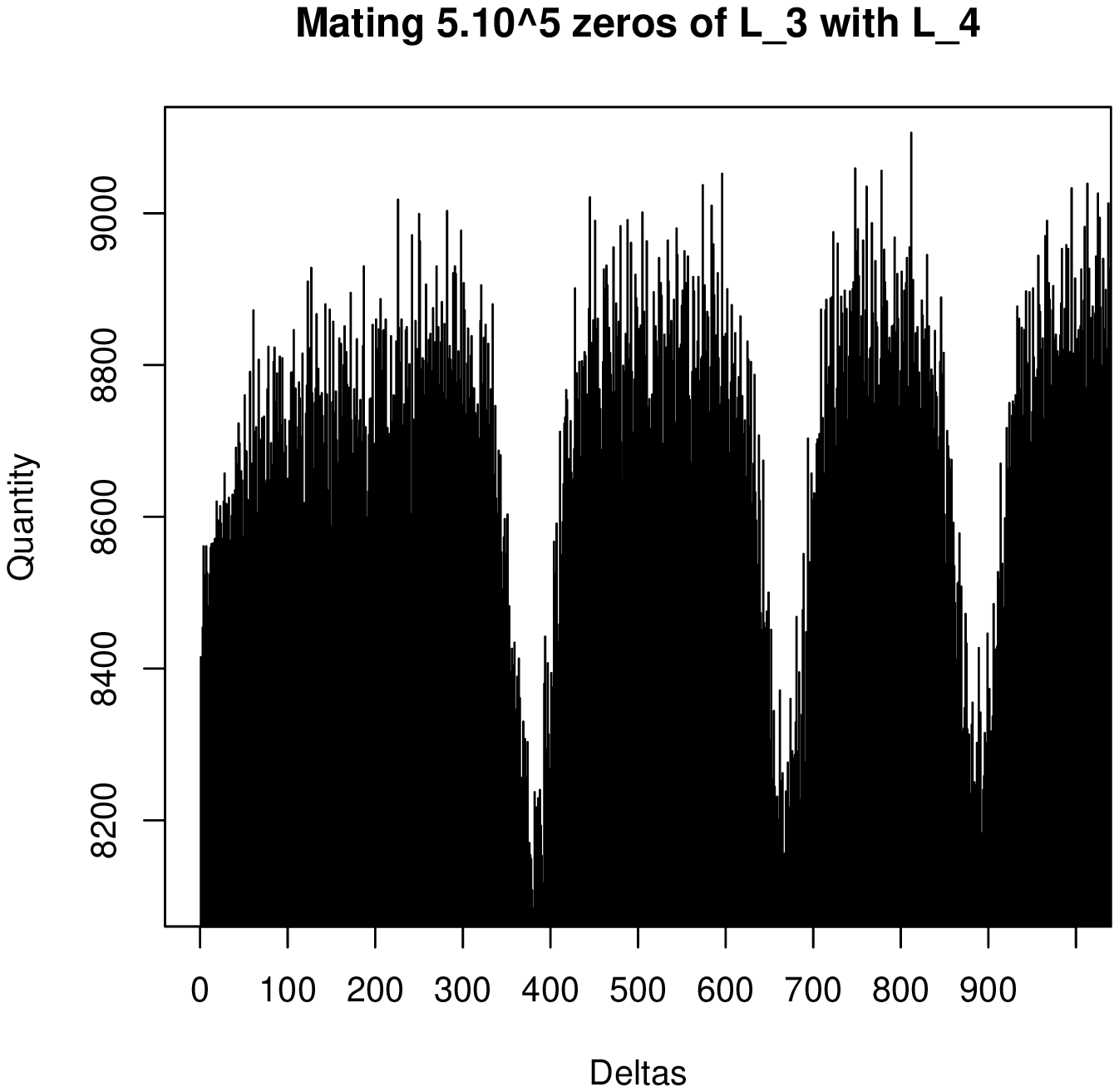}}
\end{center}

\centerline{Figure 32.}


\begin{center}
  \resizebox{6cm}{!}{\includegraphics{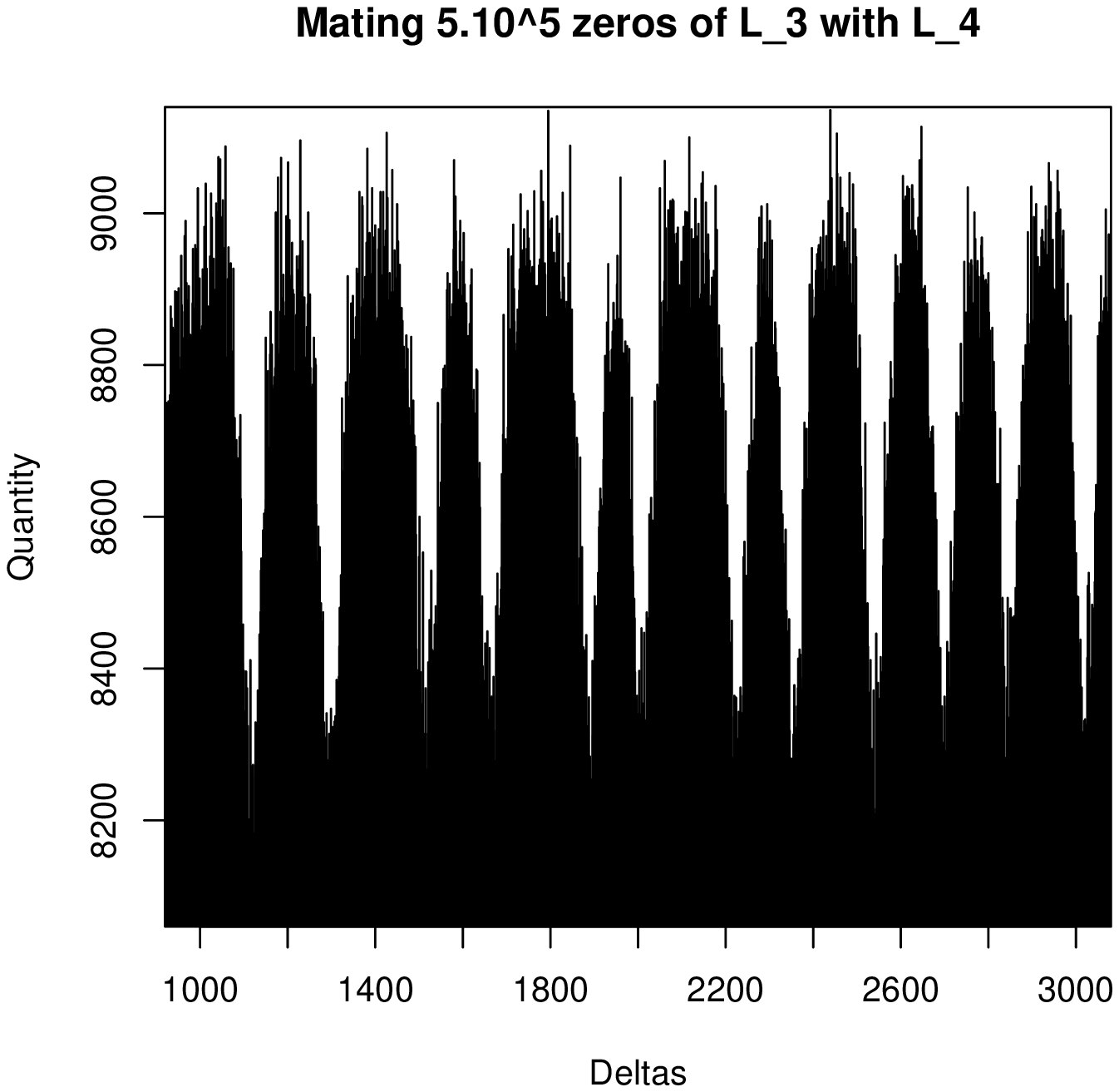}}
\end{center}

\centerline{Figure 33.}


\begin{center}
  \resizebox{6cm}{!}{\includegraphics{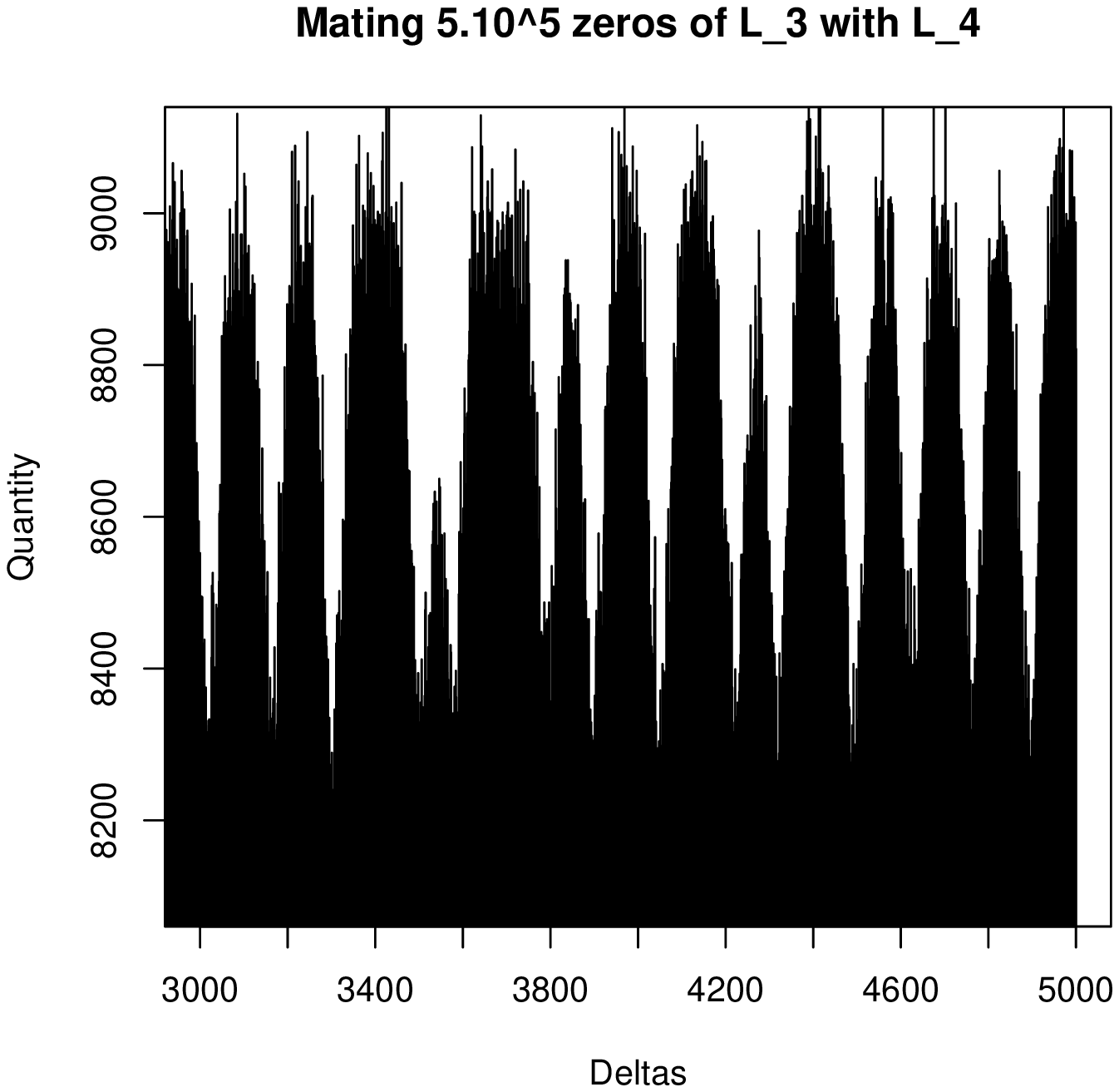}}
\end{center}

\centerline{Figure 34.}

\bigskip

 Again in this situation there is no GUE distribution near $0$
since the zeros of $L_{\chi_3}$ and $L_{\chi_4}$ are not
symmetric. Figure 35 shows the histogram of the deltas in the
range $[0,2]$ with precision $0.01$. This figure is to be compared
to figures 21.


\begin{center}
  \resizebox{6cm}{!}{\includegraphics{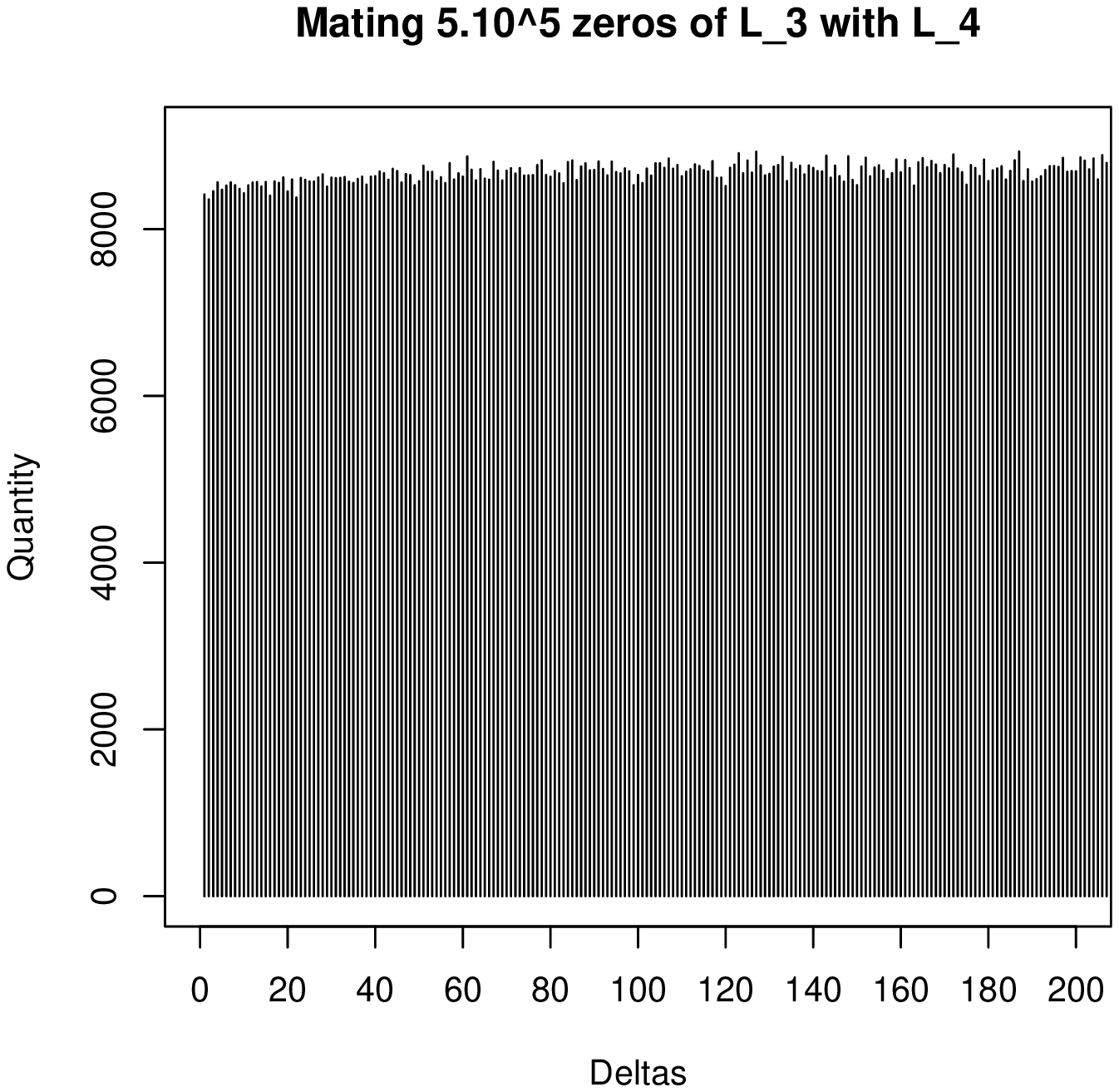}}
\end{center}

\centerline{Figure 35.}

\medskip

\textbf{E\~ne product explanation.}

\medskip

The computation follows. We have for any pair of Dirichlet
$L$-function $L_{\chi_1}$ and $L_{\chi_2}$,

\begin{align*}
L_{\chi_1} \bar \star \bar L_{\chi_2} &=L_{\chi_1} \bar
\star L_{\bar \chi_2} \\ &=\prod_p (1-\chi_1 (p)p^{-s})^{-1} \
\bar \star \ \prod_q (1-\bar \chi_2 (q) q^{-s})^{-1} \\ &=\prod_p
(1-\chi_1 (p)p^{-s})^{-1}\ \bar \star \ (1-\bar \chi_2(p)
p^{-s})^{-1} \\ &=\prod_p (1-(\chi_1 \bar \chi_2) (p)p^{-1/2}
p^{-s}) \\ &=\prod_p (1-(\chi_1 \bar \chi_2) (p) p^{-(s+1/2)})
\\ &=L_{\chi\bar \chi_2} (s+1/2)^{-1} 
\end{align*}

\section{Mating with local Euler factors.} \label{sec:Euler_factors}

In this section we study the mating of $L$-functions with 
local Euler factors. We mate the positive imaginary part of
Riemann zeros with the sequence of positive imaginary part of
zeros of Euler factor
$$
f_p(s)=(1-p^{-s}) \ ,
$$
which is the arithmetic sequence
$$
\gamma^{(p)}_k = {2\pi \over \log(p)} \ k ,
$$
with $k\in \ZZ$.

More generally we can consider the mating with general Euler
Dirichlet local factors
$$
f_{p,\chi}(s)=(1-\chi (p) p^{-s}) \ ,
$$
but we restrict the numerical statistics to Riemann Euler  local
factors.

 In order to have a substantial
number of deltas we need to work with a very large file of Riemann
zeros because this time the arithmetical sequence of Euler zeros has
constant density. For a given number of Riemann zeros we pick more
deltas if the prime $p$ is large. For this reason, with our
limited data of zeros available, we run the statistics for
$p=23$,
$$
{2\pi \over \log (23)}=2.00389\ldots \ ,
$$
and the first $10$ million Riemann zeros (statistics (a)); and
for $p=67$,
$$
{2\pi \over \log (67)}=1.494327\ldots \ ,
$$
and the first $15$ million of Riemann zeros (statistics (b)).

The histograms obtained (figures below for different ranges of
deltas) show a uniform distribution with deficit locations at the
corresponding Euler zeros, for $k\in \ZZ$,
$$
\gamma^{(p)}_k = {2\pi \over \log(p)} \ k .
$$
The deficit of deltas at these locations can also be seen for
statistics with smaller values of $p$, but because of the "small" number
of Riemann zeros at our disposition, the statistics is too poor
\footnote {Also our laptop has a RAM memory of 512 Mb
which does not allow R to read a vector with more than $20$
million Riemann zeros.}


\begin{center}
  \resizebox{6cm}{!}{\includegraphics{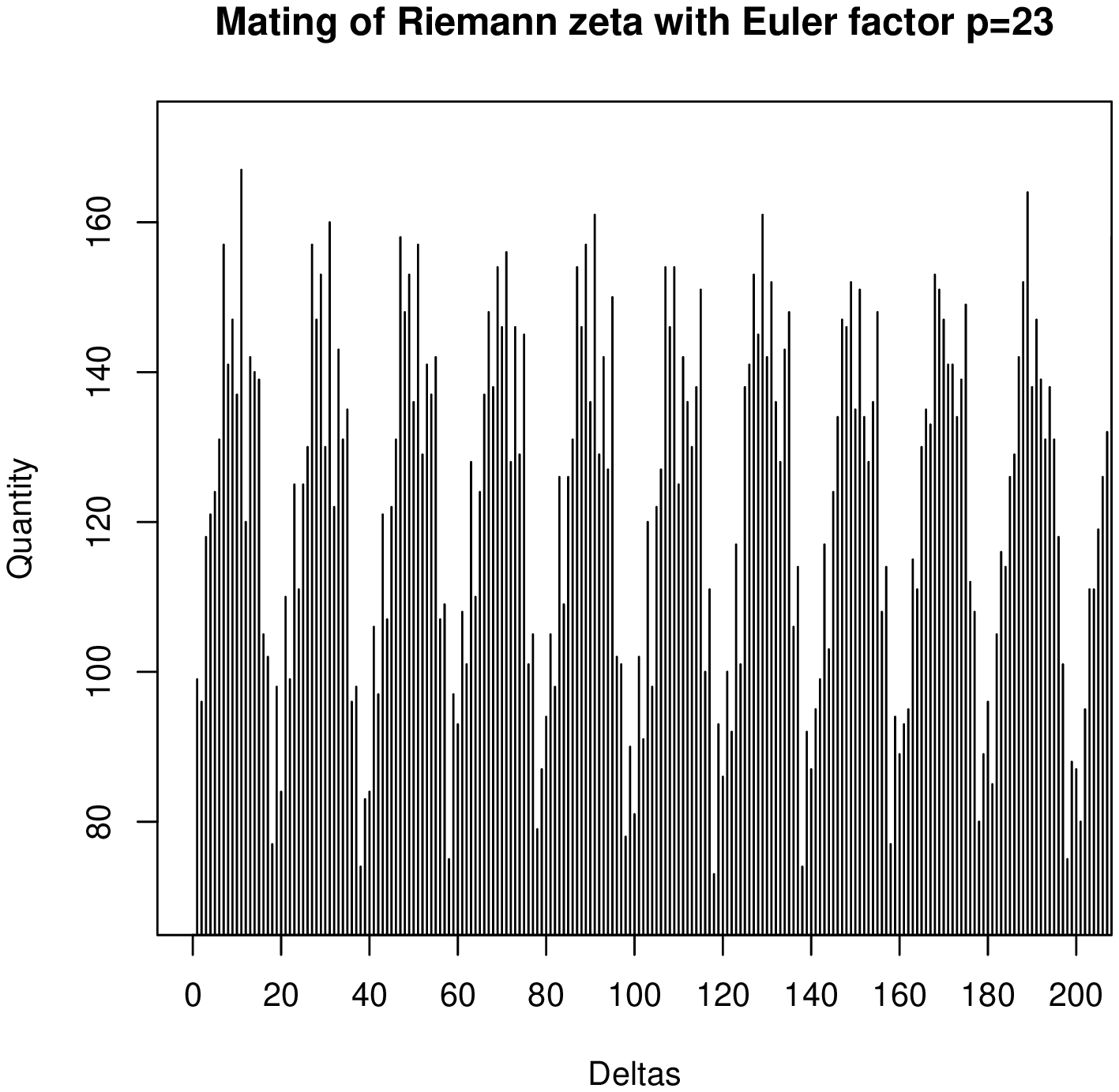}}    
  \resizebox{6cm}{!}{\includegraphics{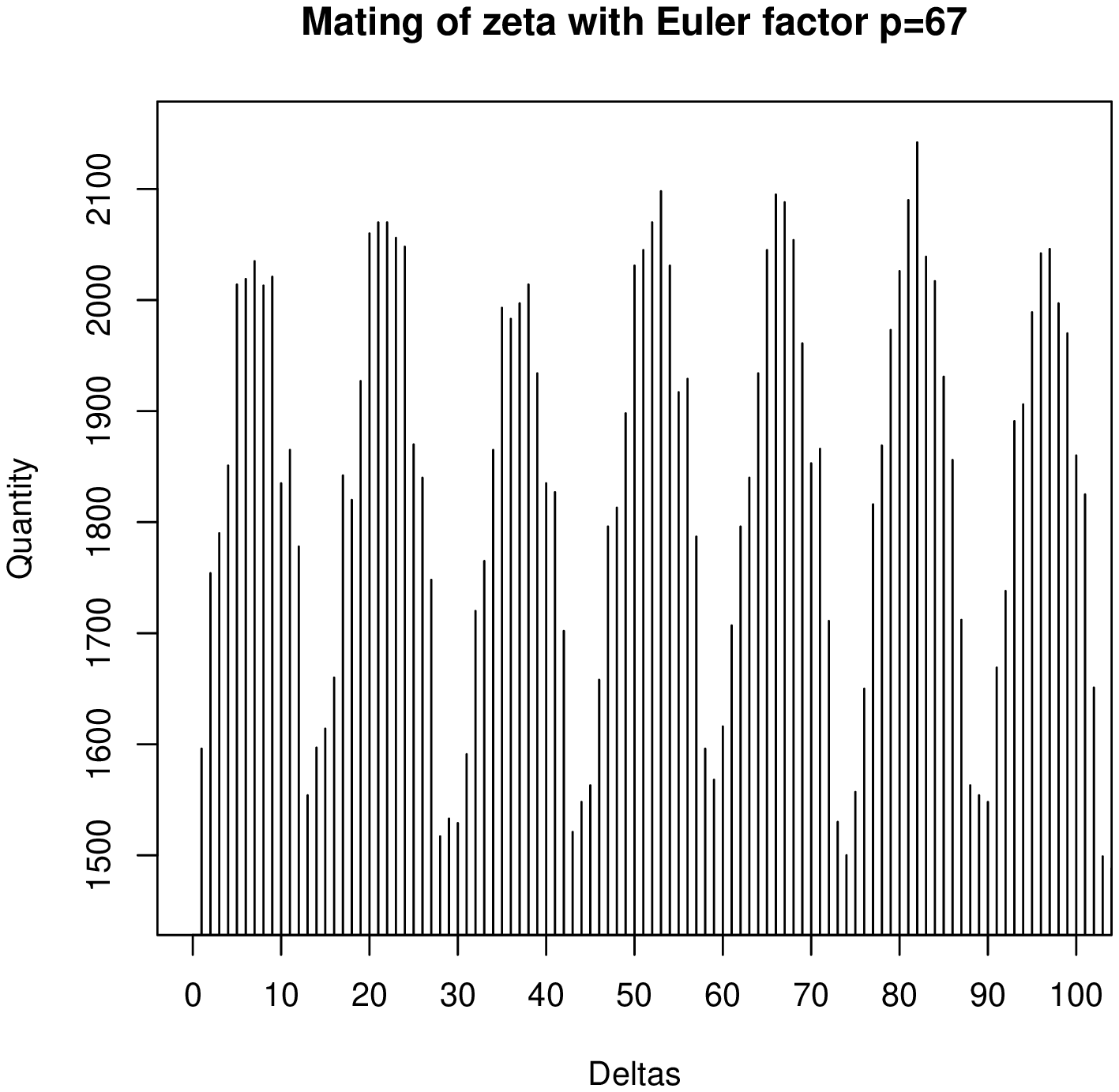}}
\end{center}

\centerline{Figure 36.a and 36.b.}


\begin{center}
  \resizebox{6cm}{!}{\includegraphics{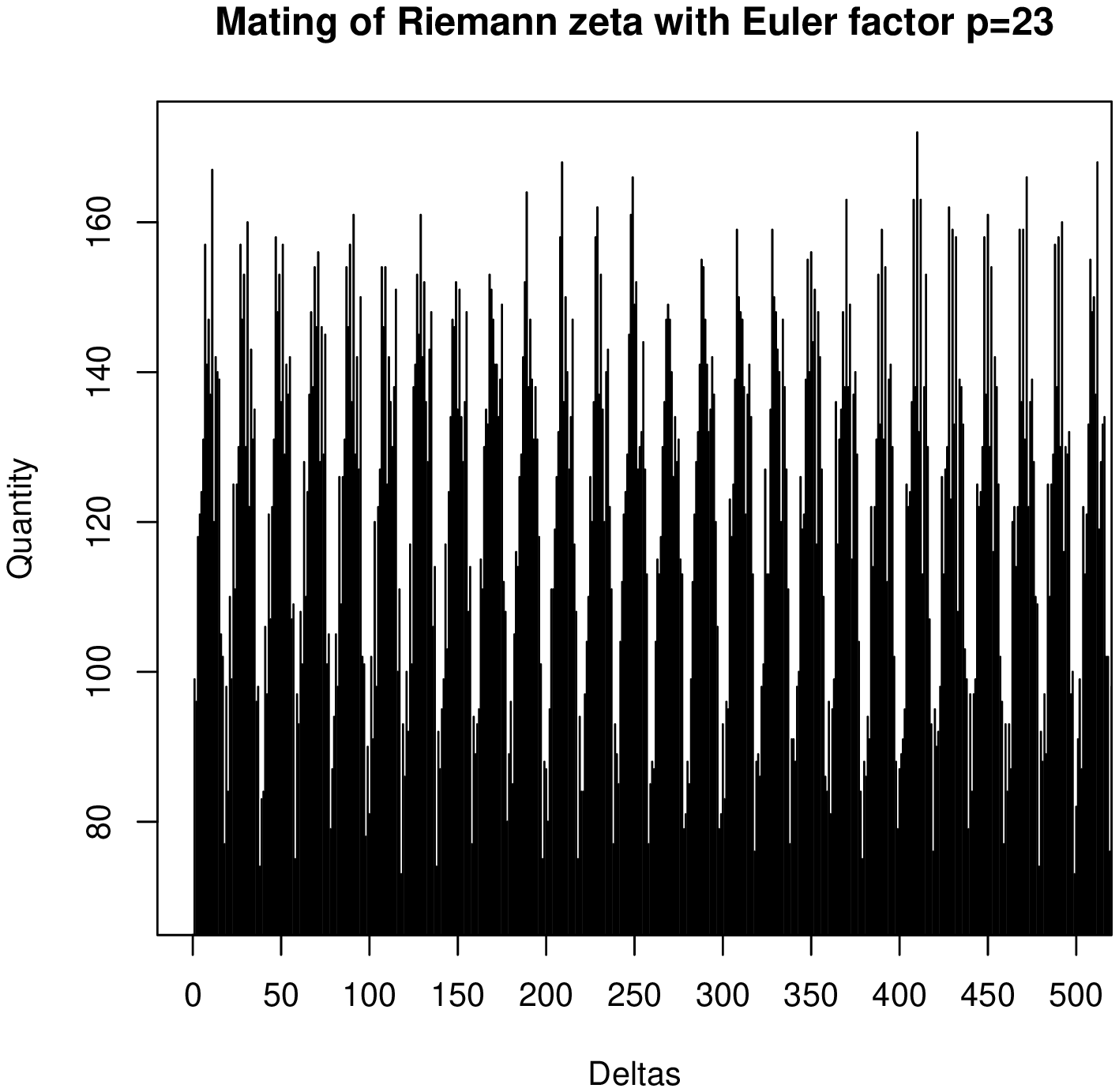}}    
  \resizebox{6cm}{!}{\includegraphics{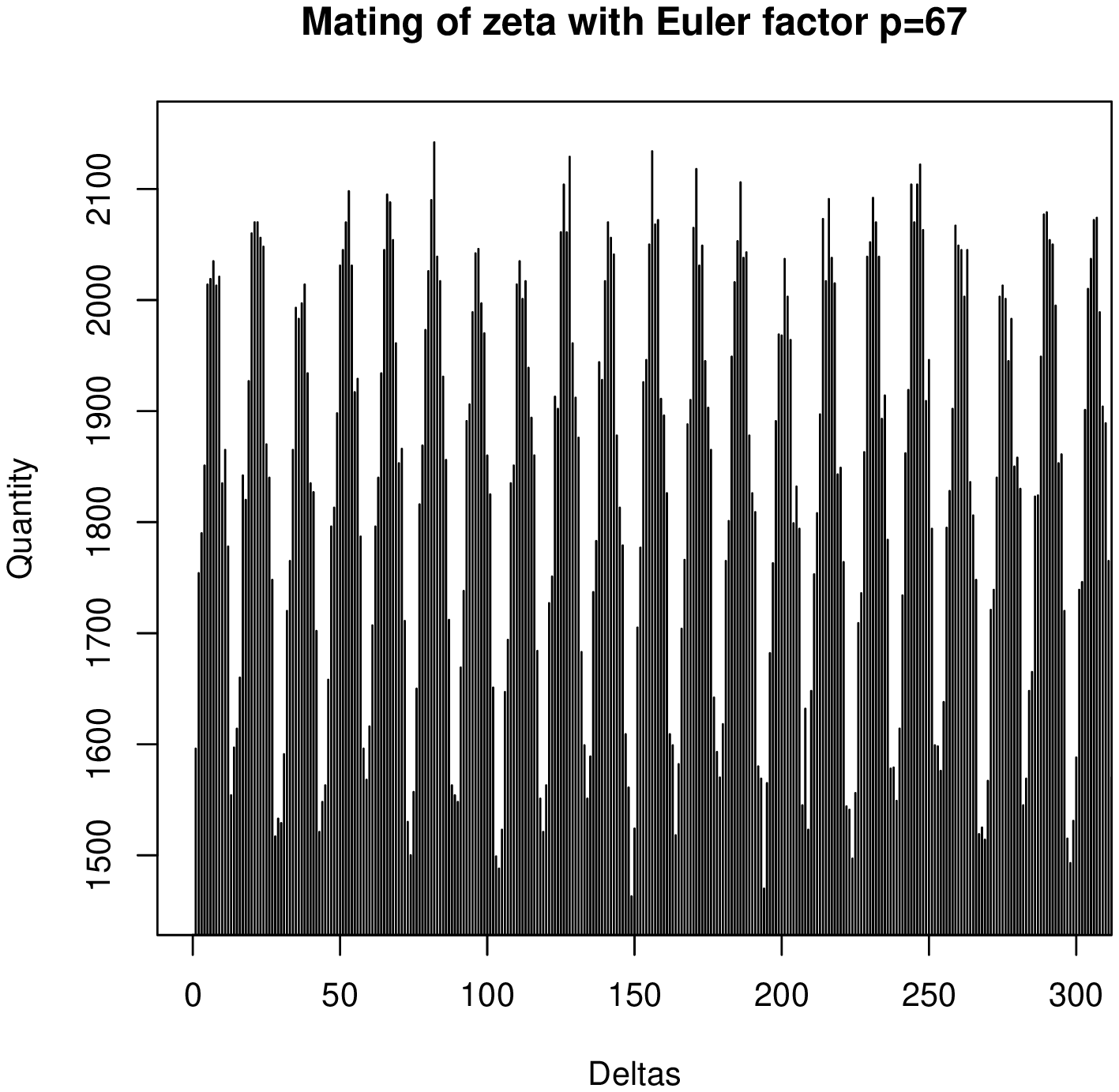}}
\end{center}

\centerline{Figure 37.a and 37.b.}


\begin{center}
  \resizebox{6cm}{!}{\includegraphics{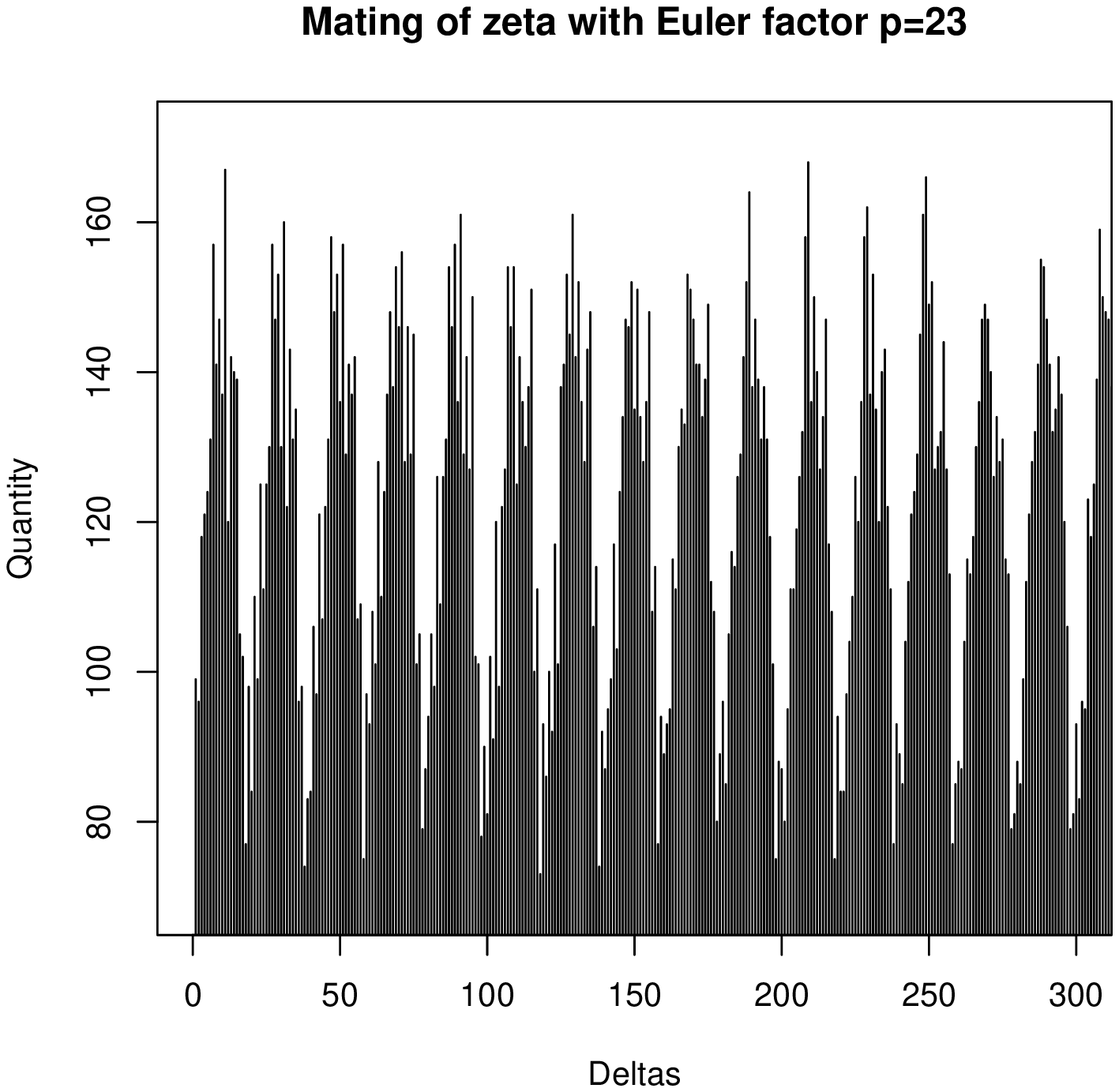}}    
  \resizebox{6cm}{!}{\includegraphics{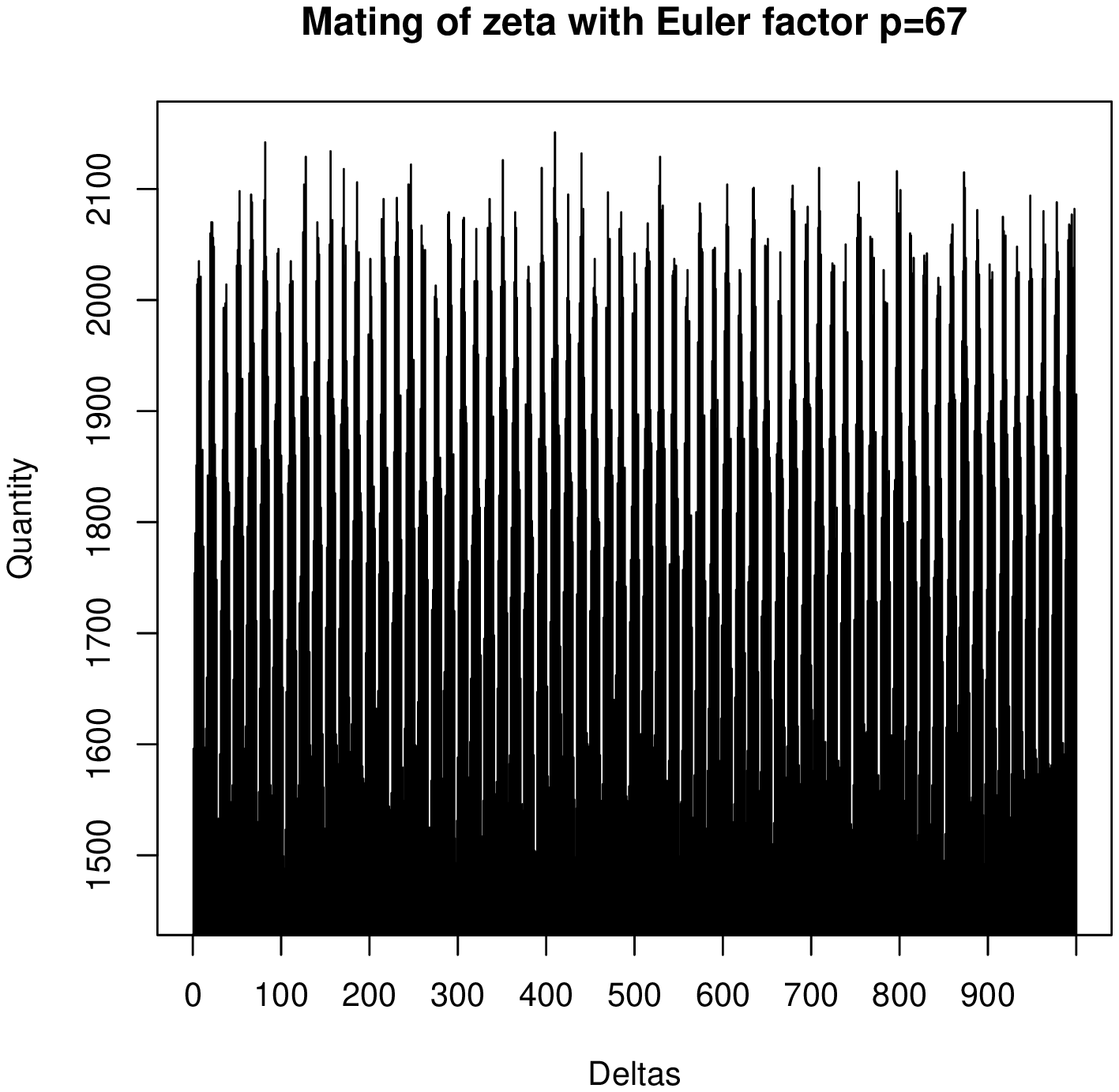}}
\end{center}

\centerline{Figure 38.a and 38.b.}

\bigskip

\textbf{E\~ne product explanation.}

\medskip

We have

\begin{align*}
\zeta \  \bar  \star \  \bar f_p &=\zeta \ \bar \star \
f_p \\ &=\prod_q (1-q^{-s})^{-1} \ \bar \star \ (1-p^{-s})^{-1}
\\ &= (1-p^{-s})^{-1} \ \bar \star \ (1-p^{-s})^{-1} \\
&=(1-p^{-s}) \ \bar \star \ (1-p^{-s})\\ &=(1-p^{-1/2} p^{-s})
\\ &=f_p(s+1/2) \ .
\end{align*}

We have in general for an arbitrary Dirichlet
$L$-function:
\begin{align*}
L_\chi \bar \star \bar f_p &=L_\chi \bar \star  f_p \\ 
&=\prod_q (1-\chi (q) q^{-s})^{-1} \ \bar \star \ (1-p^{-s})^{-1}
\\ &= (1-\chi (p) p^{-s})^{-1} \ \bar \star \ (1-p^{-s})^{-1} \\ 
&=(1-\chi (p) p^{-1/2}
p^{-s}) \\ &=f_{p,\chi}(s+1/2) \ .
\end{align*}

\section{Fine structure of deltas near $0$.} \label{sec:fine}

We come back in this section to analyze the fine structure of the statistics of
deltas near $0$.

As we have observed, Montgomery conjecture is verified
numerically for the deltas of zeros of arbitrary $L$-functions, but
not for the mating of zeros of non-conjugate $L$-functions. As 
explained, the GUE distribution arises because of the
symmetry of the zeros mated, i.e. it is a genuine real-analytic phenomenon.

For the Riemann zeta function, the e\~ne product analysis of the
distribution of deltas near zero reveals that after the first
order GUE correction, we have a second order term corresponding to
the pole of $\zeta (s+1/2)$ at $s=0$. This yields a positive
Fresnel distribution. 

We verify numerically this first order correction to the GUE
distribution. In the numerical application we consider the deltas
of $5$ million Riemann zeros. These have imaginary part less than
$T_0$ for $ T_0=2\ 630\ 122 $. We consider the frequency
$$
\omega_0= {1\over 2\pi}\log T_0 =0.02352714\ldots
$$
We correct the histogram of the deltas by adding a  GUE density
$$
t\mapsto  A \left ({\sin (\pi \omega_0 t) \over  \pi \omega_0
t}\right )^2 \ .
$$

Figure 39 shows the distribution of the deltas and figure 40 the
corrected distribution of the deltas showing the Fresnel
distribution (both in the range $[0,2]$). Figures 41 and 42 show
the same data in the range $[0,6]$ for a better view of the queue.
We have adjusted the coefficient $A$ in order to fit the best with
a Fresnel distribution looking at the first minima and the second
maximum. 

\begin{center}
  \resizebox{5.8cm}{!}{\includegraphics{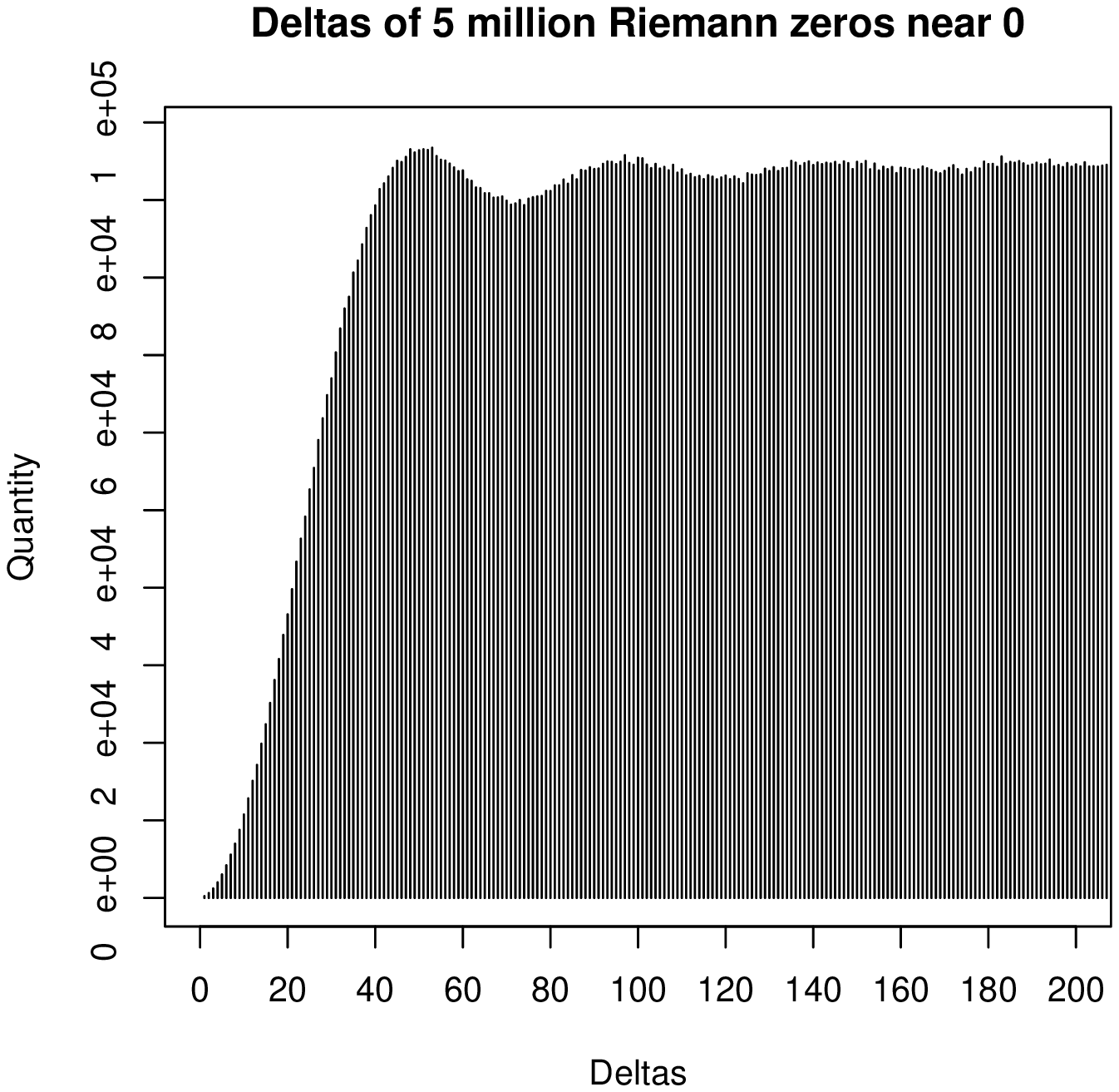}}
\end{center}

\centerline{Figure 39.}

\begin{center}
  \resizebox{6cm}{!}{\includegraphics{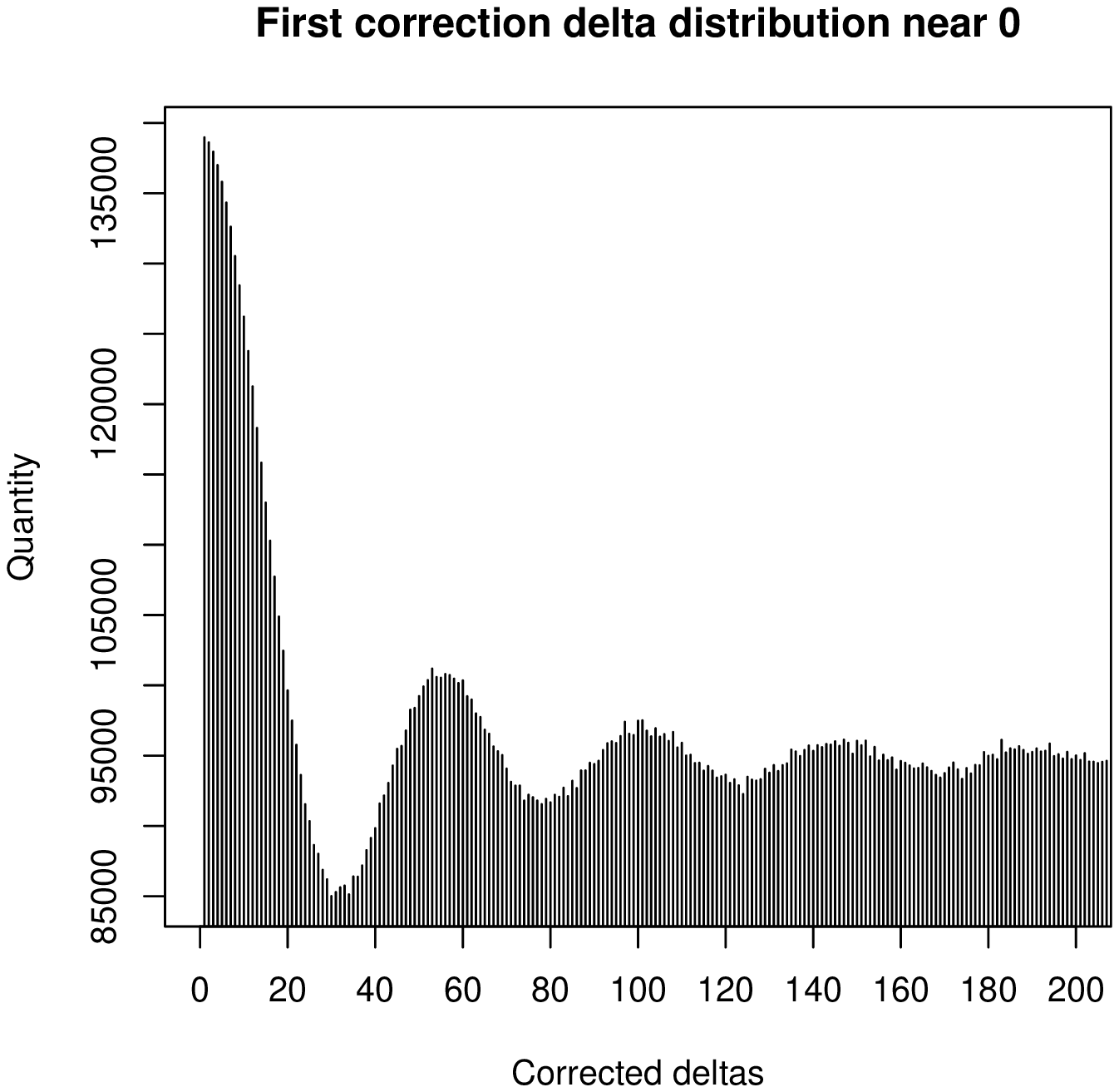}}
\end{center}

\centerline{Figure 40.}

\begin{center}
  \resizebox{6cm}{!}{\includegraphics{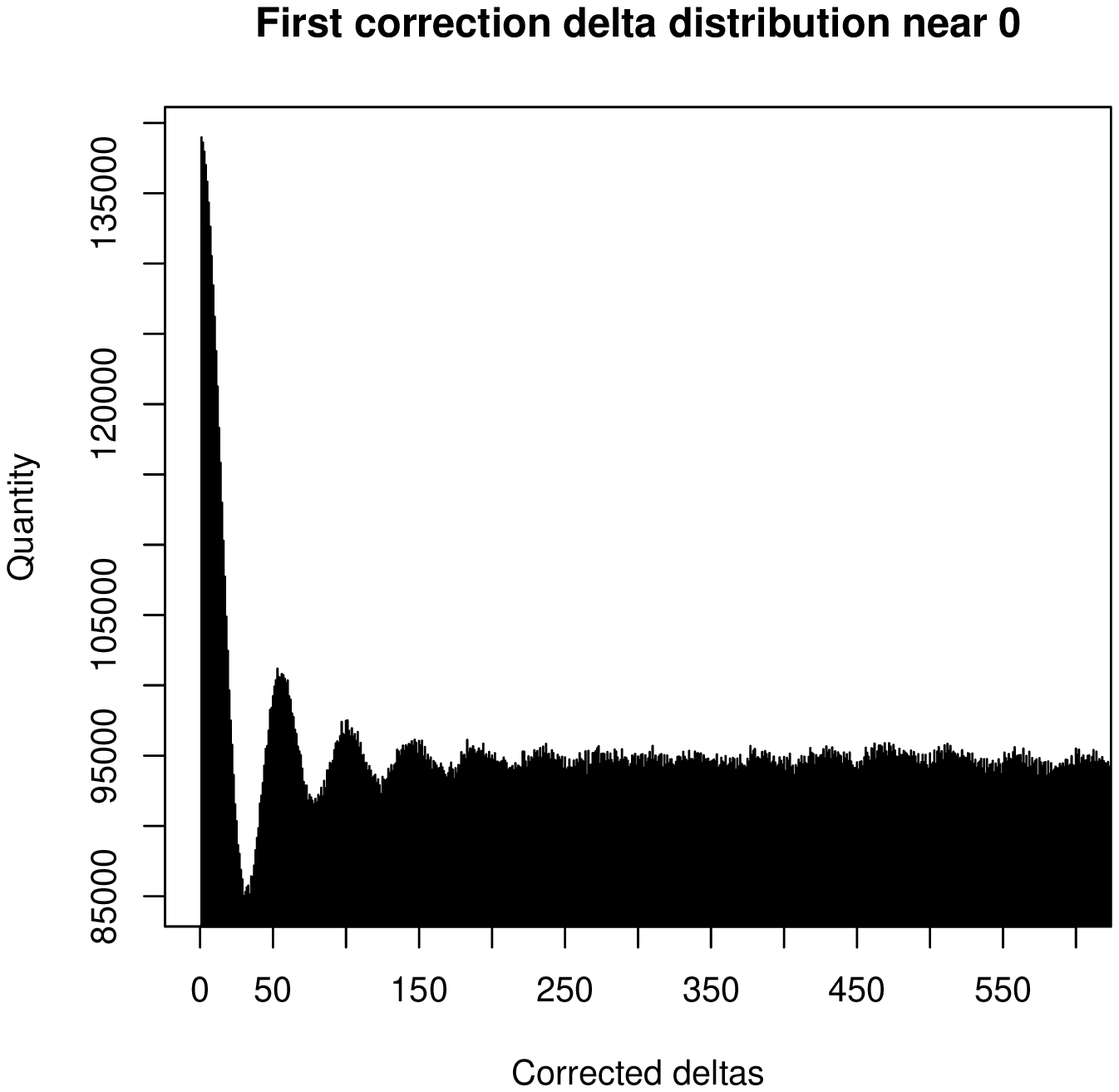}}
\end{center}

\centerline{Figure 41.}

\begin{center}
  \resizebox{6cm}{!}{\includegraphics{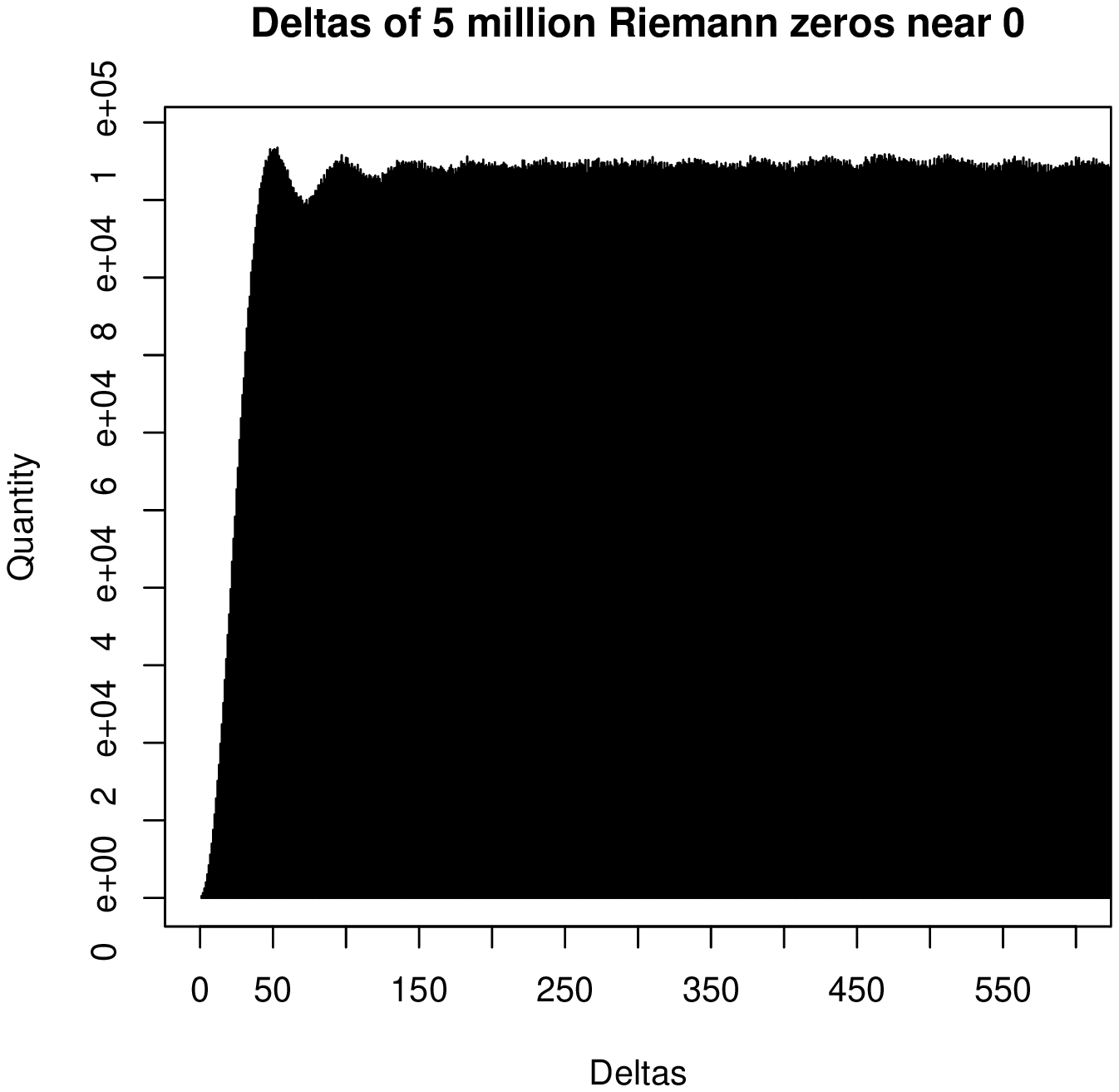}}
\end{center}

\centerline{Figure 42.}

%
%
%
%

\bigskip

\textbf{Script.}

\medskip

We have previously stored the result of computing the deltas of
$5$ million Riemann zeros in the range $[0,200]$ with precision
$0.01$. The variable $x$ contains the cumulative count of deltas
with the stated precision. Note also that $\pi .\omega_0/100=0.0739127\ldots$


{\tt

load("data 5 million")

correction=numeric()

for (k in 1:200000)
correction[k]=139000*(sin(0.0739127*k)/(0.0739127*k))\^ 2

x1<-x+correction

barplot(x1)

}


\begin{thebibliography}{1000}


\bibitem [Bo]{Bo} BOMBIERI, E., \emph{Problems of the Millenium: The Riemann
Hypothesis}, www.claymath.org, Official problem description.


\bibitem [Co]{Co} CONREY, J.B., \emph{The Riemann Hypothesis}, Notices of the
AMS, March 2003, p.341-353.




\bibitem [KS]{KS} KATZ, N.M.; SARNAK, P., \emph{Random matrices, Frobenius
eigenvalues, and monodromy}, AMS Colloquium Publications, {\bf 45},
AMS, Providence RI, 1999.



\bibitem[Mo1]{Mo1} MONTGOMERY, H.L., \emph{The pair correlations of zeros of the
zeta function}, {\it Analytic Number Theory}, editor H.G. Diamond,
Proc. Symp. Pure Math., Providence, 1973, p.181-193.


\bibitem[Mo2]{Mo2} MONTGOMERY, H.L., \emph{Distribution of the Zeros of the
Riemann Zeta function}, Proc. ICM, Vancouver, 1974, p.379-381.


\bibitem [Od]{Od} ODLYZKO, A.M. \emph {www.dtc.umn.edu/\~ \ odlyzko}, Personal
web page.

\bibitem[PM1]{PM1} P\'EREZ MARCO, R., \emph{The e\~ne product}, Manuscript.

\bibitem[PM2]{PM2} P\'EREZ MARCO, R., \emph{E\~ne product and Riemann zeta function
}, Manuscript.


\bibitem[Ri1]{[Ri1]} RIEMANN, B., \emph{Ueber die Anzahl der Primzahlen unter
einer gegebenen Gr\"osse}, Monat. der K\"onigl. Preuss. Akad. der
Wissen. zu Berlin aus der Jahre, 1859 (1860), p.671-680;
Gessammelte math. Werke und wissensch. Nachlass, 2 Aufl. 1892,
p.145-155.

\bibitem[Ri2]{[Ri2]} RIEMANN, B., \emph{Original manuscripts related to
[Ri1]}, Scan available at www.claymath.org.


\bibitem[Ru]{Ru} RUBINSTEIN, M., \emph{ pmmac03.math.uwaterloo.ca/\~\
mrubinst/L\_ function\_ public/ZEROS} , Public web page.



\end{thebibliography}
\end{document}